\documentstyle[11pt]{article}
\begin{document}
\small

\centerline {\bf MODIFIED MIXED TSIRELSON SPACES}

\bigskip

\medskip

\centerline {\sc S.A. Argyros -- I. Deliyanni}

\medskip

\centerline {\sc D.N. Kutzarova -- A. Manoussakis}

\bigskip

\bigskip

\begin{abstract}
\footnotesize

We study the modified and boundedly modified mixed Tsirelson
spaces
$T_M[({\cal F}_{k_n},\theta_n)_{n=1}^{\infty }]$
and $T_{M(s)}[({\cal F}_{k_n},\theta_n)_{n=1}^{\infty }]$
respectively, defined by a subsequence $({\cal F}_{k_n})$
of the sequence of Schreier families $({\cal F}_n)$. These
are reflexive asymptotic $\ell_1$ spaces with an
unconditional basis $(e_i)_i$ having the property that
every sequence $\{ x_i\}_{i=1}^n$ of normalized disjointly
supported vectors contained in $\langle e_i\rangle_{i=n}^{\infty }$
is equivalent to the basis of $\ell_1^n$. We show that
if $\lim\theta_n^{1/n}=1$ then the space
$T[({\cal F}_n,\theta_n)_{n=1}^{\infty }]$ and its modified
variations $T_M[({\cal F}_n,\theta_n)_{n=1}^{\infty }]$ or
$T_{M(s)}[({\cal F}_n,\theta_n)_{n=1}^{\infty }]$ are
totally incomparable by proving that $c_0$ is finitely
disjointly representable in every block subspace of
$T[({\cal F}_n,\theta_n)_{n=1}^{\infty }]$. Next, we present
an example of a boundedly modified mixed Tsirelson
space $X_{M(1),u}=T_{M(1)}[({\cal F}_{k_n},\theta_n)_{n=1}^{\infty }]$
which is arbitrarily distortable. Finally, we construct
a variation of the space $X_{M(1),u}$ which is hereditarily
indecomposable.

\end{abstract}

\bigskip

\bigskip

\medskip

\centerline { {\sc  Introduction}}

\medskip

Given a sequence $({\cal M}_k)_{k=1}^{\infty}$ of compact families of
finite subsets of ${\bf N}$ and a sequence $(\theta_k)_{k=1}^{\infty}$
of reals converging to zero, the mixed Tsirelson space 
$T[({\cal M}_k,\theta_k)_{k=1}^{\infty}]$
is defined as follows.

$T[({\cal M}_k,\theta_k)_{k=1}^{\infty}]$ is the completion of the linear
space $c_{00}$ of the sequences which are eventually zero under the norm
$\| \cdot\|$ defined by the following implicit formula:
For $x\in c_{00},$ 
$$\|x\|=\max\big\{\|x\|_{\infty},\sup_k\theta_k\sup
\{\sum_{i=1}^n\| E_ix\|: n\in {\bf N}, (E_i)_{i=1}^n \; {\rm is}\;
{\cal M}_k-{\rm admissible}\}\big\}.\leqno (1)$$
Here, for $E\subset N$, $\| Ex\|$ is the restriction of the vector $x$ on 
the set $E$ and, for a family ${\cal M}$ of subsets of ${\bf N}$, an
${\cal M}$--{\it admissible } sequence is a sequence $(E_i)_{i=1}^n$
of {\it successive} subsets of ${\bf N}$ such that the set
$\{ \min E_1,\ldots ,\min E_n\}$ belongs to ${\cal M}.$ Mixed Tsirelson
spaces were introduced in [3]. However, this class includes the
previously constructed Schlumprecht's space ([16]) which initiated
a series of results answering fundamental and long standing problems of the 
theory of Banach spaces. The remarkable nonlinear transfer by Odell
and Schlumprecht ([13]) of the biorthogonal asymptotic sets from
Schlumprecht's space to $\ell_p$, $1<p<\infty$, which settled the 
distortion problem, indicates the impact of the new spaces on the
understanding of the classical Banach spaces. On the other hand,
these new norms led to the discovery of the class of hereditarily
indecomposable (H.I.) spaces ([9]), that is, spaces with the
property that no subspace can be written as a topological direct sum
of two infinite dimensional closed  subspaces. As it
was proved by Gowers ([8]), the H.I. property is a consequence of 
the absence of unconditionality in the sense that every Banach space
which does not contain any unconditional basic sequence has an H.I.
subspace. Gowers and Maurey ([9]) have proved that the H.I. spaces
have small spaces of operators; it is a fundamental open problem 
whether  there exists such a space
with the property that every bounded linear operator
$T:X\rightarrow X$  is of the form $T=\lambda I+K$ where $K$ is a
compact operator. On the other hand, a recent result of Argyros
and Felouzis ([4]) shows that a large class of Banach spaces that 
includes $\ell_p$, $1<p<\infty$, are quotients of H.I. spaces.

\medskip

 In the present paper we study variations
of   mixed Tsirelson
spaces which we call modified mixed Tsirelson spaces.
Given a family ${\cal M}$ of finite subsets of ${\bf N}$, a sequence $(E_i)
_{i=1}^n$ of subsets of ${\bf N}$ is called ${\cal M}$--{\it allowable}
if the sets $E_i$ are {\it disjoint} and the set 
$\{\min E_1,\ldots, \min E_n\}$ belongs to ${\cal M}.$ The {\it modified mixed
Tsirelson space} $X_M$ corresponding to the mixed Tsirelson space 
$X=T[({\cal M}_k,\theta_k)_{k=1}^{\infty}]$ 
is the Banach space whose norm  $\| \cdot \|$ satisfies the implicit equation
$$\|x\|=\max\big\{\|x\|_{\infty},\sup_k\theta_k\sup
\{\sum_{i=1}^n\| E_ix\|: n\in {\bf N}, (E_i)_{i=1}^n \; {\rm is}\;
{\cal M}_k-{\rm allowable}\}\big\}.\leqno (2)$$

We also consider {\it boundedly modified mixed Tsirelson spaces} that lie
between $X$ and $X_M$. Such a space is denoted by $X_{M(s)},$ for
some $s\in {\bf N},$  and its norm is given by an implicit formula
analogous to (1) or (2) where the inner ``sup'' is taken over all
${\cal M}_k$--allowable families for $1\leq k\leq s$ and over all
${\cal M}_k$--admissible families for $k\geq s+1.$ It is clear that 
the modified and boundedly modified mixed Tsirelson spaces which are defined 
by a subsequence ${\cal M}_k={\cal F}_{n_k}$ 
of the sequence of Schreier families $({\cal F}_n)_n$ have the property that,
for every $n$, every normalized sequence $(x_i)_{i=1}^n$ of $n$ disjointly
supported vectors with supports contained in $[n,\infty)$ is
$\theta_1$-equivalent to the basis of $\ell_1^n.$

The modified Tsirelson space $T_M$ was introduced by W.B. Johnson ([10])
shortly after Tsirelson's discovery ([19]). Later, P. Casazza
and E. Odell ([6]) proved that the modified Tsirelson space
is isomorphic to the original one. The use of the modified
version of the norm in the 2--convexification of
$T$ is crucial for the proof of the fact that it is a
weak Hilbert space.
The relation between modified mixed Tsirelson norms and
the corresponding mixed Tsirelson norms is 
in general quite different from the one between $T$ and $T_M$.
 To explain the
situation we restrict our attention to the two main examples
of mixed Tsirelson norms:

The first is Schlumprecht's space $S$ ([16]) defined by 
${\cal M}_k={\cal A}_k=\{ A\subset {\bf N}: \# A\leq k\},$ and 
$\theta_k=
\frac{1}{\log_2(k+1)}$. The  
second is the space $X$ introduced by Argyros and Deliyanni
in [3], defined by a certain
 subsequence 
$({\cal F}_{n_k})_{k\in {\bf N}}$ of the  sequence of Schreier families
$({\cal F}_n)_{n\in N}$ and an appropriate sequence $(\theta_k)_{k\in{\bf N}}.$
 It is known that $c_0$ is finitely representable in every
infinite dimensional subspace of $S$ and we
show here that the same holds true for $X$. From this
we easily see that the modified versions $S_M$,
$X_M$ are totally incomparable to $S$ and $X$
respectively. Schlumprecht observed further that although
his space $S$ is reflexive, the space $S_M$
contains $\ell_1$ ([17]). On the other hand, as we show here,
the space $X_M$ remains reflexive and contains no $\ell_p$.
This is the first property where we do not have an analogy
between $S$ and $X$. The result is somehow
unexpected since $X_M$, being an asymptotic $\ell_1$ space,
has  richer $\ell_1$ structure than $S_M$.
 These results raise naturally certain questions related
to the structure of $S_M$ and $X_M$. For example,
it is not known if $S_M$ is $\ell_1$--saturated
or if $X_M$ is arbitrarily distortable.

\medskip

The results mentioned above are presented in  Section 1. More precisely,
we prove that if $\lim \theta^{1/n}_n=1$, then $c_0$ is finitely
representable in every infinite dimensional subspace of the space 
$T[({\cal F}_n,\theta_n)_{n=1}^{\infty}].$ 
Next, for an arbitrary null sequence $(\theta_n)_n$, we show 
that the modified mixed
Tsirelson space 
$T_M[({\cal F}_n,\theta_n)_{n=1}^{\infty}]$ is reflexive. As a consequence
we get that the 2-convexifications of such spaces yield weak Hilbert
spaces not containing $\ell_2$ and totally incomparable to $T^{(2)}.$

 In  Section 2 we  consider a boundedly modified mixed Tsirelson
space of the form $X_{M(1),u}=
T_{M(1)}[({\cal F}_{k_n},\theta_n)_{n=1}^{\infty}]$ for a suitable choice
of $({\cal F}_{k_n})$ and $(\theta_n)$. We show that this space is arbitrarily
distortable. This result is related to the question: Does there exist
a distortable Banach space of bounded distortion? By
 [12], [11] and  [18], such a
space must contain an asymptotic $\ell_p$
subspace with an unconditional basis which contains 
$\ell_1^n$'s uniformly; so the search turns to 
asymptotic $\ell_1 $ spaces with an unconditional basis. By [3] (also
[2]), the class of spaces $T[({\cal F}_n,\theta_n)_n]$ provides  examples
of such spaces
 which are arbitrarily distortable. However, it is not known 
whether the original representative of
this class, Tsirelson's space $T$, is arbitrarily 
distortable, or whether it contains an arbitrarily distortable subspace.
The space $X_{M(1),u}$ constructed here is closer to $T$ than  
 $T[({\cal F}_n,\theta_n)_n]$,  in the sense that it has  more homogeneous
$\ell_1$ structure.

 In Section 3 we construct a space $X$ based on 
$X_{M(1),u}$ which is hereditarily indecomposable. The basic idea for
the definition of $X$ comes from [9].

The strategy in proving these results is similar to  the one
followed in [3]. We briefly explain the idea. In order
to prove that $X_{M(1),u}$ is arbitrarily distortable,
we start with a set $K=\cup_{j=1}^{\infty }{\cal A}_j$
of functionals which define the norm of the space. Each set
${\cal A}_j$ contains functionals of the form
$\theta_j\sum_{l=1}^nf_l$ where $\{ f_l\}_{l=1}^n$ are
disjointly supported 
functionals in the dual ball and the family $\{{\rm supp}f_l\}_{l=1}^n$ is
${\cal F}_{k_j}$-- allowable if $j=1$ or ${\cal F}_{k_j}
$--admissible if $j>1$. Our goal is to show the following:

There exists $c>0$ such that for every block subspace $Y$ of $X_{M(1),u}$
and for large $j$ there exists $y_j\in Y$ with $\|y_j\|=1$
satisfying
$$\|y_j\|\approx\sup\{ f(y_j):f\in {\cal A}_j\},\leqno (3)$$
$$|f(y_j)|\leq c\theta_i\;{\rm for}\;{\rm all}
\;i< j,f\in {\cal A}_i.\leqno (4)$$

\noindent These two conditions imply that $X_{M(1),u}
$ is an arbitrarily distortable
space.

The fundamental objects that we use in order to find such vectors
$y_j$ are the $(\varepsilon,j )$--{\it basic special convex combinations}.
The $(\varepsilon,j )$--basic s.c.c. are convex combinations of the
basis $(e_n)_{n\in {\bf N}}$ of the space $X_{M(1),u}$ 
whose normalizations satisfy conditions (3) and (4) if
$\varepsilon $ is small enough.
The choice of $(\theta_n)_n, ({\cal F}_{k_n})_n$ ensures that
for every $j\geq 2$  and for every infinite
$D\subseteq {\bf N}$, there exists an  $(\varepsilon ,j)$--basic
special convex combination supported in $D$. 

Next  we show that in every block subspace $Y$ of 
$X_{M(1),u}$ and for every $j\geq 2$ we can choose 
a normalized vector $y_j$ in $Y$ with the following property:
For every $i$ and every $f\in {\cal A}_i$, there exist an
$(\varepsilon, j)$--basic special convex combination $x_f$
and a functional $g_f\in{\cal A}_i$ such that 
$$|f(y_j)|\leq Cg_f(x_f)$$ for some constant $C$.
\noindent Thus, we reduce the estimation  of the action
of $
{\cal A}_i$ on $y_j$ to the estimation of the action of
${\cal A}_i$ on  basic special
convex combinations.
 Our basic tool for this proof is the {\it analysis} of
a functional $f\in\cup_{i=1}^{\infty}{\cal A}_i$ which is the array
of functionals used for the inductive construction of $f$.

In the case of the space $X$ with no unconditional basic sequence
which is constructed in the third section, the scheme of ideas
is similar with some additional difficulties coming from the
existence of the dependent chains of functionals.

\bigskip

\bigskip

\centerline {{\sc  1. Mixed Tsirelson Spaces and their modified versions.}}

\medskip

\noindent {\bf A. Preliminaries.}

\medskip

\noindent {\bf Notation.} Let 
$(e_i)_{i=1}^{\infty}$ be the standard 
basis of the linear space $c_{00}$ of 
finitely supported sequences. For $x=\sum_{i=1}^{\infty}a_ie_i \in c_{00}$,
the {\it support} of $x$ is the set ${\rm supp}x=\{ i\in {\bf N}:
a_i\neq 0\}.$ For $E, F$ finite subsets of ${\bf N}$, $E<F$ means
$\max E < \min F$ or either $E$ or $F$ is empty. For $n\in {\bf N},$
$E\subset {\bf N},$ $n<E$ (resp. $E<n$) means $n<\min E$ (resp. $\max E <n$).
For $x, y$ in $c_{00},$  $x<y$ means ${\rm supp}x<{\rm supp}y.$
For $n\in {\bf N}$, $x\in c_{00}$
 we write $n<x$ (resp. $x<n$) if $n<{\rm supp}x$
(resp. ${\rm supp}x <n$). We say that the sets
 $E_i\subset {\bf N}$, $ i=1,\ldots,n$
are {\it successive} if $E_1<E_2\ldots <E_n.$ Similarly, the vectors $x_i,
\; i=1,\ldots n$ are {\it successive} if $x_1<x_2<\ldots <
x_n.$ For $x=\sum_{i=1}
^{\infty}
a_ie_i$ and $E$ a  subset of ${\bf N}$, we  denote by $Ex$ the vector
$Ex=\sum _{i\in E}a_ie_i.$

\bigskip

\noindent {\bf The Schreier families ${\cal F}_{\alpha}.$}
Let ${\cal M}$ be a family of finite subsets of ${\bf N}.$ We say that 
${\cal M}$ is {\it compact} if it is closed in the topology of pointwise
convergence in $2^{\bf N}.$ ${\cal M}$ is {\it heriditary} if whenever
$B\subset A$ and $A\in {\cal M}$ then $B\in {\cal M}.$ ${\cal M}$ is
{\it spreading} if whenever $A=\{ m_1,\ldots ,m_k\}\in {\cal M}$ and
$B=\{n_1, \ldots, n_k\}$ is such that $m_i\leq n_i$, $i=1,\ldots k,$
then $B\in {\cal M}.$

\medskip

\noindent {\bf Notation.} Let ${\cal M}$, ${\cal N}$ be families of
finite subsets of ${\bf N}.$ We denote by ${\cal M}[{\cal N}]$ the
family

$${\cal M}[{\cal N}]=\left\{ \cup_{i=1}^n A_i : n\in {\bf N}, A_i\in {\cal N},
A_1<A_2<\ldots <A_n \;{\rm and}\; \{\min A_1,\ldots ,\min A_n\}\in {\cal M}
\right\}.$$

\medskip

 The Schreier family ${\cal S}$  is defined as follows:

$${\cal S}=\{ A\subset {\bf N}: \# A \leq \min A\}.$$

The {\it generalized Schreier families}
 ${\cal F}_{\alpha},$ $\alpha <\omega_1,$ were 
introduced in [1]:

\medskip

\noindent {\bf 1.1 Definition.} 

$${\cal F}_0=\{ \emptyset\}\cup \{ \{ n\} : n\in {\bf N}\}$$

$$ {\cal F}_{\alpha +1}=\{ \emptyset \} \cup \{ \cup_{i=1}^nA_i :
n\in {\bf N},\;  A_i\in {\cal F}_{\alpha}, \; n\leq A_1<A_2<\ldots <A_n\}$$

\noindent {\it and for a limit ordinal $\alpha$ we choose a sequence 
$(\alpha_n)_n$, $\alpha_n\uparrow \alpha$ and set}

$${\cal F}_{\alpha}=\{\emptyset\}\cup\{ A: {\rm there \; exists \; }n\in {\bf N}
\; {\rm such \; that}\; A\in {\cal F}_{\alpha_n} \; {\rm and}\; n\leq A\}.$$

\medskip

Notice that ${\cal F}_1={\cal S}.$ Also, for $n,m<\omega$, 
${\cal F}_n[{\cal F}_m]={\cal F}_{n+m}.$

 It is  easy to see that
each ${\cal F}_{\alpha}$ is a compact, hereditary and spreading family.

\medskip

\noindent {\bf 1.2 Lemma.} {\it For $n<\omega $ define the family
${\cal F}_n^M$ inductively as follows:

${\cal F}_0^M={\cal F}_0$.

${\cal F}_{n+1}^M=\{\cup_{i=1}^kA_i: k\in {\bf N}, A_i\in {\cal F}_n^M\;
{\rm for}\,i=1,\ldots ,k,\;A_i\cap A_j=\emptyset\;{\rm for}
\;i\neq j\;{\rm and}$

\hspace*{.5in} $k\leq\min A_1<\min A_2<\ldots <\min A_k\}$.

\noindent Then, for all $n$, ${\cal F}_n^M={\cal F}_n$.}

\medskip

\noindent {\it Proof:} The proof is an immediate consequence of the
following.

\noindent {\bf Claim:} {\it Let $n\in {\bf N}$ and let $A_i\in {\cal F}_n$,
$i=1,\ldots ,k$ be such that $A_i\cap A_j=\emptyset $ for $i\neq j$
and $\min A_1<\min A_2<\ldots <\min A_k$. Then, there exist sets
$A_i^{\prime }\in {\cal F}_n$, $i=1,\ldots ,k$ such that
$A_1^{\prime }<A_2^{\prime }<\ldots <A_k^{\prime }$,
$\min A_i\leq\min A_i^{\prime }$ for $i=1,\ldots ,k$, and
$\cup_{i=1}^kA_i^{\prime }=\cup_{i=1}^kA_i$.}

\medskip

\noindent {\it Proof of the claim:} It is done by induction on $n$.
For $n=0$ it is trivial. Suppose it is true for $n$.

Let $A_i$, $i=1,\ldots ,k$ be sets in ${\cal F}_{n+1}$ such that
$A_i\cap A_j=\emptyset $ for $i\neq j$ and
$\min A_1<\min A_2<\ldots <\min A_k$. Each $A_i$ is of the form
$A_i=\cup_{j=1}^{m_i}B_j^i$ where $B_j^i\in {\cal F}_n$ and,
for each $i$, $m_i\leq B_1^i<B_2^i<\ldots <B_{m_i}^i$.
Let $\{ B_j\}_{j=1}^{m_1+\ldots +m_k}$ be a rearrangement of
the family $\{ B_j^i:\;i=1,\ldots ,k,\;j=1,\ldots ,m_i\}$, which
satisfies $\min B_1<\min B_2<\ldots <\min B_{m_1+\ldots +m_k}$.
It is easy to see that, for each $i$,
$$\min A_i=\min B_1^i\leq\min B_{m_1+\ldots +m_{i-1}+1}.\leqno (\ast )$$

By the inductive assumption, there exist sets
$B_j^{\prime },j=1,\ldots ,m_1+\ldots +m_k$, with
$B_j^{\prime }\in {\cal F}_n$,
$\cup_{j=1}^{m_1+\ldots +m_k}B_j^{\prime }=\cup_{j=1}^{m_1+\ldots +m_k}B_j$
and such that $B_1^{\prime }<B_2^{\prime }<\ldots
<B_{m_1+\ldots +m_k}^{\prime }$ and $\min B_j\leq\min B_j^{\prime }$
for all $j=1,\ldots ,m_1+\ldots +m_k$. For $i=1,\ldots ,k$, we set
$$A_i^{\prime }=\cup_{j=m_1+\ldots +m_{i-1}+1}^{m_1+\ldots +m_i}B_j^{\prime }.
$$
\noindent Then, $A_1^{\prime }<A_2^{\prime }<\ldots <A_k^{\prime }$,
$\cup_{i=1}^kA_i^{\prime }=\cup_{i=1}^kA_i$, and for each
$i=1,\ldots ,k$ we have by ($\ast $),
$$m_i\leq\min B_{m_1+\ldots +m_{i-1}+1}\leq
\min B^{\prime }_{m_1+\ldots +m_{i-1}+1},$$
\noindent so $A_i^{\prime }\in {\cal F}_{n+1}$. Moreover, using
($\ast $) again, we see that
$$\min A_i\leq\min B^{\prime }_{m_1+\ldots +m_{i-1}+1}=\min A^{\prime }_i.$$
\noindent This completes the proof of the Claim. The Lemma follows.   $\Box $ 

\medskip

\noindent {\bf Distortion.} Let $\lambda >1$. A Banach space $X$
is $\lambda $--{\it distortable} if there exists an equivalent
norm $|\;.\;|$ on $X$ such that, for every infinite dimensional
subspace $Y$ of $X$,

$$\sup\{\frac{|y|}{|z|}:y,z\in Y,\|y\|=\|z\|=1\}\geq\lambda .$$

\noindent $X$ is {\it arbitrarily distortable} if it is $\lambda $--distortable
for every $\lambda >1$.

\bigskip
%\newpage
\noindent {\bf B. Mixed Tsirelson spaces.}

\medskip

A Banach space $X$
with a basis $(e_i)_{i=1}^{\infty}$ is an {\it 	asymptotic}
 $\ell_1$ space if there
exists a constant $C$ such that, for all $n$ and all block sequences
$(x_i)_{i=1}^n$ in $X$ with $n\leq x_1<x_2<\ldots <x_n,$

$$\frac{1}{C} \sum \| x_i\| \leq
 \|\sum_{i=1}^nx_i\|
.$$

The first example of an asymptotic $\ell_1$ space not containing $\ell_1$ was
constructed by Tsirelson ([19]). Tsirelson's space is the completion of
the vector space $c_{00}$ of all eventually zero sequences under the
norm $\| \cdot \|_T$ defined implicitly as follows:

$$\|x\|_T=\max\left\{\|x\|_{\infty},\sup\{ \frac{1}{2}\sum_{i=1}^n\| E_ix\|_T:
n\in {\bf N}\; {\rm and}\; n\leq E_1<E_2<\ldots <E_n\}\right\}.$$

A sequence $(E_i)_{i=1}^n$ of finite subsets of ${\bf N}$ with 
$n\leq E_1<E_2<\ldots <E_n$ is called {\it Schreier admissible} (or 
${\cal S}$-admissible). In other words, a sequence $(E_i)_{i=1}^n$
is Schreier admissible if the $E_i$'s are successive and 
$\{ \min E_1,\ldots ,\min E_n\}\in {\cal S}.$ More generally, we give the 
following definition.

\medskip

\noindent {\bf 1.3 Definition.}  Let ${\cal M}$ be a family of finite subsets
of ${\bf N}.$ 

(a) A finite sequence $(E_i)_{i=1}^n$ of subsets of ${\bf N}$
is ${\cal M}$-{\it admissible} if $E_1<E_2<\ldots < E_n$ and 
$\{\min E_1,\ldots ,\min E_n\} \in {\cal M}.$

(b) A finite sequence $(x_i)_{i=1}^n$ of vectors in $c_{00}$ is 
${\cal M}$-{\it admissible} if the sequence $({\rm supp}x_i)_{i=1}^n$
is ${\cal M}$-admissible.

\bigskip

The mixed Tsirelson spaces are defined as follows:

\noindent {\bf 1.4 Definition.} Let $\{ {\cal M}_n\}_{n=1}^{\infty}$ be 
a sequence of compact families of finite subsets of ${\bf N}$ and let
$(\theta_n)_{n=1}^{\infty}$
 be a sequence of numbers in $(0,1)$ with $\theta_n \rightarrow
0.$ The {\it mixed Tsirelson space} $T[({\cal M}_n,\theta_n)_{n=1}^{\infty}]$
is the completion of $c_{00}$ under the norm $\| \cdot \|$ defined 
implicitly by:

$$\|x\|=\max\left\{\|x\|_{\infty},\sup_k\sup\{ \theta_k\sum_{i=1}^n\| E_ix\|:
n\in {\bf N}\; {\rm and}\; (E_i)_{i=1}^n \; {\rm is} \; {\cal M}_k 
{\rm -admissible}\}\right\}.$$

\medskip

The mixed Tsirelson spaces $T[({\cal M}_n,\theta _n)_{n=1}^{\infty}]$ where
$({\cal M}_n)_n$ is a subsequence of the sequence of Schreier families
$({\cal F}_j)_{j=1}^{\infty}$ were introduced in [3] and further studied
in [2] and [14]. Every such space is a reflexive asymptotic $\ell_1$
Banach space and the natural basis $(e_i)_i$
is a 1--unconditional basis for it. The first example of an 
arbitrarily distortable asymptotic $\ell_1$ Banach space was a space of this
type ([3]). More generally, Androulakis and Odell have proved the
following:

\medskip

\noindent {\bf 1.5 Theorem.} ([2]) {\it Suppose that the sequence
$(\theta_n)_n$ satisfies $\theta_{n+m}\geq \theta_n\theta_m$ for all $n,m$
and let $\theta=\lim {\theta^{1/n}_n}.$ If $\frac{\theta_n}{\theta^n}
\rightarrow 0$ then the space 
$T[({\cal F}_n,\theta_n)_{n=1}^{\infty}]$
is arbitrarily distortable.}   $\Box $

\medskip

In particular, this is  the case if $\lim \theta_n^{1/n}=1.$ The first
result of this section concerns mixed Tsirelson spaces $T[({\cal F}_n,\theta_n)
_n]$ corresponding to such sequences $(\theta_n)_n.$ Following [2] we
call a sequence $(\theta_n)_n$ {\it regular}, if $\theta_n\in (0,1)$ for
all $n$, $\theta_n\downarrow 0 $ and $\theta_{n+m}\geq \theta_n\theta_m$
for all $n,m \in {\bf N}.$

\medskip

\noindent {\bf 1.6 Theorem.} {\it Let $(\theta_n)_{n=1}^{\infty}$ be a regular
sequence with $\lim \theta^{1/n}_n=1.$ Let 
$X=T[({\cal F}_n,\theta_n)_{n=1}^{\infty}].$ For every $\varepsilon >0$,
every infinite dimensional block subspace $Y$ of $X$ contains for
every $n$ a sequence of disjointly supported vectors $(y_i)_{i=1}^n$ which is
$(1+\varepsilon)$ -- equivalent to the canonical basis of $\ell_{\infty}^n.$}

\medskip

Given a block subspace $Y$ of $X$ and $n\in N$ we shall construct a
sequence $(x_i)_{i=1}^n$ of disjointly supported normalized
vectors in $Y$ such that
$\| \sum _{i=1}^nx_i\|\leq 36.$ Since the basis $(e_n)_n$ of $X$ is 
1-unconditional this implies that $(x_i)_{i=1}^n$ is 36-equivalent to
the canonical basis of $\ell_{\infty}^n.$ From this the Theorem follows
by a standard argument due to R.C. James.
The building blocks of our construction  are the $(\varepsilon, j)$--
{\it rapidly increasing
special convex combinations}, the prototypes of which were used in [3].
Before proceeding to the construction we need to establish some preliminary
results most of which also have their analogues in [3].

\medskip

\noindent {\bf Notation.}  Let $X=T[({\cal F}_n,\theta_n)_{n=1}^{\infty}].$

\noindent A. Inductively, we define a subset $K=\cup_{n=0}^{\infty}K^n$
of $B_{X^*}$ as follows:

For $j=1,2,\ldots$, 

$$K^0_j=\{\pm e_n: n\in{\bf N}\}.$$

Assume that $K^n_j,\; j=1,2,\ldots$ have been defined. We set $K^n=\cup_{j=1}
^{\infty}K^n_j$ and, for $j=1,2,\ldots$, we set 

$$K^{n+1}_j=K^n_j\cup \{ \theta_j(f_1+\cdots +f_d): d\in {\bf N}, f_i\in K^n,
i=1,\ldots,n,$$
$${\rm supp f_1}<\ldots <{\rm supp}f_d \; {\rm and} \; (f_i)_{i=1}^d\;
{\rm is}\; {\cal F}_j-{\rm admissible}\}.$$

Let $K=\cup_{n=0}^{\infty}K^n.$ 

Then $K$ is a norming set for $X$, that is, for $x\in X$

$$\| x\|=\sup\{ f(x): f\in K\}.$$

\noindent B. For $j=1,2,\ldots,$ we denote by ${\cal A}_j$ the set
${\cal A}_j=\cup_{n=1}^{\infty}(K^n_j\setminus K^0).$

\noindent C. Let $m\in {\bf N},$ $\varphi \in K^m\setminus K^{m-1}.$
An {\it analysis} of $\varphi$ is a family $\{ K^s(\varphi )\}_{s=0}^m$
of subsets of $K$ such that

(1) For every $s\leq m$, $K^s(\varphi )\subset K^s,$ the elements of
$K^s(\varphi )$ are disjointly supported and $\cup_{ f\in K^s(\varphi )}
{\rm supp} f={\rm supp}\varphi.$

(2) If $f$ belongs to $K^{s+1}(\varphi )$ then either $f\in K^s(\varphi )$
or, for some $j\geq 1$, there exists a ${\cal F}_j$-admissible family 
$(f_i)_{i=1}^d$ in $K^s(\varphi )$ such that $f=\theta _j(f_1+\cdots + f_d).$

(3) $K^m(\varphi )=\{\varphi\}.$

\medskip

It is easy to see that every $\varphi \in K$ has an analysis.

\medskip

\noindent {\bf 1.7 Definition}. Let $n\geq 1$, $\varepsilon >0$ and
$F\subseteq {\bf N}$, $F\in {\cal F}_n$. A convex combination
$\sum_{k\in F}a_ke_k$ is called an $(\varepsilon ,n)$-{\it basic
special convex combination} (basic s.c.c) if, for every
$G\in {\cal F}_{n-1}$, $\sum_{k\in G}a_k<\varepsilon $.

\medskip

\noindent {\bf 1.8 Proposition.} {\it Let $D$ be an infinite subset
of ${\bf N}$. Then, for every $n\geq 1$ and $\varepsilon >0$, there
exists an $(\varepsilon ,n)$-basic special convex combination
$x=\sum_{k\in F}a_ke_k$ with $F={\rm supp}(x)\subset D$.}

\medskip

\noindent {\it Proof:} For $n=1$, we choose $m_0>\frac{1}{\varepsilon }$
and $A\subset D$ with $m_0<A$ and $|A|=m_0$. Then,
$x=\frac{1}{m_0}\sum_{k\in A}e_k$ is an $(\varepsilon ,1)$-basic
s.c.c.

For $n>1$ the proof is by induction based on the following:

\medskip

\noindent {\bf 1.9 Lemma}. {\it Let $n\geq 1$ and suppose that the
integers $m_0, m_1,\ldots ,m_{m_0}$ and the block vectors
$x_1, x_2,\ldots ,x_{m_0}$ satisfy the following: For every
$k=1,2,\ldots ,m_0-1$,

{\rm (a)} $2m_{k-1}<m_k$. 

{\rm (b)} ${\rm supp}(x_k)\subset (m_{k-1},m_k]$.

{\rm (c)} $x_k$ is a $(\frac{1}{2m_{k-1}},n)$-basic s.c.c.

\noindent Then, the vector $x=\frac{1}{m_0}\sum_{k=1}^{m_0}x_k$ is a
$(\frac{2}{m_0}, n+1)$-basic s.c.c.   }

\medskip

\noindent {\it Proof:} The proof is straightforward. (See also Lemma 1.6 
[3]).  $\Box$
\bigskip

\noindent {\bf 1.10 Definition}. Let $\varepsilon >0$, $j\in {\bf N}$
and suppose that $\{z_k\}_{k=1}^n$ is a finite block sequence
with the property that there exist integers $\{l_k\}_{k=1}^n$ with
$2<z_1\leq l_1<z_2\leq l_2<\ldots \leq l_{n-1}<z_n\leq l_n$, and such
that a convex combination $\sum_{k=1}^na_ke_{l_k}$ is an
$(\varepsilon ,j)$-basic s.c.c.
Then, the corresponding convex combination of the $z_k$'s,
$x=\sum_{k=1}^na_kz_k$, is called an $(\varepsilon ,j)$-{\it s.c.c.
of $\{z_k\}_{k=1}^n$}.

An $(\varepsilon ,j)$-s.c.c. $x=\sum_{k=1}^na_kz_k$ of {\it unit} vectors
$\{z_k\}_{k=1}^n$ is said to be {\it seminormalized} if $\|x\|\geq\frac{1}{2}$.

\medskip

\noindent {\bf Remark:} It is easy to see that if $x=\sum_{k=1}^na_kz_k$ is an
$(\varepsilon ,j)$-s.c.c and $\|z_k\|=1,k=1,\ldots ,n$, then
$\|x\|\geq\theta_{j+1}$. Indeed, if $f_k\in B_{X^{\ast }}$
are chosen so that $f_k(z_k)=\|z_k\|=1$, ${\rm supp}(f_1)\subset (2,l_1]$, and
${\rm supp}f_k\subset (l_{k-1},l_k]$ for $k=2,\ldots ,n$, then the family $\{f_k\}_k$
is ${\cal F}_{j+1}$-admissible. This implies that the functional
$\varphi =\theta_{j+1}\sum f_k$ belongs to $B_{X^{\ast }}$,
hence $\|x\|\geq\varphi (x)\geq\theta_{j+1}$. 

\medskip

The following Lemma states that
every block subspace $Y$ of $X$ contains for any $\varepsilon $ and $j$
a seminormalized $(\varepsilon ,j)$-s.c.c. The condition $\lim \theta^{1/j}_j
=1$ is essential at this point. 

\medskip

\noindent {\bf 1.11 Lemma}. {\it Let $j\in {\bf N}$, $\varepsilon >0$ and let
$\{z_k\}_{k=1}^{\infty }$ be a block sequence in $X$. There exists
$n\in {\bf N}$ and normalized blocks $y_k, k=1,\ldots ,n$ of the sequence
$\{z_k\}_{k=1}^{\infty }$ such that a convex combination
$x=\sum_{k=1}^na_ky_k$ is a seminormalized $(\varepsilon ,j)$-s.c.c.}

\medskip

\noindent {\it Proof:} We may assume that the vectors $z_k, k=1,2,\ldots $ are
normalized. Choose an infinite block sequence $\{x_l^1\}_{l=1}^{\infty }$
of $\{z_k\}_{k=1}^{\infty }$ such that, for each $l$,
$x_l^1=\sum_{k\in A_l}a_kz_k$ is an $(\varepsilon ,j)$-s.c.c of
$\{z_k\}_{k\in A_l}$.

If for some $l$, $\|x_l^1\|\geq\frac{1}{2}$, then we are done. If not,
we set $y_l^1=\frac{x_l^1}{\|x_l^1\|}$ and, as before, choose an
infinite sequence $\{x_l^2\}_l$ of $(\varepsilon ,j)$-s.c.c
of $\{y_l^1\}_{l=1}^{\infty }$.

Notice that, for each $l$, the family $\{z_k:\;{\rm supp}(z_k)\subset
{\rm supp}(x_l^2)\}$ is ${\cal F}_{2j+2}$-admissible
(since ${\cal F}_{2j+2}={\cal F}_{j+1}[{\cal F}_{j+1}]$), and so
$x_l^2$ is a combination of the form $x_l^2=\sum b_k(\lambda_kz_k)$
where $\sum b_k=1,\lambda_k\geq 2$, and $\{z_k\}$ is an
${\cal F}_{2j+2}$-admissible family. This gives that
$\|x_l^2\|\geq 2\theta_{2j+2}$.

If, for some $l$, $\|x_l^2\|\geq\frac{1}{2}$ then we are done. If not,
then we set $y_l^2=\frac{x_l^2}{\|x_l^2\|}$ and continue as before.

Continuing in this manner, if we never get some $(\varepsilon ,j)$-s.c.c
$x_l^k$ with $\|x_l^k\|\geq\frac{1}{2}$, then we can repeat
the same procedure for as many steps $s$ as we wish and always get $1\geq
\|x_l^s\|\geq 2^{s-1}\theta_{s(j+1)}$.

But the assumption that $\lim_n\theta_n^{\frac{1}{n}}=1$ implies
that $\lim_{s\rightarrow\infty }2^{s-1}\theta_{s(j+1)}=\infty $.
This leads to a contradiction which completes the proof.    $\Box $

\bigskip

\noindent {\bf 1.12 Lemma}. {\it Let $x=\sum_{l\in F}a_le_l$, 
where $F\in {\cal F}_j$,
be an $(\varepsilon ,j)$-basic s.c.c. Then, $\theta_j\leq\|x\|<\theta_j+\varepsilon $.}

\medskip

\noindent {\it Proof:} It is obvious that $\varphi =\theta_j
(\sum_{l\in F}e_l^{\ast })$ belongs to $B_{X^{\ast }}$ and $\varphi (x)=\theta_j$.
This yields the lower estimate for $\|x\|$.

It remains to prove that, for all $\psi\in K$, $|\psi (x)|\leq
\theta_j+\varepsilon $. Let $\psi\in K$; we may assume that $\psi$ is 
positive. Set
$$J=\{l\in F:\psi (e_l)\leq \theta_j\}$$
\noindent and
$$L=F\backslash J=\{l\in F:\psi (e_l)>\theta_j\}.$$

\noindent We shall prove that $L\in {\cal F}_{j-1}$ and so
$\sum_{k\in L}a_k<\varepsilon $. This is a consequence of the following:

\medskip

\noindent {\it Claim:} Let $r=1,2,\ldots ,\;f\in K$ and suppose
that $f(e_k)>\theta_r$ for all $k\in {\rm supp}(f)$. Then,
${\rm supp}(f)\in {\cal F}_{r-1}$.

\medskip 

\noindent {\it Proof of the Claim:} The proof is by induction on $s$,
for $f\in K^s,s=1,2,\ldots $.

For $s=1$, let $f\in K^1$, with $f=\theta_i\sum_{k\in A}e_k^{\ast },
\;A\in {\cal F}_i$. Since $\theta_i>\theta_r$, we get $i\leq r-1$
and so $A={\rm supp}(f)\in {\cal F}_{r-1}$.

Suppose that the claim is true for all $g\in K^s$ and let $f\in K^{s+1}$.
Then, $f=\theta_i(\sum_{l=1}^mf_l)$ where the set
$(f_l)_{l=1}^m$ is ${\cal F}_i$-admissible and, for each $l$,
$f_l\in K^s$. Suppose that $f(e_k)>\theta_r$ for all $k\in {\rm supp}(f)$.
Then, $r>i$ and, for each $l=1,\ldots ,m$, $f_l(e_k)>\frac{\theta_r}{\theta_i}
\geq\theta_{r-i}$. It follows from the inductive hypothesis that
${\rm supp}(f_l)\in {\cal F}_{r-i-1},l=1,\ldots ,m$.
So, ${\rm supp}(f)\in {\cal F}_i[{\cal F}_{r-i-1}]={\cal F}_{r-1}$.
This completes the proof of the claim.

\medskip

We conclude that $L\in {\cal F}_{j-1}$ and so
$$|\psi (\sum_{l\in F}a_le_l)|\leq \psi (\sum_{l\in J}a_le_l)+
\sum_{l\in L}a_l<\theta_j+\varepsilon.\;\;\;\;\Box $$

\medskip

\noindent {\bf 1.13 Lemma.} {\it Let $x=\sum_{k=1}^na_ky_k$ be an
$(\varepsilon ,j)$-s.c.c, where $\varepsilon <\theta_j$.
Let $i<j$ and suppose that $(E_r)_{r=1}^s$ is an ${\cal F}_i$-admissible
family of intervals. Then,}
$$\sum_{r=1}^s\|E_rx\|\leq (1+\frac{\varepsilon }{\theta_i})\max_{1\leq k\leq n}\|y_k\|\leq 2\max_{1\leq k\leq n}\|y_k\|.$$

\medskip

\noindent {\it Proof:} We can assume that the $E_r$'s are adjacent
intervals. Set
$$L=\{k: k=1,\ldots ,n\;{\rm and}\;{\rm supp}(y_k)\;{\rm is}
\;{\rm intersected}\;{\rm by}\;{\rm at}\;{\rm least}\;{\rm two}
\;{\rm different}\;E_r{\rm 's}\}.$$

\noindent For each $r=1,\ldots ,s$, define 
$$B_r=\{k:k=1,\ldots ,n\;{\rm and}\;{\rm supp}(y_k)\subset E_r\}.$$

\noindent The sets $B_r$ are mutually disjoint and $\{1,2,\ldots ,n\}
=(\cup_{r=1}^sB_r)\cup L$. So,
$$\sum_{r=1}^s\|E_rx\|\leq \sum_{r=1}^s\|E_r(\sum_{k\in B_r}a_ky_k)\|
+\sum_{k\in L}a_k\sum_{r=1}^s\|E_ry_k\|$$
$$\leq\sum_{k=1}^na_k\|y_k\|+\sum_{k\in L}a_k\frac{\|y_k\|}{\theta_i}.$$

Suppose now that $2<y_1\leq l_1<y_2\leq \ldots\leq l_{k-1}
<y_k\leq l_k$ and $\sum_{k=1}^na_ke_{l_k}$ is the basic s.c.c
which defines the s.c.c $x=\sum_{k=1}^na_ky_k$. We shall show that
$\{l_k:k\in L\}\in {\cal F}_i\subset {\cal F}_{j-1}$. This will
imply that $\sum_{k\in L}a_k<\varepsilon $ and hence complete
the proof.

To see that $\{l_k:k\in L\}\in {\cal F}_i$, for each $k\in L$ let
$r_k=\min\{r:E_r\;{\rm intersects}\;{\rm supp}(y_k)\}$. The map
$k\rightarrow r_k$ from $L$ to $\{1,2,\ldots ,s\}$ is one to one.
This gives that $\# L\leq s$. Consider now, for each $k\in L$,
$m_{r_k}=\min E_{r_k}$. Then, $m_{r_k}\leq l_k,k\in L$. Since the
set $\{m_{r_k}:k\in L\}$ belongs to ${\cal F}_i$, we conclude
(by the spreading property of ${\cal F}_i$) that $\{l_k:k\in L\}\in
{\cal F}_i$ as well.   $\Box $

\bigskip

\noindent {\bf 1.14 Definition.} A. A finite or infinite sequence
$\{z_k\}_k$ is called a {\it rapidly increasing sequence} if there
exists an increasing sequence of positive integers $\{t_k\}_k$
such that the following are satisfied:

(a) The sequence $\{
\frac{\theta_{t_k}}{\theta_{t_{k+1}}}
\}_k$ is
increasing,  $2<\theta_{t_k}/\theta_{t_{k+1}}$ for each $k$, and
$\lim_{k\rightarrow \infty}\frac{\theta_{t_k}}{\theta_{t_{k+1}}}=\infty$ if
the sequence is infinite.

(b) Each $z_k$ is a semi-normalized $(\theta_{t_k}^2,t_k)$-s.c.c.

(c) For each $k$, $\|z_k\|_{\ell_1}\leq\frac{\theta_{t_k}}{\theta_{t_{k+1}}}$.

\medskip

B. Let $j\in {\bf N}$, $\varepsilon >0$. Let $\{z_k\}_{k=1}^n$ be a
rapidly increasing sequence, where each $z_k$ is a seminormalized
$(\theta_{t_k}^2,t_k)$-s.c.c and $2<\frac{\theta_{j+1}}{\theta_{t_1}}
<\frac{\theta_{t_1}}{\theta_{t_2}}$. Suppose also that there exist
coefficients $\{a_k\}_{k=1}^n$ such that the vector
$x=\sum_{k=1}^na_kz_k$ is an $(\varepsilon ,j)$-s.c.c of
$\{z_k\}_{k=1}^n$. Then $x$ is called an $(\varepsilon ,j)$-{\it rapidly
increasing special convex combination} ($(\varepsilon ,j)$-R.I.s.c.c).

\bigskip

\noindent {\bf 1.15 Proposition.}
 {\it Let $j\in {\bf N}, 0<\varepsilon <\theta_j^2$,
and let $x=\sum_{k=1}^na_kz_k$ be an $(\varepsilon ,j)$-R.I.s.c.c
of the $z_k$'s, where each $z_k$ is a seminormalized
$(\theta^2_{t_k},t_k)$-s.c.c. Let  $t_0$ be any integer
such that $j+1\leq t_0<t_1$ and $2<\frac{\theta_{t_0}}{\theta_{t_1}}.$

Then, for every $\varphi$
in the norming set $K$ of $X$, we have the following estimates:

  {\rm (i)} $|\varphi (x)|\leq 8\theta_j$,   if $\varphi\in {\cal A}_i,i<j$

 {\rm (ii)} $|\varphi (x)|\leq 4\theta_i$,   if $\varphi\in {\cal A}_i,j\leq i<t_1$

{\rm (iii)} $|\varphi (x)|\leq 4(\theta_{t_{p-1}}+a_{t_p})$,   if
$\varphi\in {\cal A}_i,\;  t_p\leq i<t_{p+1},\; p\geq 1$.

\noindent In particular, $\frac{\theta_{j+1}}{2}\leq\|x\|\leq 8\theta_j$.}

\medskip

\noindent {\it Proof:} The lower estimate for $\|x\|$ follows by the
Remark after Definition 1.10 and the fact that $\|z_k\|\geq\frac{1}{2}$.
The upper estimate follows from the first part of the Proposition. 
The proof of this is similar to the one of Proposition
2.12 in [3]. Let $\{l_k\}_{k=1}^n$ be such that $2<z_1\leq l_1<
\ldots \leq l_{n-1}<z_n\leq l_n$ and $\sum_{k=1}^na_ke_{l_k}$ is an
$(\varepsilon ,j)$-basic s.c.c. 

\medskip

\noindent Given $\varphi\in K$, we shall construct $\psi\in {\rm co}(K)$ such that

(a) $\varphi (\sum_{k=1}^na_kz_k)\leq 4\psi (\sum_{k=1}^na_ke_{l_k})$.

(b) If $\varphi \in {\cal A}_i,i<t_1$, then $\psi\in {\rm co}({\cal A}_i)$.

(c) If $\varphi\in {\cal A}_i,t_p\leq i<t_{p+1}$ for some $p\geq 1$,
then $\psi =\frac{1}{2}(\psi_1+e_{l_p}^{\ast })$, where
$\psi_1\in {\rm co}({\cal A}_{t_{p-1}})$.

\medskip

\noindent Since, for $\psi\in {\rm co}({\cal A}_i)$ we have
$\psi (\sum a_ke_{l_k})\leq\theta_i$, estimates (ii) and (iii)
will follow immediately. For (i) we apply Lemma 1.12.

\medskip

We consider an analysis $\{K^s(\varphi )\}_{s=1}^m$ of $\varphi $, and
we cut each $z_k$ into two parts, $z_k^{\prime }$ and $z_k^{\prime\prime }$, with the following property:

\medskip

{\it {\rm ($\ast $)} For each level $K^s(\varphi )$ of the analysis of $\varphi $, and for each
$z_k^{\prime }$, either there exists a unique $f\in K^s(\varphi )$ with
${\rm supp}(z^{\prime }_k)\cap {\rm supp}(f)\neq\emptyset $ or there exists
$f\in K^s(\varphi )$ such that $\max {\rm supp}(z^{\prime }_{k-1})<{\rm supp}(f)<
\min {\rm supp}(z^{\prime }_{k+1})$.}

\medskip

\noindent The same is true for $z_k^{\prime\prime}$. This partition of the $z_k$'s
is possible, as done in [3] (Definition 2.4).

\medskip

We shall see that using property ($\ast $) we can build $\psi^{\prime }$
and $\psi^{\prime\prime }$ such that $|\varphi (z_k^{\prime })|
\leq \psi^{\prime }(e_{l_k})$ and $|\varphi (z_k^{\prime\prime })|
\leq \psi^{\prime\prime }(e_{l_k})$ for all $k$. So we may assume that
the $z_k$'s have property ($\ast $) and then multiply our
estimate by 2.

For each $f\in\cup_{s=0}^mK^s(\varphi )$ we set
$$D_f=\{k:{\rm supp}(\varphi )\cap {\rm supp}(z_k)=
{\rm supp}(f)\cap {\rm supp}(z_k)\neq\emptyset\}.$$

\noindent By induction on $s=0,\ldots ,m$ we shall define a function
$g_f\in {\rm co}(K)$, supported on $\{l_k:k\in D_f\}$ and such
that:

(a) $|f(z_k)|\leq 2g_f(e_{l_k})$ for all $k\in D_f$.

(b) If $f\in {\cal A}_q,q<t_1$, then $g_f\in {\rm co}({\cal A}_q)$. If
$f\in {\cal A}_q,t_p\leq q<t_{p+1}$, then
$g_f=\frac{1}{2}(g_f^1+e_{l_p}^{\ast })$, where $g_f^1\in {\rm co}(
{\cal A}_{t_{p-1}})$.

\medskip

For $s=0,f=e_r^{\ast }$, if $D_f=\{k\}$ we set $g_f=e_{l_k}^{\ast }$.

Let $s>0$. Suppose that $g_f$ has been defined for all
$f\in\cup_{t=0}^{s-1}K^t(\varphi )$. Let
$$f=\theta_q(f_1+\ldots +f_d)\in K^s(\varphi )\backslash K^{s-1}(\varphi ).$$

\noindent We set $I=\{i:1\leq i\leq d, D_{f_i}\neq\emptyset\}$ and
$T=D_f\backslash\cup_{i\in I}D_{f_i}$.

\medskip

\noindent {\it Case 1:} $q<t_1$.

Then, we set
$$g_f=\theta_q\left(\sum_{i\in I}g_{f_i}+\sum_{k\in T}e_{l_k}^{\ast }\right).$$

\noindent Property (a) for the case $k\in\cup_{i\in I}D_{f_i}$
follows from the inductive assumption. For $k\in T$ we get, by Lemma 1.13,
since $q<t_k$, that
$$|f(z_k)|\leq\theta_q\sum_{i=1}^d|f_i(z_k)|\leq 2\theta_q=2g_f(e_{l_k}).$$

\noindent To prove that $g_f\in {\rm co}({\cal A}_q)$ we need to show
that the set $\{g_{f_i}:i\in I\}\cup\{l_k:k\in T\}$
is ${\cal F}_q$-admissible.

Here we use property ($\ast $). According to ($\ast $), for each
$k\in T$ there exists an $i_k\in\{1,\ldots ,d\}$ such that
$\max {\rm supp}(z_{k-1})<{\rm supp}(f_{i_k})<\min {\rm supp}(z_{k+1})$.

This means that $i_k\neq i_l$ for $k\neq l\in  T$ and $i_k\notin I$. It
follows that $|T|+|I|\leq d$. Since also, for each $k\in T$,
$\min {\rm supp}(f_{i_k})\leq l_k$, by the spreading property of
${\cal F}_q$ we get that
$$\{\min {\rm supp}(f_i):i\in I\}\cup\{l_k:k\in T\}\in {\cal F}_q,$$

\noindent hence the family $\{g_{f_i}\}_{i\in I}\cup\{e_{l_k}^{\ast }\}_{k\in I}$ is ${\cal F}_q$-admissible.

\medskip

\noindent {\it Case 2:} $q\geq t_1$.

Suppose that $t_p\leq q<t_{p+1}$. If $p\notin D_f$ or $p\in\cup_{i\in I}D_{f_i}$,
then we set
$$g_f=\theta_{t_{p-1}}\left(\sum_{i\in I}g_{f_i}+\sum_{k\in T}e_{l_k}^{\ast }\right).$$

Since ${\rm supp}(g_f)\subset\{l_k:k=1,\ldots ,n\}\in {\cal F}_j$
and $j<t_{p-1}$, it is clear that $g_f\in {\rm co}({\cal A}_{t_{p-1}})$.

\medskip

For $k\in\cup_{i\in I}D_{f_i}$ we get
$$|f(z_k)|=\theta_q|f_i(z_k)|<2\theta_qg_{f_i}(e_{l_k})
<\theta_{t_{p-1}}g_{f_i}(e_{l_k})=g_f(e_{l_k})$$

\noindent by the inductive assumption and the fact that
$2\theta_{t_p}<\theta_{t_{p-1}}$.

\medskip

For $k\in T$, $k<p$, we have
$$|f(z_k)|\leq\theta_q\sum_{i=1}^d|f_i(z_k)|\leq\theta_q\|z_k\|_{\ell_1}
\leq\theta_q\frac{\theta_{t_k}}{\theta_{t_{k+1}}}
\leq\theta_{t_p}\frac{\theta_{t_{p-1}}}{\theta_{t_p}}=
\theta_{t_{p-1}}=g_f(e_{l_k})$$

\noindent by the property of the R.I.S $\{z_k\}_k$.

\medskip

For $k\in T$, $k>p$, we have $q<t_{p+1}\leq t_k$, so
$$|f(z_k)|=\theta_q\sum_{i=1}^d|f_i(z_k)|\leq 2\theta_q
<\theta_{t_{p-1}}=g_f(e_{l_k})$$

\noindent by Lemma 1.13.

\medskip

Suppose now that $p\in T$. Then we set
$$g_f=\frac{1}{2}\left[\theta_{t_{p-1}}\left(\sum_{i\in I}g_{f_i}+
\sum_{k\in T\backslash\{p\}}e_{l_k}^{\ast }\right)+
e_{l_p}^{\ast }\right].$$

\noindent As before, we get
$$|f(z_k)|<2g_f(e_{l_k})$$

\noindent for $k\neq p$, and
$$|f(z_p)|\leq 1=2g_f(e_{l_p}).$$

\noindent This completes the inductive step of the construction
and the proof of the Proposition.   $\Box $

\bigskip

In what follows, a finite tree of sequences ${\cal T}$ will be a finite
set of finite sequences of positive integers, partially ordered by the
relation: $\alpha \prec\beta $ iff $\alpha $ is an initial part of $\beta $,
and satisfying the following properties:

(a) For each $\alpha\in {\cal T}$, the set
$\{\beta :\beta \;{\rm is}\;{\rm an}\;{\rm initial}\;{\rm part}\;
{\rm of}\;\alpha\}$ is contained in ${\cal T}.$

(b) If $\alpha =(k_1,\ldots ,k_{m-1},k_m)\in {\cal T}$ and
$1\leq l\leq k_m$, then $(k_1,\ldots ,k_{m-1},l)\in {\cal T}$.

(c) The maximal (under $\prec $) elements of ${\cal T}$ are all of the same
length.

\medskip

It follows that ${\cal T}$ has a unique {\it root}, the empty
sequence which we denote by 0. The length of the sequence $\alpha $
is denoted by $|\alpha |$. The {\it height } of ${\cal T}$ is the length
of the maximal elements of ${\cal T}$. For each $\alpha\in {\cal T}$
which is not maximal we set $S_{\alpha }=\{\beta\in {\cal T}:
\alpha \prec \beta\;{\rm and}\;|\beta |=|\alpha |+1\}$. We also
consider the lexicographic order, denoted by $<$, on ${\cal T}$.
For $\alpha =(k_1,\ldots,k_{m-1},k_m)\in{\cal T}$ we denote by  $\alpha^+$
the sequence $\alpha^+=(k_1,\ldots, k_{m-1},k_m+1).$

\medskip

\noindent {\bf 1.16 Definition}. Let $r\in {\bf N}$. 
Let $j_1,\ldots ,j_r$ be positive
integers, and $\varepsilon >0.$ 
 An $(\varepsilon ,(j_1,\ldots ,j_r)
)$ -{\it tree in} $X$ is a set of vectors ${\cal T}^X=
\{u_{\gamma }\}_{\gamma\in {\cal T}}$ indexed by a finite
tree ${\cal T}$ of height $r$, and satisfying the following
properties:

(a) The terminal nodes $\{u_{\alpha }\}_{|\alpha |=r}$ of the tree
are elements of the basis $\{e_n\}_{n=1}^{\infty }$, i.e,
for $|\alpha |=r$, $\alpha\in {\cal T}$, $u_{\alpha }=e_{l_{\alpha }}$.
Moreover,  
for $\alpha,  \beta \in {\cal T}$  with $|\alpha |=|\beta |=r,$
if $\alpha <\beta $ (in the lexicographic order), then $l_{\alpha}<l_{\beta}.$

(b) There exist positive coefficients $\{a_{\beta }\}_{\beta\in {\cal T}\backslash\{0\}}$ such that, for each $\gamma\in {\cal T},|\gamma |=t<r$,
we have $\sum_{\beta\in S_{\gamma }}a_{\beta }=1$ and
$u_{\gamma }=\sum_{\alpha\in {\cal T},|\alpha |=r,\gamma \prec\alpha}
\left(\prod_{\gamma \prec\beta
\preceq\alpha }a_{\beta }\right)e_{l_{\alpha }}$ is an
$(\varepsilon,j_{t+1}+j_{t+2}+\ldots +j_r)$-basic s.c.c of
$\{e_{l_{\alpha }}\}_{\alpha\in {\cal T},|\alpha |=r}$.

\medskip

It is clear that, given an infinite subset $L$ of ${\bf N}$, $j_1,\ldots ,j_r$
positive integers, and $\varepsilon >0$,
 one can construct an $(\varepsilon,(j_1,\ldots ,j_r)
)$-tree in $X$, supported in $L$, by repeatedly applying Lemma 1.9.
It is also not hard to see in the same manner that the following
construction is possible:

\medskip

\noindent {\bf 1.17 Lemma.} {\it Let $L$ be an infinite subset of ${\bf N}$,
$n\in {\bf N}$, $\varepsilon >0$ and
$j_1,\ldots ,j_n$ be positive integers. There exist a tree of
sequences ${\cal T}$, subsets ${\cal T}_1^X,\ldots , {\cal T}_n^X$
of $X$, and positive coefficients $\{a_{\beta }\}_{\beta\in {\cal T}\backslash \{0\}}$ such that:

{\rm (a)} For $r\leq n$, set ${\cal T}_r=\{\alpha \in {\cal T}: |\alpha |\leq
r\}$.
Then, ${\cal T}_r^X=\{u_{\alpha }^{r}\}_{\alpha\in {\cal T}_r}$ is an
$(\varepsilon, (j_1,\ldots ,j_r))$-tree in $X$
with coefficients $\{a_{\beta }\}_{\beta\in {\cal T}_r \backslash \{0\}}$, 
supported
in $L$.

{\rm (b)} Let $\{e_{l^r_{\alpha }},\alpha\in {\cal T},|\alpha |=r\}$
be the terminal nodes of the tree ${\cal T}_r^X$. Then, if
$\alpha ,\beta \in {\cal T}$, $|\alpha |=r<n$ and $\beta\in S_{\alpha}$,
we have $l_{\alpha}^r<l_{\beta }^{r+1}<l_{\alpha ^+}^r.$} 

\bigskip

\noindent {\bf 1.18 Definition.} A finite family ${\cal T}_1^X,\ldots ,{\cal T}_n^X$
as described in Lemma 1.17 is called an $(\varepsilon,(j_1,\ldots,j_n))$ 
{\it family
of nested trees} in $X$.

\bigskip

\noindent {\bf Proof of Theorem 1.6.}

\medskip

Given $n\in {\bf N}$, and a block subspace $Y$ of $X$ we shall construct
a sequence $x_1,\ldots ,x_n$ of disjointly
supported unit vectors in $Y$ which is $36$-equivalent
to the canonical basis of $\ell_{\infty }^n$.

\medskip

The construction is as follows:

First, choose $\eta >0$ with $\eta <\frac{1}{60n}.$ Choose $j_0$ such that
$64\theta_{j_0}<\eta.$ Let $s_0\in {\bf N}$ be such that $\theta^{s_0}_{1}
<\eta .$ Choose $j_1$ such that

\centerline { $s_0j_0<j_1$    and    $\frac{\theta_{j_1+1}}{\theta
_{j_1}}\geq \frac{1}{1+\eta}.$}

Inductively, choose $j_2,\ldots,j_n$ so that, for each $k=2,\ldots n$,

$$j_1+\ldots +j_{k-1}<j_k,\;\; \;
\frac{8\theta_{j_k}}{\theta_{j_1+\ldots +j_{k-1}+1}}<\eta,
\;\;\; {\rm and }\;\;\; \frac{\theta_{j_1+\ldots +j_k+1}}{\theta_{j_k}}
\geq \frac{1}{1+\eta}.$$
The latter is possible, since $\lim _{n\rightarrow \infty}\theta^{1/n}_n=1$.

\medskip

Next, we choose an infinite R.I.S. $\{z_i\}_{i=1}^{\infty}$ in $Y$ where
each $z_i$ is a $(\theta^2_{t_i},t_i)$-seminormalized s.c.c. For each $i$,
let $l_i=\max ({\rm supp}z_i).$ Let $i_0$ be such that 

$$t_{i_0}>j_1+\ldots +j_n+1 \;\;\;\; {\rm and} \;\;\;\; \frac{\theta_{t_{i_0}}}
{\theta_{j_1+\ldots +j_n+1}}< \frac{\eta}{16}.$$

We set $L_0=\{ l_i\}_{i>i_0}.$

\medskip

Now let $0< \varepsilon <\min\{\theta^2_{j_1+\ldots +j_n+1},\eta (1-\theta_1)\}
.$

We choose an $(\varepsilon,(j_1,\ldots,j_n))$--
family of nested trees $({\cal T}_1^X,\ldots , {\cal T}_n^X)$
in $X$, indexed by a tree ${\cal T}$, supported in $L_0$. 
 Let $\{a_{\beta }\}_{\beta\in {\cal T}}$
be the corresponding coefficients. Then, for each $r\leq n$, there
exists a set $\{l_{\alpha }^r\}_{\alpha\in {\cal T},|\alpha |=r}$, contained in
$L_0$, and such that for all $t<r$ and $\gamma\in {\cal T}$ with
$|\gamma |=t$,
$$u_{\gamma }^r=\sum_{\gamma\prec\alpha , |\alpha |=r}(\prod_{\gamma\prec
\beta\preceq\alpha }a_{\beta })e_{l^r_\alpha }$$
\noindent is an $(\varepsilon,j_{t+1}+\ldots +j_r)$- basic s.c.c. of
$\{e_{l_{\alpha }^r}\}_{\alpha\in {\cal T},|\alpha |=r}$.

\noindent For each $\alpha \in {\cal T}$ with $| \alpha |=r$, 
denote by $z_{\alpha }^r$ the element of $\{z_i\}_{i\in {\bf N}
 }$ with
$\max {\rm supp}(z_{\alpha }^r)=l_{\alpha }^r$.

\noindent Then, for $\gamma\in {\cal T}$ with $|\gamma |=t<r$, the vector
$$y_{\gamma }^r=\sum_{\gamma\prec\alpha ,|\alpha |=r}(\prod_{\gamma
\prec\beta\preceq\alpha }a_{\beta })z_{\alpha }^r$$
\noindent is an $(\varepsilon,j_{t+1}+\ldots +j_r)$- R.I.s.c.c.

For each $r=1,\ldots ,n$, we set $x^r=\frac{y_0^r}{\|y_0^r\|}$. If $r\geq 2$
then  for
each $\alpha\in {\cal T}$, $1\leq |\alpha |\leq r-1$, we set
$x_{\alpha }^r=\frac{1}{\|y_0^r\|}y_{\alpha }^r$, so that, for each
$t\leq r-1$,
$$x^r=\sum_{\alpha\in {\cal T}, |\alpha |=t}(\prod_{0\prec\beta\preceq
\alpha }a_{\beta })x_{\alpha }^r.$$

\medskip

\noindent {\bf 1.19 Lemma.} {\it For each $r\leq n$, $t<r$, and $\alpha\in
{\cal T}$ with $|\alpha |=t$,}
$$\frac{1}{16}\leq\|x_{\alpha }^r\|\leq 16(1+\eta ).$$

\medskip

\noindent {\it Proof:} By the construction, for each $t\leq r-1$
and $\alpha\in {\cal T}$ with $|\alpha |=t$, $y_{\alpha }^r$ is an
$(\varepsilon, j_{t+1}+\ldots +j_r)$- R.I.s.c.c. It follows from
Proposition 1.15 that
$$\frac{\theta_{j_{t+1}+\ldots +j_r+1}}{2}\leq\|y_{\alpha }^r\|
\leq 8\theta_{j_{t+1}+\ldots +j_r}.$$
\noindent Hence, for $0<|\alpha |=t$,
$$\frac{1}{16}\leq\frac{1}{16}\frac{\theta_{j_{t+1}+\ldots +j_r+1}}
{\theta_{j_1+\ldots +j_r}}\leq\|x_{\alpha }^r\|=
\frac{\|y_{\alpha }^r\|}{\|y_0^r\|}\leq
\frac{16\theta_{j_{t+1}+\ldots +j_r}}{\theta_{j_1+\ldots +j_r+1}}\leq
16(1+\eta ).\;\;\;\Box $$

\bigskip

\noindent {\bf 1.20 Lemma.} {\it Let $r\geq 2$ and
 $\alpha\in {\cal T}$ with $|\alpha |=t<r-1$.
If $i<j_{t+1}+\ldots +j_{r-1}$ and $(E_p)_{p=1}^k$ is an
${\cal F}_i$-admissible family of sets, then}
$$\sum_{p=1}^k\|E_px_{\alpha }^r\|\leq 32(1+\eta).$$

\medskip

\noindent {\it Proof:} By the construction,
$$y_{\alpha }^r=\sum_{|\gamma |=r-1,\alpha\prec\gamma }(\prod_{\alpha\prec\beta
\preceq\gamma }a_{\beta })y_{\gamma }^r,$$
\noindent where $l_{\gamma }^{r-1}<y_{\gamma }^r<l_{\gamma^+}^{r-1}$ for
every $\gamma\in {\cal T}$ with $|\gamma |=r-1$ and $\alpha\prec\gamma $.
(Recall that $y_{\gamma }^r$ is a convex combination of
$(z_{\beta }^r)_{|\beta |=r}$ and that $\max {\rm supp}(z_{\beta }^r)=
l_{\beta }^r$. By the definition of $({\cal T}_1^X,\ldots , {\cal T}_n^X)$,
we have $l_{\gamma }^{r-1}<l_{\beta }^r<l_{\gamma^+}^{r-1}$.)

Also, the corresponding basic convex combination
$$u_{\alpha }^{r-1}=\sum_{|\gamma |=r-1,\;\alpha\prec\gamma }
(\prod_{\alpha\prec\beta\preceq\gamma }a_{\beta })e_{l^{r-1}_{\gamma }}$$
\noindent is an $(\varepsilon,j_{t+1}+\ldots +j_{r-1})$--basic s.c.c.

An argument similar to the one in Lemma 1.13 yields
$$\sum_{p=1}^k\|E_py_{\alpha }^r\|\leq 2\max_{|\gamma |=r-1,\;\alpha\prec\gamma }
\|y_{\gamma }^r\|.$$
\noindent Dividing by $\|y_0^r\|$ we obtain the conclusion.   $\Box $

\bigskip

\noindent {\bf 1.21 Proposition.}
 {\it The sequence $\{x^r\}_{r=1}^n$ is $36$-equivalent
to the standard basis of $\ell_{\infty }^n$.}

\medskip

\noindent {\it Proof:} We need to prove that
$$\|\sum_{r=1}^nx^r\|\leq 36.$$
\noindent To do this we estimate $\varphi (\sum_{r=1}^nx^r)$ for
$\varphi\in K$, distinguishing two cases for $\varphi $:

\medskip

{\bf Case I:} $\varphi\in {\cal A}_i,\;i\geq j_0$.

\medskip

\noindent Let $r_0\in \{ 0,\ldots ,n\}$ be such that
$$j_{r_0}\leq i<j_{r_0+1}.$$
\noindent Then,

\medskip

(a) For $r\geq r_0+2$ we get $i<j_{r-1}<j_1+\ldots +j_{r-1}$. Using Lemma 1.20,
we see that
$$|\varphi (x^r)|\leq 32\theta_i(1+\eta) \leq 64\theta_{j_0}<\eta .$$

\medskip

(b) Let now $1\leq r\leq r_0-1$. We know that $y_0^r$ is an
$(\varepsilon,j_1+j_2+\ldots +j_r)$- R.I.s.c.c. of the $z_i$'s. Also,
$\varphi\in {\cal A}_i$, where $j_1+j_2+\ldots +j_r<j_{r+1}\leq i$.

Let $z_{i_1},\ldots ,z_{i_k}$ be the semi-normalized s.c.c.'s which
compose $y_0^r$ where, for $p=1,\ldots ,k$, $z_{i_p}$ is a
$(\theta^2_{t_p^r},t_p^r)$-seminormalized s.c.c. Set $t_0^r=t_{i_0}$
where by construction $t_{i_0}$ is such that
$\frac{\theta_{t_{i_0}}}{\theta_{j_1+\ldots +j_n+1}}<\frac{\eta }{16}$
and $t_{i_0}=t_0^r<t_p^r$ for all $p=1,\ldots ,k$.

\medskip

From Proposition 1.15 we get
$$|\varphi (y_0^r)|\leq 4\theta_i\leq 4\theta_{j_{r+1}}\;\;\;
{\rm if}\; i<t_1^r,$$ 
\noindent and
$$|\varphi (y^r_0)|\leq 4(\theta_{t^r_0}+\varepsilon)\;\;\;
{\rm if}\; i\geq t_1^r.$$ 
\noindent Dividing by $\|y_0^r\|$ and by the choice of
the $j_k$'s we obtain
$$|\varphi (x^r)|\leq\frac{8\theta_{j_{r+1}}}{\theta_{j_1+\ldots +j_r+1}}
<\eta \;\;\;{\rm if}\;i<t_1^r,$$
\noindent and
$$|\varphi (x^r)|\leq 8\frac{\theta_{t_{i_0}}}{\theta_{j_1+\ldots +j_r+1}}+
8\frac{\theta_{j_1+\ldots +j_r+1}^2}{\theta_{j_1+\ldots +j_r+1}}
<\frac{\eta }{2}+8\theta_{j_1+\ldots +j_r+1}<\eta \;\;\;{\rm if}\;i\geq t_1^r.$$

We conclude that, in this case,

$$|\varphi (\sum_{r=1}^n x^r)| \leq |\varphi (\sum _{r\neq r_0,r_{0+1}}x^r)|
+|\varphi (x^{r_0})|+|\varphi (x^{r_0+1})| \leq n\eta +2 <3.$$

\medskip

{\bf Case II:} $\varphi\in {\cal A}_i,\;i<j_0$. 

\medskip

Consider an analysis $\{K^s(\varphi )\}_{s=1}^q$ of $\varphi $.
For $s\leq q$ and $f\in K^s(\varphi )$, let $f^{+ }\in K^s(\varphi )$
be the successor of $f$ in $K^s(\varphi )$; that is, $f^{+ }$ is
such that ${\rm supp}f<{\rm supp}f^{+ }$ and if $g\in K^s(\varphi )$
with ${\rm supp}f<{\rm supp}g$ then either $g=f^{+}$ or
${\rm supp}f^{+ }<{\rm supp}g$.

For $f\in \bigcup_sK^s(\varphi )$, we set
$$E^f=[\min ({\rm supp}f),
\min ({\rm supp}f^{+ }))\subset {\bf N}$$
\noindent ($E^f=[\min ({\rm supp}f),\max ({\rm supp}x^n)]$ if $f$ does not
have a successor).

Recall that $x^1=\sum_{k=1}^ma_kz_k^1$ and, for $k=1,\ldots ,m$,
$l_k^1=\max ({\rm supp}z_k^1)$. We set 

\centerline {$I_k=[l_k^1,l_{k+1}^1)\subset
{\bf N}$, $k=1,\ldots ,m-1$ and $I_m=[l_m^1,\max ({\rm supp}x^n)]$.}

Notice that for $r\geq 2$ we have ${\rm supp}(x_k^r)\subset I_k$.

For $k=1,\ldots ,m$ and $f\in\bigcup_sK^s(\varphi )$, we say
that $f$ {\it covers} $I_k$ if $I_k\subset E^f$.

We may assume without loss of generality that $\min ({\rm supp}\varphi )
\leq l_1^1$. Therefore, for fixed $s$, any $I_k$ is either
covered by some $f$ in $K^s(\varphi )$ or intersected by $E^f$ for at least
two different $f$'s in $K^s(\varphi )$. Also, every $I_k$ is covered by
$\varphi .$

\medskip

Set now
$$J_1=\{k=1,\ldots ,m:I_k\;{\rm is}\;{\rm covered}
\;{\rm by}\;{\rm some}\;{\rm functional}$$
$$\;\;\;\;\;\;\;\;\;\;\;\;\;\;{\rm in}\;\cup K^s(\varphi )\;{\rm belonging}\;{\rm to}\;
{\rm some}\;{\rm class}\;{\cal A}_l\;{\rm with}\;l\geq j_0\},$$
\noindent and
$$J_2=\{k=1,\ldots ,m:I_k\;{\rm is}\;{\rm covered}\;{\rm only}\;{\rm by}\;
{\rm functionals}$$
$$\;\;\;\;\;\;\;\;\;\;\;\;\;\;{\rm in}\;\cup K^s(\varphi )\;{\rm which}\;{\rm belong}\;{\rm to}
\;\cup_{l<j_0}{\cal A}_l\}.$$

Consider any $k\in J_1$. Let $f\in \cup K^s(\varphi )$ be a functional
which covers $I_k$ and such that $f\in {\cal A}_l$ for some $l\geq j_0$.
Then, exactly as in Case I we can get
$$|\varphi (x_k^r)|\leq |f(x_k^r)|<\eta$$
\noindent for all but two $r\in\{2,\ldots ,n\}$. This gives
$|\varphi (\sum_{r=2}^nx_k^r)|\leq n\eta +32(1+\eta)<34$, and we
conclude that
$$|\varphi (\sum_{r=2}^n\sum_{k\in J_1}a_kx_k^r)|<34.$$

\bigskip

We turn now to $J_2$. Let $\varphi =\theta_i\sum_{p=1}^sf_p$
where $i<j_0$. Consider the set
$$R_1=\{k\in J_2: I_k\;{\rm is}\;{\rm intersected}\;{\rm by}
\;{\rm at}\;{\rm least}\;{\rm two}\;f_p{\rm s}\}.$$
\noindent Since the family $(f_p)_{p=1}^s$ is ${\cal F}_i$-admissible,
the set $\{l^1_k:k\in R_1\backslash\min R_1\}$ belongs to
${\cal F}_i\subset {\cal F}_{j_0}$ and so,
$\{l_k^1:k\in R_1\}\in {\cal F}_{j_1-1}$. Therefore, $\sum_{k\in R_1}a_k<
\varepsilon $.

\medskip

Let $L_1=J_2\backslash R_1$ and, for $p=1,\ldots ,s$, let
$$L_1^p=\{k\in L_1:I_k\subset E^{f_p}\}.$$ 
For any $r\geq 2$, we get
$$|\varphi (\sum_{k\in J_2}a_kx_k^r)|\leq
\theta_i\left (\sum_{p=1}^s|f_p(\sum_{k\in L_1^p}a_kx_k^r)|\right )+
(\sum_{k\in R_1}a_k)\max_k\|x_k^r\|$$
$$\leq\theta_{1}\left (\sum_{p=1}^s|f_p(\sum_{k\in L_1^p}a_kx_k^r)|\right )+
\varepsilon\max_k\|x_k^r\|.$$

\noindent Consider now any $p$, $1\leq p\leq s$, with $L_1^p\neq\emptyset $.
By the definition of $J_2$ this implies that
$f_p=\theta_{i_p}\sum_{t=1}^{l_p}g^p_t$ where $i_p<j_0$ and
$(g^p_t)_{t=1}^{l_p}$ is ${\cal F}_{i_p}$-admissible. (It is clear that we 
cannot have $f_p\in K^0$ and $L^p_1\neq\emptyset.$)
We will partition $L_1^p$ in the same way that we partitioned $J_2$:
We set
$$R_2^p=\{k\in L_1^p: I_k\;{\rm is}\;{\rm intersected}\;{\rm by}
\;{\rm at}\;{\rm least}\;{\rm two}\;g^p_t\;{\rm s}\}$$
\noindent and for each $t=1,\ldots ,l_p$,
$$L_2^t(p)=\{k\in L_1^p: I_k\subset E^{g^p_t}\}.$$
The family $\{ g^p_t:  p$ such that $L^p_1\neq\emptyset,\; t=1,\ldots l_p\}$
 is ${\cal F}_{i+j_0}$-admissible
and so the set 
 $\{l_k^1:k\in\cup_{p=1}^sR_2^p\}$ belongs to
${\cal F}_{i+j_0+1}
\subset{\cal F}_{2j_0}\subset {\cal F}_{j_1-1}$. We conclude that
$$\sum_{k\in\cup_p R_2^p}a_k<\varepsilon .$$
\noindent So, for each $r\geq 2$ we get the estimate
$$|\varphi (\sum_{k\in J_2}a_kx_k^r)|\leq\theta_{1}\sum_p\theta_{i_p}\sum_t
|g^p_t(\sum_{k\in L_2^t(p)}a_kx_k^r)|+\theta_{1}\sum_pf_p
(\sum_{k\in R_2^p}a_kx_k^r)+\varepsilon\max_k\|x_k^r\|$$
$$\leq\theta_{1}^2\sum_{p,t}|g^p_t(\sum_{k\in L_2^t(p)}a_kx_k^r)|+\theta_{1}
(\sum_{k\in\cup_pR_2^p}a_k)\max_k\|x_k^r\|+\varepsilon\max_k\|x_k^r\|$$
$$\leq\theta_{1}^2\sum_{p,t}|g^p_t(\sum_{k\in L_2^t(p)}a_kx_k^r)|+(\theta_{1}+1)
\varepsilon .$$

\noindent 
We can now partition each $L_2^t(p)$ and continue in this manner for $s_0$
steps, where $\theta_{1}^{s_0}<\eta$. By the choice of $j_1$,
$j_0s_0<j_1$. Recall that $\varphi\in K^q$. If $q>s_0$ then for $r\geq 2$,
$$|\varphi (\sum_{k\in J_2}a_kx_k^r)|\leq
\theta_{1}^{s_0}\sum_{f\in K^{q-s_0}(\varphi )}f(\sum_{I_k\subset E^f}a_kx_k^r)
+(1+\theta_{1}+\ldots +\theta_{1}^{s_0-1})\varepsilon\max_k\|x_k^r\|$$
Of course, if $q\leq s_0$ then we have only the second term at the right hand
side. Finally, for $r\geq 2$, we get
$$|\varphi (\sum_{k\in J_2}a_kx_k^r)|<
\max_k\|x_k^r\|(\eta +\frac{\varepsilon}{1-\theta_1} )<60\eta .$$
\noindent We conclude that
$$|\varphi (\sum_{r=1}^nx^r)|\leq |\varphi (x^1)|+
|\varphi (\sum_{r=2}^n\sum_{k\in J_1}a_kx_k^r)|+
\sum_{r=2}^n|\varphi (\sum_{k\in J_2}a_kx_k^r)|$$
$$\leq 1+34+60n\eta <36.$$
\noindent This completes the proof of the Proposition. Theorem 1.6
now follows.   $\Box $ 

\bigskip

\noindent {\bf C. Modified Mixed Tsirelson spaces.}

\medskip

The modified Tsirelson space $T_M$ was introduced by W.B. Johnson in [10].
Later, P. Casazza and E. Odell ([6]) proved that $T_M$ is naturally 
isomorphic to $T$. Analogously, given a sequence of compact families
$\{ {\cal M}_k\}_{k=1}^{\infty}$ in $[{\bf N}]^{<\omega}$ and a sequence
of positive reals $\{\theta_k\}_{k=1}^{\infty}$, we define the {\it modified
mixed Tsirelson space} 
$T_M[({\cal M}_k, \theta_k)_{k=1}^{\infty}].$

\medskip

\noindent {\bf 1.22 Definition.} Let ${\cal M}$ be a  family of finite
subsets of ${\bf N}$.

(a) A finite sequence $(E_i)_{i=1}^k$ of finite non-empty subsets
of ${\bf N}$ is said to be ${\cal M}$-{\it allowable} if the set
$\{\min E_1,\min E_2,\ldots ,\min E_k\}$ belongs to ${\cal M}$
and $E_i\cap E_j=\emptyset $ for all $i,j=1,\ldots ,k,\;i\neq j$.

(b) A finite sequence $(x_i)_{i=1}^k$ of vectors in $c_{00}$ is
${\cal M}$-{\it allowable} if the sequence
$({\rm supp}(x_i))_{i=1}^k$ is ${\cal M}$-allowable.

\medskip

\noindent {\bf 1.23 Definition of the space
 $T_M[({\cal M}_k,\theta_k)_{k=1}^{\infty }]$}. 
Let $({\cal M}_k)_k$ be a sequence of compact, hereditary
and spreading families of finite
subsets of ${\bf N}$ and let $(\theta_k)_k$ be a sequence of positive reals
with $\theta_k<1$ for every $k$ and $\lim_k\theta_k=0$.
Inductively, we define a subset $K$ of $B_{\ell_{\infty }}$ as follows:

\noindent We set $K^0=\{\pm e_n:n\in {\bf N}\}$.

\noindent For $s\geq 0$, given $K^s$ we define for each $k\geq 1$,
$$K_k^{s+1}=\{\theta_k (\sum_{i=1}^nf_i):n\in {\bf N},f_i\in K^s,i\leq n,\;
{\rm and}\;{\rm the}\;{\rm sequence}\;(f_i)_{i=1}^n\;{\rm is}\;{\cal M}_k-{\rm allowable}\}.$$
\noindent We set
$$K^{s+1}=K^s\bigcup\left(\bigcup_{k=1}^{\infty }K_k^{s+1}\right).$$
\noindent Finally, we define
$$K=\bigcup_{s=0}^{\infty }K^s.$$
Note that $K$ is the smallest subset of $B_{\ell_{\infty }}$ which
contains $\pm e_n$ for all $n\in {\bf N}$ and has the property
that $\theta_k(f_1+\ldots +f_n)$ is in $K$ whenever $f_1,\ldots ,f_n\in K$
and the sequence $(f_i)_{i=1}^n$ is ${\cal M}_k$-allowable.

\medskip

\noindent We now define a norm on $c_{00}$ by
$$\|x\|=\sup_{f\in K}\langle x,f\rangle \;\;\;{\rm for}\;{\rm all}\;x\in c_{00}.$$

The space $T_M[({\cal M}_k,\theta_k)_{k=1}^{\infty }]$ is the completion of
$(c_{00},\|.\|)$. We call $K$ the {\it norming set} of\\
$T_M[({\cal M}_k,\theta_k)_{k=1}^{\infty }]$.

\bigskip

The following Proposition is an easy consequence of the definition:

\medskip

\noindent {\bf 1.24 Proposition.} {\it Let $X=T_M[({\cal M}_k,\theta_k)_{k=1}^{\infty }]$.

{\rm (a)} The norm of $X$ satisfies the following implicit equation: 
For all $x\in X$,
$$\|x\|=\max\{\|x\|_{\infty },\sup_k\theta_k\sup \{\sum_{i=1}^n\|E_ix\|
: (E_i)_{i=1}^n \; {\rm is }\; {\cal M}_k -{\rm allowable}\}\}.$$

{\rm (b)} The sequence $(e_n)_{n=1}^{\infty }$ is a 1-unconditional basis
for $X$.   $\Box $}

\bigskip

We also consider {\it boundedly modified mixed Tsirelson spaces} denoted by
$$T_{M(m)}[({\cal M}_k, \theta_k)_{k=1}^{\infty}],$$
\noindent for some $m\in {\bf N}$.
The definition of 
$T_{M(m)}[({\cal M}_k, \theta_k)_{k=1}^{\infty}] $ is similar to that of
$T_M[({\cal M}_k, \theta_k)_{k=1}^{\infty}],$ the only difference being
that at the inductive step $s+1$ we set
$$K_k^{s+1}=\{\theta_k (\sum_{i=1}^nf_i):n\in {\bf N},f_i\in K^s,i\leq n,\;
{\rm and}\;{\rm the}\;{\rm sequence}\;(f_i)_{i=1}^n\;{\rm is}\;{\cal M}_k-{\rm allowable}\},$$
\noindent for $k\leq m$, while
$$K_k^{s+1}=\{\theta_k (\sum_{i=1}^nf_i):n\in {\bf N},f_i\in K^s,i\leq n,\;
{\rm and}\;{\rm the}\;{\rm sequence}\;(f_i)_{i=1}^n\;{\rm is}\;{\cal M}_k-
{\rm admissible}\},$$
\noindent for $k\geq m+1.$

\medskip

\noindent {\bf 1.25 Proposition.} 
{\it Let $Y=
T_{M(m)}[({\cal M}_k, \theta_k)_{k=1}^{\infty}]. $

{\rm (a)} The norm $\| \cdot \|$ of $Y$ satisfies the following implicit
equation:
\begin{eqnarray*}
\|x\|=\max\Big
\{\|x\|_{\infty },& &\max_{k\leq m}\theta_k\sup \{\sum_{i=1}^n\|E_ix\|
: (E_i)_{i=1}^n \; {\rm is }\; {\cal M}_k -{\rm allowable}\},\\
& &\sup_{k\geq m+1} \theta_k \sup\{ \sum_{i=1}^n\| E_ix\|: (E_i)_{i=1}^n\;
{\rm is }\; {\cal M}_k -{\rm admissible}\} \Big\}.\end{eqnarray*}

 {\rm (b)} The sequence $(e_n)_n$
is a 1-unconditional basis for $Y.$   $\Box $}

\bigskip

In the sequel we consider spaces 
$T_{M}[({\cal M}_k, \theta_k)_{k=1}^{\infty}] $ or 
$T_{M(m)}[({\cal M}_k, \theta_k)_{k=1}^{\infty}] $ where $({\cal M}_k)_k$
is a subsequence of the Schreier sequence $({\cal F}_n)_{n=1}^{\infty}.$
In this case, by Proposition 1.24(a) (resp. Proposition 1.25(a)) we
have that for all sequences $(x_i)_{i=1}^n$ of disjointly supported vectors
with ${\rm supp}x_i\subset [n,\infty ),$ 
 $$ \| \sum _{i=1}^nx_i\|\geq \theta_1\sum_{i=1}^n\|x_i\|$$
in 
$T_{M}[({\cal M}_k, \theta_k)_{k=1}^{\infty}] $  (resp. 
$T_{M(m)}[({\cal M}_k, \theta_k)_{k=1}^{\infty}]. $)  It is clear from this
inequality that $c_0$ is not finitely disjointly representable in
any block subspace of
$T_{M}[({\cal M}_k, \theta_k)_{k=1}^{\infty}] $  or 
$T_{M(m)}[({\cal M}_k, \theta_k)_{k=1}^{\infty}]. $   
Combining this with Theorem 1.6 we get the following.

\medskip

\noindent {\bf 1.26 Corollary.} {\it Let $(\theta_n)_{n=1}^{\infty}$ be a
regular sequence with $\lim \theta^{1/n}_n=1.$ Let 
$X=T_{M}[({\cal F}_k, \theta_k)_{k=1}^{\infty}] $  or 
$X=T_{M(m)}[({\cal F}_k, \theta_k)_{k=1}^{\infty}] .$   Then the spaces
$X$ and
$T[({\cal F}_k, \theta_k)_{k=1}^{\infty}] $  are totally
incomparable. $\Box $   }

\bigskip

\noindent {\bf 1.27 Theorem.} {\it Suppose that the
sequence $(\theta_k)_k$ decreases to 0 and that the Schreier family ${\cal S}
$ is contained in  ${\cal M}_1 $.
Then,  the spaces
 $T_M[({\cal M}_k,\theta_k)_{k=1}^{\infty }]$ and
 $T_{M(m)}[({\cal M}_k,\theta_k)_{k=1}^{\infty }]$, $m=1,2,\ldots$
 are reflexive.}

\medskip

\noindent {\it Proof:} Let $X=T_M[({\cal M}_k,\theta_k)_{k=1}^{\infty }]$.
The proof for $T_{M(m)}[({\cal M}_k,\theta_k)_{k=1}^{\infty }]$ is the
same. We shall prove that the basis $(e_n)_{n=1}^{\infty }$
is boundedly complete and shrinking in $X$.

\medskip

\noindent (a) {\it $(e_n)_{n=1}^{\infty }$ is boundedly complete:} 
Suppose on the contrary that there exist $\varepsilon >0$
and a block sequence $\{x_i\}_{i=1}^{\infty }$ of $\{e_n\}_{n=1}^{\infty }$ such that
$\sup_n\|\sum_{i=1}^nx_i\|\leq 1$ while $\|x_i\|\geq \varepsilon $ for $i=1,2,\ldots $.

Choose $n_0\in {\bf N}$ such that $n_0\theta_1>\varepsilon $. Then, the
finite sequence $(x_i)_{i=n_0+1}^{2n_0}$ is ${\cal S}$-allowable and since
${\cal S}\subseteq {\cal M}_1$ it is ${\cal M}_1$-allowable. Using Proposition 
1.24(a) (resp. 1.25(a))
we get
$$\|\sum_{i=n_0+1}^{2n_0}x_i\|\geq\theta_1\sum_{i=n_0+1}^{2n_0}\|x_i\|\geq n_0
\theta_1\varepsilon >1,$$
\noindent a contradiction which completes the proof.

\medskip

\noindent (b) {\it $(e_n)_{n=1}^{\infty }$ is a shrinking basis:} For
$f\in X^{\ast }$, $m\in {\bf N}$, we denote by $Q_m(f)$ the restriction
of $f$ to the space spanned by $(e_k)_{k\geq m}$. We need to prove that,
for every $f\in B_{X^{\ast }}$,
$Q_m(f)\rightarrow 0\;\;{\rm as}\;m\rightarrow\infty .$

Let $K$ be the norming set of $X$.
 Then $B_{X^{\ast }}=\overline{{\rm co}(K)}$
 where the closure is in the topology of  pointwise
convergence.  We shall show that
for all $f\in B_{X^{\ast }}$ there is $l\in {\bf N}$ such that
$Q_l(f)\in\theta_1B_{X^{\ast }}$. By standard arguments it
suffices to prove this for $f\in\overline{K}$.

Let $f\in\overline{K}$. Let $(f^n)_{n=1}^{\infty }$ be a sequence
in $K$ converging pointwise to $f$. If $f^n\in K^0$ for an infinite
number of $n$, then there is nothing to prove. So, suppose that for
every $n$ there are $k_n\in {\bf N}$, a set
$M_n=\{m_1^n,\ldots ,m_{d_n}^n\}\in {\cal M}_{k_n}$ and vectors
$f_i^n\in K$, $i=1,\ldots ,d_n$ such that $f^n=\theta_{k_n}\sum_{i=1}^{d_n}
f_i^n$, $m_i^n=\min {\rm supp}(f_i^n)$, $i=1,\ldots ,d_n$ and
${\rm supp}(f_i^n)\cap {\rm supp}(f_j^n)=\emptyset $. If there is
a subsequence of $(\theta_{k_n})_n$ converging to 0, then $f=0$.
So we may assume that there is a $k$ such that $k_n=k$ for
all $n$, that is, $\theta_{k_n}=\theta_k$ and
$M_n=\{m_1^n,\ldots ,m_{d_n}^n\}\in {\cal M}_k$.

Since ${\cal M}_k$ is compact, substituting $\{f^n\}$ with a
subsequence we get that there is a set $M=\{m_1,\ldots ,m_d\}\in {\cal M}_k$
such that the sequence of indicator functions of $M_n$ converges
to the indicator function of $M$. So, for large $n$, $m_i^n=m_i$,
$i=1,2,\ldots ,d$ and $m_{d+1}^n\rightarrow\infty $ as $n\rightarrow\infty $.
Since $\min {\rm supp}f_{d+1}^n=m_{d+1}^n\rightarrow\infty $, the
sequence ${\tilde{f}}^n=\theta_k\sum_{i=1}^df_i^n$ tends to $f$
pointwise and we may assume that $f^n=\theta_k\sum_{i=1}^df_i^n$.
Passing again to a subsequence of $\{f^n\}$ we have that,
for each $i=1,\ldots ,d$ there exists $f_i\in\overline{K}$ with
$f_i^n\rightarrow f_i$ pointwise and $f=\theta_k(f_1+\ldots +f_d)$.

Now, for each $i=1,\ldots,d$, either $f_i^n=e^*_{m_i}$ for all $n$ (eventually)
or 
$$f_i^n=\theta_{k_i^n}\sum_{m=1}^{l_i^n}g_{m}^{n,i}\;,\;\;i=1,\ldots ,d$$
\noindent where for every $n\in {\bf N}$ and $m=1,\ldots ,l_i$, $g_m^{n,i}\in K$
and the family $\{g_m^{n,i}\}_{m=1}^{l_i^n}$ is
${\cal M}_{k_i^n}$-allowable. Let $A\subset \{ 1,\ldots, d\}$ be the set
of indices $i$ for which $f^n_{i}$ is of the second type for all $n$.
As before, forgetting those $i$'s for which $f_i^n\rightarrow 0$, we may
assume that, for each $i\in A$, there 
is $k_i$ such that $k_i^n=k_i$ and a set $M_i=\{m_1^i,\ldots ,m_{l_i}^i\}$
such that $m_r^i=\min{\rm supp}(g_r^{n,i})$ for all
$n=1,2,\ldots $, $r=1,\ldots ,l_i$, and $\min{\rm supp}(g_{l_i+1}^{n,i})\rightarrow\infty $ as $n\rightarrow\infty $. So, for $i\in A$, the
sequence $\tilde{f_i^n}=\theta_{k_i}\sum_{m=1}^{l_i}g_m^{n,i}$
tends to $f_i$ pointwise.

Let $l=\max (\{\sum_{i\in A}l_i\}\cup \{ m_i:i=1\ldots d\})$
 and  $h_m^{n,i}=Q_l(g_m^{n,i})\in K$, $i\in A, m=1,\ldots l_i, n=1,2,\ldots$.
Then, the sequence 
$\theta_k\sum_{i\in A}\theta_{k_i}\sum_{m=1}^{l_i}h_m^{n,i}=Q_l(\theta_k\sum
_{i=1}^d\tilde{f_i^n})$
tends to $Q_l(f)$ as $n\rightarrow\infty $.

On the other hand, since, for each $n$, $\#\{h_m^{n,i},\;i\in A,\;
m=1,\ldots ,l_i\}\leq l$, $l\leq\min{\rm supp}(h_m^{n,i})$
for every $i$ and $m$, and the sets ${\rm supp}(h_m^{n,i})$, $i\in A$,
$m=1,\ldots ,l_i$
are mutually disjoint, we get that the family $\{h_m^{n,i}\}_{i,m}$
is Schreier-allowable. Since the Schreier family ${\cal S}$ is contained
in ${\cal M}_{1}$, $0<\frac{\theta_{k_i}}{\theta_1}\leq 1$, $\{ h_m^{n,i}\}_
{i,m}$ is ${\cal S}$-allowable for every $n$ 
and $h_m^{n,i}\in K$, it is easy to see that
$\frac{1}{\theta_k} Q_l(\theta _k\sum_{i=1}^d \tilde{f_i^n})  =
\theta_1(\sum _{i\in A}\frac{\theta_{k_i}}{\theta_1}\sum_{m=1}^{l_i}h_m^{n,i})
\in {\rm co}(K)$ for all $n$. We conclude  that 
$ Q_l(\theta _k\sum_{i=1}^d \tilde{f_i^n}) \in \theta_k {\rm co}(K)$, and so,
$Q_l(f)\in \theta_k\overline{{\rm co}(K)}\subseteq\theta_1\overline{{\rm co}
(K)}.$  $\Box$

\bigskip

We note that the 2-convexifications  
$T^{(2)}_{M}[({\cal F}_k, \theta_k)_{k=1}^{\infty}] $   and
$T^{(2)}_{M(m)}[({\cal F}_k, \theta_k)_{k=1}^{\infty}] $  of  
$T_{M}[({\cal F}_k, \theta_k)_{k=1}^{\infty}] $ and 
$T_{M(m)}[({\cal F}_k, \theta_k)_{k=1}^{\infty}] $  are weak Hilbert spaces.
The proof of this is similar to the proof of the analogous statement for
the 2-convexifications 
$T^{(2)}_{\delta}$ of 
the Tsirelson spaces $T_{\delta}$ 
as presented in [15] (Lemma 13.5). It is an immediate consequence 
of Theorem 1.27 that 
$T^{(2)}_{M}[({\cal F}_k, \theta_k)_{k=1}^{\infty}] $   (and
$T^{(2)}_{M(m)}[({\cal F}_k, \theta_k)_{k=1}^{\infty}] $)   does not
contain $\ell_2.$ Moreover, we can show that for sequences $(\theta_n)_n$
with $\lim_n\theta^{1/n}_n=1$, no subspace of 
$T^{(2)}_{M}[({\cal F}_k, \theta_k)_{k=1}^{\infty}] $   
(or $T^{(2)}_{M(m)}[({\cal F}_k, \theta_k)_{k=1}^{\infty}] $)   can
be isomorphic to a subspace of $T_{\delta}^{(2)}$.  It suffices to prove the
following.

\medskip

\noindent {\bf 1.28 Proposition.} {\it Let $0<\delta<1$ and let
 $(\theta_n)_n$ be a regular sequence
with 
$\lim \theta^{1/n}_n=1.$  Let
$X=T_{M}[({\cal F}_k, \theta_k)_{k=1}^{\infty}] $  or 
$X=T_{M(m)}[({\cal F}_k, \theta_k)_{k=1}^{\infty}] .$   
Then the spaces $X$ and  $T_{\delta}$ are totally incomparable.}

\medskip

\noindent {\it Proof:} Let $X=
T_{M}[({\cal F}_k, \theta_k)_{k=1}^{\infty}] $   
$X=T_{M(m)}[({\cal F}_k, \theta_k)_{k=1}^{\infty}]. $   Suppose on the 
contrary that there exist normalized block sequences $\{ x_i\}_i $ in
$X$ and $\{y_i\}_i$ in $T_{\delta}$ which are equivalent as basic sequences.
Let $l_i=\min {\rm supp}y_i,$ $i=1,2,\ldots.$ From [5] Theorem 13
we get that $\{ x_i\}_X$ is equivalent to $\{ e_{l_i}\}_{T_{\delta}}.$
Let $m_i=\min {\rm supp}x_i,$ $i=1,2,\ldots$. We choose a 
subsequence $\{ i_k\}_k$ of indices such that either
$l_{i_1}\leq m_{i_1}<l_{i_2}\leq m_{i_2}<\ldots$
or
$m_{i_1}<l_{i_1}<m_{i_2}<l_{i_2}<\ldots$
In either case, using Theorem 13 [5] once more, we get that
the basic sequences $\{ e_{l_{i_k}}\}$ and $\{ e_{m_{i_k}}\}$ are
equivalent in $T_{\delta}.$ We conclude that $\{ e_{m_{i_k}}\}_{T_{\delta}}$
is equivalent to $\{ x_{i_k}\}_X.$

Let now $j\in N$ and let $\sum _{k\in A}a_ke_{m_{i_k}}$ be a $(\theta^2_j, j)$
- special convex combination. As in Lemma 1.12 we get that 
$\|\sum_{k\in A}a_ke_{m_{i_k}}\|_{T_{\delta}}\leq \delta^j +\theta^2_j.$
On the other hand, since the sequence $(x_{i_k})_{k\in A}$ is ${\cal F}_j$-
admissible, we have that $\| \sum_{k\in A}a_kx_{i_k}\|_{X}\geq \theta_j.$
But the assumption $\lim \theta^{1/j}_j=1$ implies that 
$\delta^{-j}\theta_j\rightarrow \infty.$ This leads to a contradiction
which completes the proof. $\Box$

\bigskip

\bigskip

\centerline { {\sc 2. The Space $X_{M(1),u}$}}

\medskip

We give an example of a boundedly modified mixed Tsirelson space
space of the form\\
$T_{M(1)}[({\cal F}_{k_j},\theta_j)_{j=1}^{\infty}]$
which is arbitrarily distortable.

\noindent {\bf Definition of $X_{M(1),u}.$}
We choose a sequence of integers $(m_j)_{j=1}^{\infty }$ such that
$m_1=2$ and for $j=2,3\ldots $, $m_j>m_{j-1}^{m_{j-1}}$. 

We choose inductively a subsequence 
$({\cal F}_{k_j})_{j=0}^{\infty }$ of $({\cal F}_n)_n$ :

We set $k_1=1$.
Suppose that $ k_j,\; j=1,\ldots ,n-1$ have been chosen.
Let $t_n$ be such that $2^{t_n}\geq m_n^2.$ We set
$ k_n=t_n(k_{n-1}+1)+1.$ 

\medskip

For $j=0,1,\ldots $, we set ${\cal M}_j={\cal F}_{k_j}.$
We define 
$$X_{M(1),u}=T_{M(1)}[({\cal M}_j,\frac{1}{m_j})_{j=1}^{\infty}]$$

\bigskip

\noindent {\bf Notation:} Let ${\cal F}$ be a family of finite
subsets of ${\bf N}.$  We set
$${\cal F}^{\prime }=\{A\cup B: A\in {\cal F}, B\in {\cal F},
A\cap B=\emptyset \}.$$

\medskip

\noindent {\bf 2.1 Definition.} Given $\varepsilon >0$ and
$j=2,3,\ldots $, an $(\varepsilon ,j)$-{\it basic special convex
combination ($(\varepsilon ,j)$- basic s.c.c.) relative
to $X_{M(1),u}$} is a vector of the form
$\sum_{k\in F}a_ke_k$ such that: $F\in {\cal M}_j,a_k\geq 0,
\sum_{k\in F}a_k=1$, $\{a_k\}_{k\in F}$ is decreasing,
and, for every $G\in {\cal F}^{\prime }_{t_j(k_{j-1}+1)}$, $\sum_{k\in G}a_k<
\varepsilon $.

\medskip

\noindent {\bf 2.2 Lemma.} {\it Let $j\geq 2,\varepsilon >0$, $D$ be an
infinite subset of ${\bf N}$. There exists an $(\varepsilon , j)$-basic
special convex combination relative to $X_{M(1),u},$ $x=\sum_{k\in F}a_ke_k$,
with $F={\rm supp}x\subset D.$}

\medskip

\noindent {\it Proof:} Since ${\cal M}_j={\cal F}_{t_j(k_{j-1}+1)+1},$
by Proposition 1.8 there exists a convex combination $x=\sum_{k\in F} a_ke_k$
with $F\in {\cal M}_j, $ $F\subset D$ and such that $\sum_{k\in G}a_k<
\frac{\varepsilon}{2}$ for all $G\in {\cal F}_{t_j(k_{j-1}+1)}.$
It is clear that this $x$ is an $(\varepsilon , j)$-basic s.c.c. relative
to $X_{M(1),u}.$  $\Box $

\medskip

In the sequel, when we refer to $(\varepsilon , j)$-special convex 
combinations
 we always imply ``relative to $X_{M(1),u}$".

\bigskip

\noindent {\bf Notation.} Let
$X'_{(n)}=T_{M(1)}[({\cal M}^{\prime}_l,\frac{1}{m_l})_{l=1}^{n}]$
and  let $K'{(n)}$ be the norming set of $X'_{(n)}.$ We denote by
$|\cdot|_{n}$ the norm of $X_{(n)}$ and by $|\cdot |^*_{n}$ the
corresponding dual norm.

 We set
$${\cal G}_{(n)}=\{{\rm supp}f:f\in K'{(n)}{\rm \; and \;}
{\rm for}\;{\rm every}\;k\in
{\rm supp}f,\;f(e_k)>\frac{1}{m_{n+1}^2}\}.$$

\noindent {\bf Remark.} Using Lemma 1.2 it is easy
to see that  ${\cal G}_{(n-1)}\subset  {\cal F}_{t_n(k_{n-1}+1)}.$ 
It follows that if $x=\sum_{k\in F}a_ke_k$ is an $(\varepsilon, n)$-
basic s.c.c. then, for all $G\in {\cal G}^{\prime}_{(n-1)},$
$\sum_{k\in G}a_k <\varepsilon.$

\bigskip

We give the definition of the set $K$ of functionals that define
the norm of the space $X_{M(1),u}$:

\medskip

We set $K_j^0=\{\pm e_n:n\in {\bf N}\}$ for $j=1,2,\ldots $.

\medskip

Assume that $\{K_j^n\}_{j=1}^{\infty } $ have been defined. Then, we set
$K^n=\bigcup_{j=1}^{\infty }K^n_j$, and for $j=2,3\ldots $ we set
$$K_j^{n+1}=K_j^n\cup\{\frac{1}{m_j}(f_1+\ldots +f_d):
{\rm supp}f_1<\ldots <{\rm supp}f_d, \; (f_i)_{i=1}^d
\;{\rm is}\;{\cal M}_j-{\rm admissible}$$
$${\rm and}\;f_1,\ldots ,f_d\;{\rm belong}\;{\rm to}\;K^n\},$$
\noindent while for $j=1$, we set
$$K_1^{n+1}=K_1^n\cup\{\frac{1}{2}(f_1+\ldots +f_d):f_i\in K^n,d\in {\bf N},
d\leq\min {\rm supp}f_1<\ldots <\min {\rm supp}f_d,$$
$${\rm and}\;{\rm supp}f_i\cap {\rm supp}f_j=\emptyset \}.$$

Set $K=\bigcup_{n=0}^{\infty }K^n$. Then, the norm $\|.\|$ of
$X_{M(1),u}$ is
$$\|x\|=\sup\{f(x):f\in K\}.$$ 

\noindent {\bf Notation.} For $j=1,2,\ldots $, we denote by ${\cal A}_j$ the
set ${\cal A}_j=\bigcup_{n=1}^{\infty }(K_j^n\backslash K^0)$. 
Then, $K=K^0\cup (\cup_{j=1}^{\infty }
{\cal A}_j)$. 

We will also consider the space
$T_{M(1)}[({\cal M}^{\prime}_{j},\frac{1}{m_j})_{j=1}^{\infty }]$.
We denote by $K'$ the norming set of this space and by $K^{\prime n}$,
$K^{\prime n}_j$, ${\cal A}'_j$ the subsets of $K'$ corresponding to
$K^n$, $K^n_j$ and ${\cal A}_j$ respectively.

\bigskip

\noindent {\bf 2.3 Definition.} A. Let $m\in {\bf N}$, $\varphi\in K^m\backslash
K^{m-1}$. An {\it analysis}
 of $\varphi $ is a sequence $\{K^s(\varphi )\}_{s=0}^m$
of subsets of $K$ such that:

(1) For every $s$, $K^s(\varphi )$ consists of  disjointly
supported elements of $K^s$, and
$\bigcup_{f\in K^s(\varphi )}{\rm supp}f={\rm supp}\varphi $.

(2) If $f$ belongs to $K^{s+1}(\varphi )$, then either $f\in K^s(\varphi )$
or there exists an ${\cal S}$-allowable family $(f_i)_{i=1}^d$ in 
$K^s(\varphi )$ such that $f=\frac{1}{2}(f_1+\cdots +f_d)$, or,
for some $j\geq 2$, there exists an ${\cal M}_j$-admissible family
$(f_i)_{i=1}^d$ in $K^s(\varphi )$ such that $f=\frac{1}{m_j}(f_1+\cdots
+f_d).$

(3) $K^m(\varphi )=\{\varphi \}$.

\medskip

\noindent B. For $g\in K^{s+1}(\varphi)\backslash K^0(\varphi )$, 
 the set of functionals $\{ f_1,\ldots , f_l\}
\subset K^s(\varphi )$ such that
 $g=\frac{1}{m_j}(\sum_{i=1}^lf_i)$  is called the {\it decomposition}
of $g.$

\bigskip

\noindent {\bf 2.4 Lemma.} {\it Let $j\geq 2$, 
$0<\varepsilon \leq\frac{1}{m_j^2}$, $M>0$
and let $x=\sum_{k=1}^mb_ke_{n_k}$ be 
an $(\varepsilon ,j)$-basic s.c.c. Suppose that the vectors  
 $x_k=\sum_{i=1}^{l_k}a_{i,k}e_{n_{i,k}}$ are such that $a_{i,k}\geq 0$
for all $i,k$ and
$\sum_{i=1}^{l_k}a_{i,k}\leq M$, $k=1,2,\ldots ,m$, and
$n_1\leq n_{1,1}<n_{2,1}<\ldots <n_{l_1,1}<n_2\leq
n_{1,2}<n_{2,2}<\ldots n_{3}\leq\ldots <n_{l_m,m}$. Then

\noindent {\rm (a)} For $\varphi\in\cup_{s=1}^{\infty}{\cal A'}_s$,

$|\varphi (\sum_{k=1}^mb_kx_k)|\leq\frac{M}{m_s}$,  if $\varphi\in {\cal A'}_s,
 s\geq j$

$|\varphi (\sum_{k=1}^mb_kx_k)|\leq\frac{2M}{m_sm_j}$,  if $\varphi\in {\cal A'}
_s, s<j$.

\noindent {\rm (b)} If $\varphi $ belongs to the norming set $K'(j-1)$ of
$T_{M(1)}[({\cal M}^{\prime}_{l},\frac{1}{m_l})_{l=1}^{j-1}]$, then}
$$|\varphi (\sum b_kx_k)|\leq\frac{2M}{m_j^2}.$$

\medskip

\noindent {\it Proof:} (1) If $s\geq j$, then the estimate is obvious.

Let $s<j$ and $\varphi =\frac{1}{m_s}\sum_{l =1}^df_{l }$.
Without loss of generality we assume that $\varphi(e_{n_{i,k}})\geq 0$
for all $n_{i,k}$. We set
$$D=\{n_{i,k} :\sum_{l =1}^df_{l }(e_{n_{i,k}})
>\frac{1}{m_j}\}.$$
\noindent We set $g_l=f_l|_D$. Then, $\frac{1}{m_s}\sum_{l=1}^dg_l\in
K'{(j-1)}$, and for every $k\in {\rm supp}(\frac{1}{m_s}\sum_{l=1}^dg_l)$
we have $\frac{1}{m_s}\sum_{l=1}^dg_l(e_k)>
\frac{1}{m_sm_j}>\frac{1}{m_j^2}$. Therefore,
$D={\rm supp}(\frac{1}{m_s}\sum_{l=1}^dg_l)\in {\cal G}_{(j-1)}$.
Let $B=\{ k: {\rm \; there \; exists \;} i \; {\rm with }\; n_{i,k}\in 
D\}.$ Then $B\in {\cal G}'_{(j-1)}$ and so, by the Remark after Lemma 2.2,
$\sum_{k\in B}b_k<\varepsilon\leq\frac{1}{m^2_j}.$ We get

$$\frac{1}{m_s}\sum_{l=1}^dg_l(\sum_{k=1}^m b_kx_k)
\leq\sum_{k\in B}b_k(\sum_{i=1}^{l_k}a_{i,k})\leq
M\sum_{k\in B}b_k\leq\frac{M}{m_j^2}.$$
\noindent  On the other hand,
$$(\frac{1}{m_s}\sum_{l =1}^df_{l |_{D^c}})(\sum b_kx_k)\leq\frac{M}{m_sm_j}.$$
\noindent Hence,
$$\varphi (\sum b_kx_k)\leq\frac{M}{m_sm_j}+\frac{M}{m_j^2}\leq\frac{2M}{m_sm_j}.$$

\medskip 

(b) We assume again that $\varphi$ is positive.
 We set $L=\{n_{i,k}:\varphi (e_{n_{i,k}})>\frac{1}{m_j^2}\}$. Then,
$$\varphi|_{L^c}(\sum b_kx_k)\leq\frac{M}{m_j^2}.$$
\noindent On the other hand, ${\rm supp}(\varphi|_L)\in {\cal G}_{(j-1)}$ and
as before we get
$\varphi|_L(\sum b_kx_k)\leq\frac{M}{m_j^2}$. 
 Hence,
$$|\varphi (\sum b_kx_k)|\leq\frac{2M}{m_j^2}.\;\;\;\Box $$

\bigskip

\noindent {\bf 2.5 Definition.} (a) Given a  block
sequence $(x_k)_{k\in {\bf N}}$ in $X_{M(1),u}$ and $j\geq 2$,
 a convex combination
$\sum_{i=1}^na_ix_{k_i}$ is said to be an
$(\varepsilon ,j)$-{\it special convex combination} of $(x_k)_{k\in {\bf N}}$
($(\varepsilon ,j)$-s.c.c), if there exist $l_1<l_2<\ldots <l_n$
such that $2<{\rm supp}x_{k_1}\leq l_1<{\rm supp}x_{k_2}\leq l_2<
\ldots <{\rm supp}x_{k_n}\leq l_n$, and $\sum_{i=1}^na_ie_{l_i}$
is an $(\varepsilon ,j)$-basic s.c.c.

(b) An $(\varepsilon ,j)$-s.c.c. $\sum_{i=1}^n a_ix_{k_i}$
 is called {\it seminormalized} if  $\| x_{k_i}\|=1,\; i=1,\ldots ,n$ and
$$\|\sum_{i=1}^na_ix_{k_i}\|\geq\frac{1}{2}.$$

\bigskip

\noindent {\bf 2.6 Lemma.} {\it Let $(x_k)_{k=1}^{\infty }$ be
a  block sequence in $X_{M(1),u}$ and $j=2,3\ldots $, $\varepsilon
>0$. Then, there exists a normalized finite block sequence
$\{y_k\}_{k=1}^n$ of $\{x_k\}_{k=1}^{\infty }$ and
 a convex combination $\sum_{k=1}^n a_ky_k$ which is
a  seminormalized $(\varepsilon ,j)$-s.c.c. }

\medskip

\noindent {\it Proof:} Using that ${\cal M}_j={\cal F}_{t_j(k_{j-1}+1)+1}$
where $2^{t_j}\geq m^2_j$,
the proof  
is similar to the proof of Lemma 1.11.   $\Box $

\bigskip

\noindent {\bf 2.7 Lemma.}  {\it Let $j\geq 3$ and
let  $x=\sum_{k=1}^n a_kx_k$ be
a $(\frac{1}{m_j^4},j)$-s.c.c  where $\|x_k\|\leq 1,$ $k=1,\ldots n.$ Suppose 
$\varphi =\frac{1}{m_r}\sum_{i=1}^df_i\in {\cal A}_r,\;2\leq r<j$. 
 Let $$L=\{k\in\{1,2,\ldots,n\}:\; {\rm there\; exist\; at\; least\;
two}\; i_1\neq i_2\in\{1,\ldots ,d\}\;$$
$$ {\rm with} \;{\rm supp}f_{i_l}\cap
{\rm supp}x_k\neq\emptyset,\;l=1,2\}.$$
Then,

{\rm (a)} $|\varphi (\sum_{k\in L}a_kx_k)|\leq\frac{1}{m_j^4}$.

{\rm (b)} $|\varphi (\sum_{k=1}^n a_kx_k)|\leq\frac{2}{m_r}$.}

\medskip

\noindent {\it Proof:} (a) Let $\{l_1,\ldots,l_n\}\in{\cal M}_j $ be such that
$2\leq x_1<l_1<x_2\leq l_2<\ldots \leq l_n .$
Let $n_i=\min {\rm supp}f_i,\; i=1,\ldots d.$
 Then $\{ n_i: i=1,\ldots , d\}\in
{\cal M}_r .$ For each $k\in L$, let $i_k=\min \{i:{\rm supp}f_i 
\;{\rm intersects}\; {\rm supp}x_k \}.$ The map $k\rightarrow n_{i_k}$ from
$L$ to $\{ n_i : i=1,\ldots d \}$ is 1--1, so $\#L\leq d.$ Moreover, $n_{i_k}
\leq l_k$ for each $k\in L$, so $\{l_k :k\in L\}$ belongs to ${\cal M}_r .$
It follows that $\sum _{k\in L}a_k <\frac{1}{m^4_j}$ and so,

$$ |\varphi (\sum _{k\in L}a_kx_k)|\leq \sum_{k\in L}a_k\|x_k\| <\frac{1}{m^4_j}.$$

(b) Let $P=\{ 1,\ldots ,n\}\backslash L$ and, for each $i=1,\ldots ,d$, let

\noindent $P_i=\{ k\in P : {\rm supp}x_k \cap {\rm supp}f_i \neq \emptyset \}.$
Then

$$ |\varphi (\sum _{k=1}^n a_kx_k)|\leq
\frac{1}{m_r}\sum_{i=1}^d|f_i(\sum_{k\in P_i}a_kx_k)|+
 \sum_{k\in L}a_k\|x_k\| <\frac{1}{m_r}+\frac{1}{m^4_j}<\frac{2}{m_r}.\;\;\;
\Box $$

\bigskip

In the sequel we shall write $\tilde{K}\prec K$ if $\tilde{K}$ is a
subset of $K$ satisfying the following.

(i) For every $f\in\tilde{K}$ there exists an analysis
$\{K^s(f)\}$ such that $\cup K^s(f)\subset \tilde{K}$.

(ii) If $f\in\tilde{K}$ then $-f\in \tilde{K}$  and $f|{[m,n]}\in\tilde{K}$
for all $m<n\in {\bf N}$.

(iii) If $(f_i)_{i=1}^d$ is an ${\cal S}$-allowable family in
  $\tilde{K}$  then $\frac{1}{2}\sum_{i=1}^d f_i$ belongs to $\tilde{K}.$

\medskip

For $\tilde{K}\prec K$ we denote by $\| \cdot \|_{\tilde{K}}$ the norm
induced by $\tilde{K}$:
$$\| x\|_{\tilde{K}}=\sup\{ f(x):f\in\tilde{K}\}$$

The results that follow involve a subset $\tilde{K}$ of $K$ having the 
properties mentioned above. For the purposes of this section we only
need these results with $\tilde{K}=K.$ However, we find it convenient
to present them now in the more general formulation that we will need in
Section 3.

\bigskip

\noindent {\bf 2.8 Definition.} Let $\tilde{K}\prec K$. A finite
block sequence $(x_k)_{k=1}^n$ is said to be a rapidly increasing
sequence (R.I.S.) with respect to $\tilde{K}$ if there exist
integers $j_1,\ldots ,j_n$ satisfying the following:

\medskip

(i) $2\leq j_1<j_2<\ldots <j_n$.

\medskip

(ii) Each $x_k$ is a seminormalized $(\frac{1}{m_{j_k}^4},j_k)$-s.c.c.
with respect to $\tilde{K}$. That is, $x_k$ is a $(\frac{1}{m^4_{j_k}},
j_k)$-s.c.c. of the form $x_k=\sum_{t}a_{(k,t)}x_{(k,t)}$ where
$\|x_{(k,t)}\|_{\tilde{K}}=1$ for each $t$, and 
 $\|x_k\|_{\tilde{K}}\geq\frac{1}{2}.$

\medskip

(iii) For $k=1,2,\ldots ,n$, let $l_k=\max {\rm supp}x_k$
and let $n_k\in {\bf N}$ be such that
$$\frac{l_k}{2^{n_k}}<\frac{1}{m_{j_k}}.$$
\noindent We set
$$O_{x_k}=\{f\in K:{\rm supp}f\subset [1,l_k]
\;{\rm and}\;|f(e_m)|>\frac{1}{2^{n_k}}\;{\rm for}\;{\rm every }\;
m\in {\rm supp}f\}.$$
\noindent Then $j_{k+1}$ is such that $m_{j_{k+1}}>\#O_{x_k}$ and
$x_{k+1}$ satisfies $\min {\rm supp}x_{k+1}>\#O_{x_k}$.

\medskip

(iv) $\|x_k\|_{\ell_1}\leq m_{j_{k+1}}/m_{j_{k+1}-1}$.

\bigskip

\noindent {\bf Notation.} If $\varphi\in K\backslash K^0$ then
$\varphi $ is of the form $\varphi =\frac{1}{m_r}\sum_{i=1}^df_i$,
where either $r=1$ and $(f_i)_{i=1}^d$ is an ${\cal S}$-allowable family
of functionals in $K$, or $r\geq 2$ and $(f_i)_{i=1}^d$ is a
${\cal M}_r$-admissible family of functionals in $K$. In either
case we set $w(\varphi )=\frac{1}{m_r}$ (the {\it weight} of
$\varphi $). That is, $w(\varphi )=\frac{1}{m_r}$ if and only if $\varphi \in
{\cal A}_r.$

\bigskip

The following Proposition is the central result of this section:

\medskip

\noindent {\bf 2.9 Proposition.} {\it Let $\tilde{K}\prec K$.
Let $(x_k)_{k=1}^n$ be a R.I.S. with respect to $\tilde{K}$ and
let $\varphi\in\tilde{K}$. There exists a functional $\psi\in K'$
with $w(\varphi )=w(\psi )$ and vectors $u_k,k=2,\ldots ,n$,
 with $\| u_k\|_{\ell_1}\leq 16$ and 
${\rm supp}u_k\subset {\rm supp}x_k$ for each $k$, such that
$$|\varphi (\sum_{k=1}^n\lambda_kx_k)|\leq \max_{1\leq k\leq n}|\lambda_k|+
\psi (\sum_{k=2}^n|\lambda_k|u_k)+6\sum_{k=1}^n|\lambda_k|\frac{1}{m_{j_k}}$$
for every choice of coefficients 
$\lambda_1,\ldots ,\lambda_n\in {\bf R}$. }

\medskip

As it follows from the above statement, we reduce the estimation
of the action of $\varphi$ on the R.I.S. $\{ x_k\}_k$ to the estimation
of the action of the functional $\psi$ on a finite block sequence 
$\{ u_k\}_k$ of subconvex combinations of the basic vectors.
The construction of the functional $\psi$ and the finite block sequence
$\{ u_k\}_k$ will be done in several steps. We describe this process
briefly:

We fix an analysis $\{ K^s(\varphi )\}$ of the functional $\varphi$.
We first replace each vector $x_k$ by its `essential part' relative to
$\varphi$,
denoted by $\overline{x}_k$. Next, for each $\overline{x}_k$ we consider
certain families of functionals in $\cup K^s(\varphi )$ which fall
under two types (families of type I and type II, Definition 2.11). These
families yield a partition of the 
support of $\overline{x}_k$. The restriction from 
$x_k$ to $\overline{x}_k$ gives us a control on the number of families
of type I and type II which act on each $\overline{x}_k$ (Lemma 2.13).
Fixing $k$, to each such family of functionals acting on $\overline{x}_k$,
we correspond a subconvex combination of the basis and the sum of these
combinations is the vector $u_k$. The functional $\psi$ is defined
inductively, following the analysis of the functional $\varphi$.

\bigskip

 From now on we fix the R.I.S.
$(x_k)_{k=1}^n$ and the functional $\varphi$ of Proposition 2.9. We also
fix an analysis
 $\{K^s(\varphi )\}$  of $\varphi $ contained in $\tilde{K}$. We first partition
each vector $x_k$ into three disjointly supported vectors $x'_k,\; x''_k
$ and $\overline{x}_k;$ this partition depends on the analysis
$\{ K^s(\varphi )\}$.

Let
$$F_k=\{ f\in\cup K^s(\varphi ): {\rm supp}f\cap {\rm supp}x_k\neq
\emptyset ,\; {\rm supp}f\cap {\rm supp}x_j\neq\emptyset \;{\rm for}
\;{\rm some}\;j>k$$
$${\rm and}\;w(f )\leq 1/m_{j_{k+1}}\}.$$
We set $A_k=\cup _{f\in F_k}{\rm supp}f$ and $x'_k=x_k|A_k.$

Let now
$$F_k^{\prime }=\{ f\in\cup K^s(\varphi ):\; |f(e_m)|\leq
1/2^{n_k}\;{\rm for}\;
{\rm every}\;m\in {\rm supp}f\cap {\rm supp}(x_k-x_k^{\prime })$$ $$
{\rm and}\;{\rm supp}f\cap {\rm supp}(x_j-x_j^{\prime })\neq\emptyset
\;{\rm for}\;{\rm some}\;j>k\}.$$

We set $A'_k=\cup _{f\in F'_k}{\rm supp}f$ and $x''_k=(x_k-x'_k)|A'_k.$

Finally, $\overline{x}_k={x_k-x'_k-x''_k}$

\bigskip

\noindent {\bf 2.10 Lemma.}
{\it For $\varphi (x'_k)$ and $\varphi (x''_k)$ we have the following estimates:}

$${\rm (1)}\;\; 
|\varphi (x_k^{\prime })|\leq
\frac{1}{m_{j_{k+1}-1}}\;\;\;\;\; {\rm  and}\;\;\;\; {\rm (2)}\;\;
|\varphi (x_k^{\prime\prime })|
<1/m_{j_k}.$$
{\it Proof:} 
To see (1), let us call an $f\in F_k$ maximal if there is no $f'\neq f$ in
$F_k$ such that ${\rm supp}f\subset{\rm supp}f'.$ The maximal elements
of $F_k$ have disjoint supports. So,
$$|\varphi (x_k^{\prime })|\leq\sum_{f\;{\rm maximal}\;{\rm in}\;F_k}
|f(x_k^{\prime })|\leq\sum_f
\frac{1}{m_{j_{k+1}}}\|x_k|_{{\rm supp}f}\|_{\ell_1}
\leq\frac{1}{m_{j_{k+1}}}\frac{m_{j_{k+1}}}{m_{j_{k+1}-1}}=
\frac{1}{m_{j_{k+1}-1}},$$
by property (iv) of the R.I.S.

For (2), we notice that 
for every $n\in {\rm supp}x_k^{\prime\prime }$ we have
$|\varphi (e_n)|\leq 1/2^{n_k}$, hence
$$|\varphi (x_k^{\prime\prime })|\leq\frac{\|x_k\|_{\ell_1}}{2^{n_k}}\leq
\frac{\max
{\rm supp}x_k}{2^{n_k}}<\frac{1}{m_{j_k}}.$$

\medskip

\noindent {\bf Remarks:} (1) By the definition of $x_k^{\prime }$
and $x_k^{\prime\prime }$ we have $x_n^{\prime }=x_n^{\prime\prime }=0$,
since $x_n$ is the last element of $(x_k)_1^n$.

(2) If $f\in\cup K^s(\varphi )$ and $1\leq k<l\leq n$ are such that
${\rm supp}f\cap{\rm supp}\overline{x}_k\neq\emptyset$ 
and ${\rm supp}f\cap{\rm supp}\overline{x}_l\neq\emptyset$ then
$w(f)>\frac{1}{m_{j_{k+1}}}$ and there exists $m\in {\rm supp}\overline{x}_k$
such that $|f(e_m)|>\frac{1}{2^{n_k}}.$

\medskip

\noindent {\bf 2.11 Definition of the families of type--I and type--II w.r.t. $
\overline{x}_k$:}

Without loss of generality, we assume that ${\rm supp}\varphi\cap{\rm supp}
\overline{x}_1\neq\emptyset.$
Let $k\in\{2,\ldots n\}$ be fixed.

(A) A set of functionals $F=\{f_1,\ldots ,f_l\}$  contained in
some level $K^s(\varphi)$ of the analysis
of $\varphi $ is said to be a family of type--I with respect to $\overline{x}_k$ if

(A1) ${\rm supp}f_i\cap {\rm supp}\overline{x}_k\neq\emptyset $ and
${\rm supp}f_i\cap {\rm supp}\overline{x}_j=\emptyset $ for every
$j\neq k$ and every $i=1,2,\ldots ,l$.

(A2) There exists $g\in K^{s+1}(\varphi )$ such that
$f_1,\ldots ,f_l$ belong to the decomposition of $g$
and ${\rm supp}g\cap {\rm supp}\overline{x}_j\neq\emptyset $
for some $j<k$. 
Moreover, $F$ is the maximal subset of the decomposition
of $g$ with property A1; that is, $g=\frac{1}{m_r}(\sum _{i=1}^d h_i +
\sum_{i=1}^l f_i),$ where, for each $i=1,\ldots , d$, either 
${\rm supp}
h_i\cap {\rm supp}\overline{x}_k=\emptyset$ 
or
${\rm supp}
h_i\cap {\rm supp}\overline{x}_j\neq \emptyset$ for some $j\neq k.$

\medskip

(B) A set of functionals $F=\{f_1,\ldots ,f_m\}$ contained in some level
$K^s(\varphi )$ of  the analysis
of $\varphi $ is said to be a family of type--II with respect to
$\overline{x}_k$ if

(B1) ${\rm supp}f_i\cap {\rm supp}\overline{x}_k\neq\emptyset $,
${\rm supp}f_i\cap {\rm supp}\overline{x}_j=\emptyset $
for every $j<k$ and every $i=1,2,\ldots ,m$,
and for every $i=1,2,\ldots ,m$ we can find $j_i>k$ such that
${\rm supp}f_i\cap {\rm supp}\overline{x}_{j_i}\neq\emptyset $.

(B2) There exists $g\in K^{s+1}(\varphi )$ such that
$f_1,\ldots ,f_m$ belong to the  decomposition
of $g$ and ${\rm supp}g\cap {\rm supp}\overline{x}_j\neq\emptyset $
for some $j<k$.
Moreover, $F$ is the maximal subset of the decomposition
of $g$ with property B1; that is, $g=\frac{1}{m_r}(\sum _{i=1}^d h_i +
\sum_{i=1}^m f_i),$ where, for each $i=1,\ldots , d$, either 
${\rm supp}
h_i\cap {\rm supp}\overline{x}_k=\emptyset$ 
or
${\rm supp}
h_i\cap {\rm supp}\overline{x}_j\neq \emptyset$ for some $j<k$ or
${\rm supp}h_i\cap {\rm supp}\overline{x}_j= \emptyset$ for all $j\neq k.$

\bigskip

\noindent {\bf Remarks:} (1) It is easy to see that for $k=2,3,\ldots ,n$,
$${\rm supp}\overline{x}_k\cap {\rm supp}\varphi =
{\rm supp}\overline{x}_k\cap\bigcup\{ \cup_{f\in F}{\rm supp}f : F\;{\rm  is\; a
\; family\; of\; type\; I\; or\; type\; II\; w.r.t.\;} \overline{x}_k\}.$$

 (2) Let $k$ be fixed. If each of the families
 $\{f_1,\ldots ,f_l\}$ 
 and $\{f'_1,\ldots ,f'_m\}$   is   of type I or  of type II
w.r.t. $\overline{x}_k$  and they are not identical, 
then, for all $i\leq l$, $j\leq m$, ${\rm supp}f_i 
\cap {\rm supp}f'_j=\emptyset$.

 (3) Let $F$ be a family of type I or type II w.r.t. $\overline{x}_k$
and let $g_F$ be the functional in $\cup K^s(\varphi )$ which contains $F$
in its decomposition. Then $g_F$ intersects $\overline{x}_j$ for some $j<k$.
By Remark (2) after Lemma 2.10 this implies that $w(g_F)>\frac{1}{m_{j_k}}$.

\medskip

\noindent {\bf 2.12 Lemma.} {\it Let $2\leq k\leq n$. 
 If $f$ is a member of a family of
type--I or type--II with respect to $\overline{x}_k$,  then there
exist  sets $A_{k,f}, A'_{k,f}\subset {\rm supp}f$ satisfying}
$$ |f(x_k^{\prime })|\leq
\frac{1}{m_{j_{k+1}}}\|x_k|_{A_{k,f}}\|_{\ell_1}\;\;\;\;
{\rm and}$$
$$ |f(x_k^{\prime\prime })|\leq
\frac{1}{2^{n_k}}\|x_k|_{A'_{k,f}}\|_{\ell_1}.$$
{\it Moreover, if $f$ and $f'$ are two distinct such functionals then
$A_{k,f}\cap A_{k,f'}=\emptyset$ and 
$A'_{k,f}\cap A'_{k,f'}=\emptyset.$  }

\medskip

\noindent {\it Proof:} If $f(x'_k)\neq 0$ then, by the definition of $x'_k$,
either there exists $g\in F_k$ with ${\rm supp}f \subset {\rm supp}g$ or
there exists $g\in F_k$ with ${\rm supp}g \subset {\rm supp}f$. But the 
first case is impossible because then we would have ${\rm supp}f\cap
{\rm supp}x_k \subset {\rm supp}x'_k$ and so ${\rm supp}f\cap{\rm supp}\overline
{x}_k=\emptyset.$ So, if we set
$$A_{k,f}=\bigcup\{ {\rm supp}g\cap{\rm supp}x_k: g\in F_k \; {\rm and}\;
{\rm supp}g\subset {\rm supp}f \},$$
then $f(x'_k)=f(x_k|A_{k,f})$. This gives 
$$|f(x'_k)|\leq \frac{1}{m_{j_{k+1}}}\| x_k|A_{k,f}\|_{\ell_1}.$$

In the same way, if $f(x''_k)\neq 0$ we set 
$$A'_{k,f}=\bigcup \{ {\rm supp}g\cap {\rm supp}(x_k-x'_k): g\in F'_k \; {\rm and}
\; {\rm supp}g\subset {\rm supp }f \}.$$
Then $f(x''_k)=f(x_k|A'_{k,f})$, so
$$|f(x''_k)|\leq \frac{1}{2^{n_k}}\|x_k|A'_{k,f}\|_{\ell_1}.\;\;\;$$
The disjointness follows from the preceding Remark (2).$\;\;\;\Box $

\medskip

\noindent {\bf 2.13 Lemma.} {\it Let $k=2,3,\ldots ,n$. Then,

{\rm (a)} The number of families of type I w.r.t $\overline{x}_k$ is less
than $\min {\rm supp}x_k$.

{\rm (b)} The number of families of type II w.r.t $\overline{x}_k$ is less
than $\min {\rm supp}x_k$.}

\medskip

\noindent {\it Proof:} (a) For each family $F$ of type I w.r.t $\overline{x}_k$
let $g_F$ be the (unique) functional in $\cup K^s(\varphi )$
which contains $F$ in its decomposition. 

By the maximality of $F$
in the decomposition of $g_F$, it is clear that if $F\neq F^{\prime }$ are
two families of type I then $g_F\neq g_{F^{\prime }}$. Since both
$g_F$ and $g_{F^{\prime }}$ are elements of the analysis of
$\varphi $, it follows that either ${\rm supp}g_F\subset {\rm supp}g_{F^{\prime }}$ or ${\rm supp}g_{F^{\prime }}\subset {\rm supp}g_F$ or
${\rm supp}g_F\cap {\rm supp}g_{F^{\prime }}=\emptyset $. In either case
$g_F(e_k)\neq g_{F^{\prime }}(e_k)$ for all $k$.
Moreover, for each $F$, $g_F$ has the property that
${\rm supp}g_F\cap {\rm supp}\overline{x}_i\neq\emptyset $ for
some $i<k$. Let $i_F=\min\{i: {\rm supp}g_F\cap {\rm supp}\overline{x}_i
\neq\emptyset \}$. It follows from  Remark 2 after Lemma 2.10
that there exists $m_F$ in ${\rm supp}\overline{x}_{i_F}$ with
$|g_F(e_{m_F})|>1/2^{n_{i_F}}$.

So, for each family $F$ of type I w.r.t $\overline{x}_k$, we set
$h_F=g|_{\{m_F\}}\in K$. The map $F\rightarrow h_F$ is one to one;
moreover, each $h_F$ belongs to $O_{x_k}$ (see Definition 2.8).

It follows that
$$\#\{F: F\;{\rm is}\;{\rm a}\;{\rm family}\;{\rm of}\;{\rm type}
\;{\rm I}\;{\rm w.r.t}\;\overline{x}_k\}\leq\#O_{x_k}<\min {\rm supp}x_{k+1}.$$

\noindent (b) The proof is the same as that of part (a).   $\Box $

\bigskip

\noindent {\bf Notation:} For each $k=2,3,\ldots ,n$, we classify
the families of type--I and type--II into four classes according
to the weight $w(g_F)$ of the functional $g_F$ which contains
each family $F$ in its decomposition. We set:

$$A_{\overline{x}_k}=\{F: F\;{\rm is}\;{\rm a}\;{\rm family}\;{\rm of}
\;{\rm type}\;{\rm I}\;{\rm w.r.t} \;\overline{x}_k\;{\rm and}
\;w(g_F)=\frac{1}{2}\},$$

$$B_{\overline{x}_k}=\{F: F\;{\rm is}\;{\rm a}\;{\rm family}\;{\rm of}
\;{\rm type}\;{\rm I}\;{\rm w.r.t}\;\overline{x}_k\;{\rm and}
\;w(g_F)<\frac{1}{2}\},$$

$$C_{\overline{x}_k}=\{F: F\;{\rm is}\;{\rm a}\;{\rm family}\;{\rm of}
\;{\rm type}\;{\rm II}\;{\rm w.r.t}\;\overline{x}_k\;{\rm and}
\;w(g_F)=\frac{1}{2}\},$$

$$D_{\overline{x}_k}=\{F: F\;{\rm is}\;{\rm a}\;{\rm family}\;{\rm of}
\;{\rm type}\;{\rm II}\;{\rm w.r.t}\;\overline{x}_k\;{\rm and}
\;w(g_F)<\frac{1}{2}\}.$$

\medskip

\noindent {\bf Remarks:} (1) If $F\in D_{\overline{x}_k}$, then $F$
is a singleton, i.e. $F=\{f\}$. Because, if
$g_F=\frac{1}{m_s}(\sum h_i+\sum_{i=1}^mf_i)$ where $s\geq 1$
and $F=\{f_1,\ldots ,f_m\}$, then $f_1<f_2<\ldots <f_m$, and
each ${\rm supp}f_i$ intersects ${\rm supp}\overline{x}_k$
and ${\rm supp}\overline{x}_{j_i}$ for some $j_i>k$. This is
impossible unless  $m=1$.

\medskip

(2) If $f^{\prime }<f<f^{\prime\prime }$ belong to $\cup K^s(\varphi )$
and there exists a family of type II w.r.t $\overline{x}_k$
which is contained in the analysis of $f$, then
${\rm supp}f^{\prime }\cap {\rm supp}\overline{x}_k=\emptyset $
and ${\rm supp}f^{\prime\prime }\cap {\rm supp}\overline{x}_k=\emptyset $.

\bigskip

\noindent {\bf Notation.}

\medskip

A. Each $x_k$ is a seminormalized $(\frac{1}{m_{j_k}^4},j_k)$-s.c.c of
the form 
$$x_k=\sum_{t=1}^{r_k}a_{(k,t)}x_{(k,t)},$$
\noindent where $a_{(k,t)}\geq 0$, $\sum_ta_{(k,t)}=1$ and
$\|x_{(k,t)}\|_{\tilde{K}}=1$.

For each $k=1,\ldots ,n$, $t=1,\ldots ,r_k$, we set
$$\overline{x}_{(k,t)}=x_{(k,t)}|_{{\rm supp}\overline{x}_k}.$$

B. Fix $k\in\{2,\ldots ,n\}$. If $f\in\cup K^s(\varphi )$ is
a member of a family of type I or type II w.r.t $\overline{x}_k$,
we set
$$n_f=\min ({\rm supp}\overline{x}_k\cap {\rm supp}f)\;\;\;\;{\rm and}
\;\;\;\;e_f=e_{n_f}.$$

\noindent Also, if $F=\{f_1,\ldots ,f_l\}$ is a family of type I or type II
w.r.t $\overline{x}_k$, then we set
$$n_F=\min ({\rm supp}\overline{x}_k
\cap (\cup_{i=1}^l{\rm supp}f_i))\;\;\;\;{\rm and}\;\;\;\; e_F=e_{n_F}.$$

\noindent For $F=\{f_1,\ldots ,f_l\}\in A_{\overline{x}_k}
\cup C_{\overline{x}_k}$
we set 
$$h_F=\frac{1}{2}(f_1+\ldots +f_l)\;\;\;\;{\rm and}\;\;\;\;
a_F=|2h_F(\overline{x}_k)|.$$

\noindent For $\{f\}\in D_{\overline{x}_k}$ we set
$$a_f=|f(\overline{x}_k)|.$$

\noindent Finally, if $F\in B_{\overline{x}_k}$, for every $f\in F$
we set
$$\Omega_{f}=\{t:{\rm supp}f\cap {\rm supp}\overline{x}_{(k,t)}\neq
\emptyset\;{\rm and}\;{\rm supp}h\cap {\rm supp}\overline{x}_{(k,t)}=
\emptyset\;{\rm for}\;{\rm every}\; h\neq f {\rm \; in\; }F\}$$
\noindent and
$$a_{f}=\sum_{t\in\Omega_{f}}a_{(k,t)}|f(\overline{x}_{(k,t)})|.$$

C. For each $k=2,3,\ldots ,n$ we define
$$u_k=\sum_{\{f\}\in D_{\overline{x}_k}}a_fe_f+
\sum_{F\in A_{\overline{x}_k}\cup C_{\overline{x}_k}}a_Fe_F+
\sum_{F\in B_{\overline{x}_k}}\sum_{f\in F}a_fe_f.$$

\noindent {\bf 2.14 Lemma.} {\it For $k=2,3,\ldots ,n$, }
$$\| u_k\|_{\ell_1}=\sum_{\{f\}\in D_{\overline{x}_k}}a_f
+\sum_{F\in A_{\overline{x}_k}\cup C_{\overline{x}_k}}a_F+
\sum_{F\in B_{\overline{x}_k}}\sum_{f\in F}a_f\leq 16.$$

\newcommand{\oxk}{\overline{x}_k}
\newcommand{\fd}{\{f\}\in D_{\overline{x}_k}}
\newcommand{\sfd}{\sum_{\fd}}

\noindent {\it Proof:} For each $f$ with $\fd $, set $\varepsilon_f=
{\rm sign}(f(\oxk ))$. Then,
$$\sfd a_f=\sfd |f(\oxk )|=\sfd\varepsilon_ff(\oxk )$$
$$=\sfd\varepsilon_ff(x_k)-\sfd\varepsilon_ff(x^{\prime }_k)-
\sfd\varepsilon_ff(x^{\prime\prime }_k)$$
$$\leq\sfd\varepsilon_kf(x_k)+\sfd |f(x^{\prime }_k)|+
\sfd |f(x^{\prime\prime }_k)|$$
$$\leq\sfd\varepsilon_ff(x_k)+\frac{1}{m_{j_{k+1}}}\sfd
\|x_k|_{A_{k,f}}\|_{\ell_1}+\frac{1}{2^{n_k}}\sfd\|x_k|_{A^{\prime }_{k,f}}\|_{\ell_1},$$
\noindent where the last inequality follows from Lemma 2.12. From the same
lemma and Definition 2.8 we get
$$\frac{1}{m_{j_{k+1}}}\sfd\|x_k |_{A_{k,f}}\|_{\ell_1}\leq
\frac{1}{m_{j_{k+1}}}\|x_k \|_{\ell_1}<\frac{1}{m_{j_k}}$$
\noindent and
$$\frac{1}{2^{n_k}}\sfd\|x_k |_{A^{\prime }_{k,f}} \|_{\ell_1}\leq
\frac{1}{2^{n_k}}\|x_k \|_{\ell_1}\leq\frac{l_k}{2^{n_k}}<
\frac{1}{m_{j_k}}.$$
For every $f\in\tilde{K}$ we have that $\varepsilon_ff|_{[\min {\rm supp}
x_k,\infty )}\in\tilde{K}$. Also, by Remark (2) following Definition 2.11,
we have that if $f\neq f^{\prime }$ and both $\{f\}$ and $\{f^{\prime }\}$
are families of type II w.r.t $\oxk $, then
${\rm supp}f\cap {\rm supp}f^{\prime }=\emptyset $. By Lemma 2.13 we have
$\# D_{\oxk }<\min {\rm supp}x_k$.
It follows that the set $\{\varepsilon_ff|_{[\min ({\rm supp}x_k),\infty )}:
\{f\}\in D_{\oxk }\}$ is ${\cal S}$-allowable, and so the
functional $\frac{1}{2}\sfd\varepsilon_ff_{[\min ({\rm supp}x_k),\infty )}$
belongs to $\tilde{K}$. We conclude that $|\frac{1}{2}\sum\varepsilon_ff(x_k)|\leq\|x_k\|_{\tilde{K}}\leq 1$, and so,
$$\sfd a_f\leq 2+\frac{2}{m_{j_k}}<3.\leqno (1)$$

\newcommand{\sfc}{\sum_{F\in C_{\oxk}}}

For $F\in C_{\oxk }$ we set $\varepsilon_F={\rm sign}h_F(\oxk )$. Then,
$$\sfc a_F=\sfc |2h_F(\oxk )|=2\sfc\varepsilon_Fh_F(\oxk )$$
$$=2\left ( \sum\varepsilon_Fh_F(x_k)-\sum\varepsilon_Fh_F(x^{\prime }_k)
-\sum\varepsilon_Fh_F(x^{\prime\prime }_K)\right )$$
$$\leq 2\sfc\varepsilon_Fh_F(x_k)+2\sfc\sum_{f\in F}|f(x^{\prime }_k)|+
2\sfc\sum_{f\in F}|f(x^{\prime\prime }_k)|$$
$$2\sfc\varepsilon_Fh_F(x_k)+\frac{2}{m_{j_{k+1}}}\sfc\sum_{f\in F}
\|x_k|_{A_{k,f}}\|_{\ell_1}+\frac{2}{2^{n_k}}\sfc\|x_k|_{A^{\prime }_{k,f}}\|_{\ell_1}$$
$$\leq 2\sfc\varepsilon_Fh_F(x_k)+\frac{4}{m_{j_k}},$$
\noindent again by Lemma 2.12. On the other hand, for
$F=\{f_1,\ldots ,f_l\}\in C_{\oxk }$, $h_F=\frac{1}{2}(f_1+\ldots +f_l)\in
\tilde{K}$ and $\varepsilon_Fh_F\in\tilde{K}$. By Lemma 2.13 we
have that $\#C_{\oxk }<\min {\rm supp}x_k$ and by Remark (2) after
2.11 we have that the functionals $h_F$, $F\in C_{\oxk },$ are
disjointly supported. We conclude that the set
$\{h_F|_{[\min{\rm supp}x_k,\infty )}:F\in C_{\oxk }\}$ is
${\cal S}$-allowable and so, the functional
$\frac{1}{2}\sfc\varepsilon_Fh_F|_{[\min {\rm supp}x_k,\infty )}$
belongs to $\tilde{K}$ and
$$|\sfc\varepsilon_Fh_F(x_k)|\leq 2\|x_k\|\leq 2.$$
\noindent We conclude that
$$\sfc a_F\leq 4+\frac{4}{m_{j_k}}<5.\leqno (2)$$
\noindent In the same way we get
$$\sum_{F\in A_{\oxk }}a_F<5.\leqno (3)$$

Finally, we have

$$\sum_{F\in B_{\overline{x}_k}}\sum_{f\in F}a_f=
\sum_{F\in B_{\overline{x}_k}}\sum_{f\in F}\sum_{t\in\Omega_f}a_{(k,t)}
|f(\overline{x}_{(k,t)})|$$
$$\leq\sum_{F\in B_{\overline{x}_k}}\sum_{f\in F}\sum_{t\in\Omega_f}
a_{(k,t)}\left (|f(x_{(k,t)})|+|f(x^{\prime }_{(k,t)})|+
|f(x^{\prime\prime }_{(k,t)})|\right ).$$

\noindent For each $F\in B_{\overline{x}_k}$ and $f\in F$ we have
$$\sum_{t\in\Omega_f}a_{(k,t)}|f(x^{\prime }_{(k,t)})|\leq
\frac{1}{m_{j_{k+1}}}\|x_k|_{A_{k,f}}\|_{\ell_1}$$
\noindent and
$$\sum_{t\in\Omega_f}a_{(k,t)}|f(x^{\prime\prime }_{(k,t)})|\leq
\frac{1}{2^{n_k}}\|x_k|_{A^{\prime }_{k,f}}\|_{\ell_1}.$$

\noindent Since the sets $A_{k,f}$, $f\in\cup_{F\in B_{\overline{x}_k}}F$
are disjoint, we get

$$\sum_{F\in B_{\overline{x}_k}}\sum_{f\in F}
\sum_{t\in\Omega_f}a_{(k,t)}|f(x^{\prime }_{(k,t)})|\leq
\frac{1}{m_{j_{k+1}}}\|x_k\|_{\ell_1}<\frac{1}{m_{j_k}}.\leqno ({\rm i})$$

\noindent In a similar way,

$$\sum_{F\in B_{\overline{x}_k}}\sum_{f\in F}
\sum_{t\in\Omega_f}a_{(k,t)}|f(x^{\prime\prime }_{(k,t)})|
\leq\frac{1}{2^{n_k}}\|x_k\|_{\ell_1}<\frac{1}{m_{j_k}}.\leqno ({\rm ii})$$

\noindent It remains to estimate
$$\sum_{F\in B_{\overline{x}_k}}\sum_{f\in F}\sum_{t\in\Omega_f}
a_{(k,t)}|f(x_{(k,t)})|.$$
\noindent For each $F\in B_{\overline{x}_k}$ and $t\in\cup_{f\in F}
\Omega_f$, let $f_t^F$ be the unique element of $F$ with
$f_t^F(\overline{x}_{(k,t)})\neq 0$. Let also,
$\Omega_F=\cup_{f\in F}\Omega_f$ and $\Omega =\cup_{F\in B_{\overline{x}_k}}
\Omega_F$. Then,
$$\sum_{F\in B_{\overline{x}_k}}\sum_{f\in F}\sum_{t\in\Omega_f}
a_{(k,t)}|f(x_{(k,t)})|=
\sum_{F\in B_{\overline{x}_k}}\sum_{t\in\Omega_F}a_{(k,t)}|f_t^F(x_{(k,t)})|$$
$$=\sum_{t\in\Omega }a_{(k,t)}\sum_{F\in B_{\overline{x}_k}}
|f_t^F(x_{(k,t)})|.$$
\noindent Fix $t\in\Omega $. For each $F\in B_{\overline{x}_k}$, we set
$\varepsilon_F={\rm sign}f_t^F(x_{(k,t)})$.
Since $\# B_{\overline{x}_k}<\min {\rm supp}x_k$, the functional
$$h=\frac{1}{2}\sum_{F\in B_{\overline{x}_k}}\varepsilon_F
f_t^F|_{[\min {\rm supp}x_k,\infty )}$$
\noindent belongs to $\tilde{K}$. So, we get
$$\sum_{F\in B_{\overline{x}_k}}|f_t^F(x_{(k,t)})|
=2h(x_{(k,t)})\leq 2\|x_{(k,t)}\| =2.$$
\noindent We conclude that

$$\sum_{t\in\Omega }a_{(k,t)}\sum_{F\in B_{\overline{x}_k}}
|f_t^F(x_{(k,t)})|\leq 2\sum_{t\in\Omega }a_{(k,t)}\leq 2.
\leqno ({\rm iii})$$
\noindent Finally, by (i), (ii) and (iii),
$$\sum_{F\in B_{\overline{x}_k}}\sum_{f\in F}a_f\leq 2+\frac{2}{m_{j_k}}<3.
\leqno (4)$$

\noindent Combining (1), (2), (3), (4) we get the desired
estimate for $\|u_k\|_{\ell_1}$. $\; \Box$

\bigskip

\noindent {\bf 2.15 Lemma.} {\it There exists a functional
$\psi\in K^{\prime }$ with $w(\psi )=w(\varphi )$ and such that,
for $k=2,\ldots ,n$,}
$$|\varphi (\overline{x}_k)|\leq\psi (u_k)+\frac{2}{m_{j_k}}.$$

\medskip

\noindent {\it Proof:} We build the functional $\psi $
inductively, following the way $\varphi $ is built by
the analysis $\cup K^s(\varphi )$.

We first introduce some more notation: For $f\in\cup K^s(\varphi )$,
we set
$$K(f)=\{ f^{\prime }\in\cup K^s(\varphi ):
{\rm supp}f^{\prime }\subset {\rm supp}f\},$$
\noindent that is, $K(f)$ is the analysis of $f$ induced by
$\cup K^s(\varphi )$.

For $f=\frac{1}{m_s}\sum_{i=1}^df_i$ and each $k=2,\ldots ,n$, we set
$$I_k^f=\{ i\in\{1,\ldots ,d\}: f_i\;{\rm is}\;{\rm an}
\;{\rm element}\;{\rm of}\;{\rm a}\;{\rm family}
\;{\rm of}\;{\rm type}\;{\rm I}\;{\rm w.r.t}\;\overline{x}_k\},$$
$$J_k^f=\{ i\in\{1,\ldots ,d\}:f_i\;{\rm is}\;{\rm an}
\;{\rm element}\;{\rm of}\;{\rm a}\;{\rm family}\;{\rm of}
\;{\rm type}\;{\rm II}\;{\rm w.r.t.}\;\overline{x}_k\},$$
\noindent and
$$\Lambda_k^f=\{ i\in\{1,\ldots ,d\} : K(f_i)\;{\rm contains}
\;{\rm a}\;{\rm family}\;{\rm of}\;{\rm type}\;{\rm I}
\;{\rm or}\;{\rm type}\;{\rm II} {\rm \; w.r.t. \;}\overline{x}_k\}.$$
\noindent We also set
$$I^f=\bigcup _{k=2}^nI_k^f,\;\;\; J^f=\bigcup _{k=2}^nJ_k^f,\;\;\;
\Lambda^f=\bigcup_{k=2}^n\Lambda_k^f$$
\noindent and
$$D_f=\bigcup_{k=2}^n\bigcup \{ \cup_{f'\in F} {\rm supp}f^{\prime }:
F \;{\rm is\; a\; family \; of \; type \; I \; or \; type \; II} $$
$${\rm w.r.t.}\;\overline{x}_k
\;{\rm and}\;F\subset K(f)\}.$$

Let $k=2,\ldots ,n$ and let $F$ be a family in $B_{\overline{x}_k}$.
We set
$$L_F=\{ t:\;{\rm there}\;{\rm exist}\;{\rm at}\;{\rm least}
\;{\rm two}\;{\rm functionals}\;h,h^{\prime }\in F\;{\rm such}
\;{\rm that}$$
$${\rm supp}h\cap{\rm supp}\overline{x}_{(k,t)}\neq\emptyset\;{\rm and}\;
{\rm supp}h^{\prime }\cap{\rm supp}\overline{x}_{(k,t)}\neq\emptyset\}.$$

Let $g_F$ be the functional in $\cup K^s(\varphi )$
which contains the family $F$ in its decomposition.
We set
$$C_F=w(g_F)|\;\sum_{t\in L_F}a_{(k,t)}\sum_{f\in F}
f(\overline{x}_{(k,t)})|.$$
\noindent Finally, for $f\in\cup K^s(\varphi )$ we set
$B_k(f)=\{ F\in B_{\overline{x}_k}:F\subset K(f)\}$.

\medskip

By induction on $s=0,\ldots ,m$, for every $f\in K^s(\varphi )$
we shall construct a functional $\psi_f\in K^{\prime }$
such that:

\medskip

If $D_f=\emptyset $, then $\psi_f=0$.

\medskip

If $D_f\neq\emptyset $, then $\psi_f$ has the following properties:

(a) ${\rm supp}\psi_f\subset D_f\subset {\rm supp}f$.

(b) For each $k=2,\ldots ,n$,
$$|f(\overline{x}_k|_{D_f})|\leq\psi_f(u_k)+
\sum_{F\in B_k(f)}C_F.$$ 

(c) $w(\psi_f)=w(f)$.

\medskip

Suppose that $\psi_f$ has been defined for all $f\in\cup_{t=1}^{s-1}
K^t(\varphi )$. Let $f=\frac{1}{m_q}\sum_{i=1}^df_i\in K^s(\varphi )
\backslash K^{s-1}(\varphi )$ be such that $D_f\neq\emptyset $.

\medskip

\noindent {\bf Case 1:} $w(f)=\frac{1}{m_q}<\frac{1}{2}$.

\medskip

Then we set
$$\psi_f=\frac{1}{m_q}\left (\sum_{i\in\Lambda^f}\psi_{f_i}+
\sum_{i\in I^f}e^{\ast }_{f_i}+
\sum_{i\in J^f}e^{\ast }_{f_i}\right ).$$
By the inductive assumption, property (a) is satisfied.

 We note that the sets $\Lambda^f$ and $J^f$ are not disjoint.
If $i\in J_k^f$ then $i\in \Lambda_m^f$ for some $m>k.$ In this case,
${\rm supp} \psi_{f_i}\subset D_{f_i}\subset [\min{\rm supp}\overline{x}_{k+1},
\infty)$, while ${\rm supp}e^*_{f_i}=\{ n_{f_i}\}\subset {\rm supp}
\overline{x}_k.$
It follows that $e^*_{f_i}<\psi_{f_i}.$ 

Fix now $k\in\{ 2,\ldots ,n\}$.
Since $w(f)<\frac{1}{2}$, we have $f_1<f_2<\ldots <f_d$, so each
of the sets $J_k^f$ and $\Lambda_k^f$ is either empty or a
singleton. Suppose that $\Lambda_k^f=\{ i_1\}$ and
$J_k^f=\{ i_2\}$. Then,
$$ |f(\overline{x}_k|_{D_f})|=\frac{1}{m_q}
|f_{i_1}(\overline{x}_k|_{D_f})+\sum_{i\in I_k^f}f_i(\overline{x}_k)+
f_{i_2}(\overline{x}_k)|$$
$$\leq\frac{1}{m_q}|f_{i_1}(\overline{x}_k|_{D_f})|+
\frac{1}{m_q}|\sum_{i\in I_k^f}f_i(\overline{x}_k)|+
\frac{1}{m_q}|f_{i_2}(\overline{x}_k)|.$$

\noindent We have
$$\frac{1}{m_q}|f_{i_1}(\overline{x}_k|_{D_f})|\leq
\frac{1}{m_q}\left (\psi_{f_{i_1}}(u_k)+
\sum_{F\in B_k(f_{i_1})}C_F\right ) \leqno (1)$$
\noindent by the inductive assumption. Also,
$$\frac{1}{m_q}|f_{i_2}(\overline{x}_k)|=
\frac{1}{m_q}|f_{i_2}(\overline{x}_k)|e_{f_{i_2}}^{\ast }
(e_{f_{i_2}})=\frac{1}{m_q}e_{f_{i_2}}^{\ast }
(a_{f_{i_2}}e_{f_{i_2}})=\frac{1}{m_q}e^{\ast }_{f_{i_2}}(u_k).
\leqno (2)$$

Finally, let $G=\{ f_i:i\in I_k^f\}$ be the family of type I w.r.t. 
$\overline{x}_k$ contained
in the decomposition of $f$. Then,
$$\frac{1}{m_q}|\sum_{i\in I_k^f}f_i(\overline{x}_k)|
=\frac{1}{m_q}|\sum_{i\in I_k^f}f_i(\sum_ta_{(k,t)}\overline{x}_{(k,t)})|$$
$$=\frac{1}{m_q}|\sum_{f_i\in G}\sum_{t\in\Omega_{f_i}}
a_{(k,t)}f_i(\overline{x}_{(k,t)})+
\sum_{t\in L_G}a_{(k,t)}(\sum_{f_i\in G}f_i)(\overline{x}_{(k,t)})|$$
$$\leq\frac{1}{m_q}\sum_{f_i\in G}\sum_{t\in\Omega_{f_i}}
a_{(k,t)}|f_i(\overline{x}_{(k,t)})|+
\frac{1}{m_q}|\sum_{t\in L_G}a_{(k,t)}(\sum_{f_i\in G}
f_i(\overline{x}_{(k,t)})|$$
$$=\frac{1}{m_q}\sum_{i\in I_k^f}a_{f_i}+C_G=
\frac{1}{m_q}\sum_{i\in I_k^f}a_{f_i}e_{f_i}^{\ast }(e_{f_i})+C_G
=\frac{1}{m_q}\sum_{i\in I_k^f}e_{f_i}^{\ast }(u_k)+C_G.$$

\noindent So,
$$\frac{1}{m_q}|\sum_{i\in I_k^f}f_i(\overline{x}_k)|\leq
\frac{1}{m_q}\sum_{i\in I_k^f}e_{f_i}^{\ast }(u_k)+C_G.\leqno (3)$$

From (1), (2) and (3) we conclude that property (b) holds for
$\psi_f$, that is,
$$|f(\overline{x}_k|_{D_f})|\leq\psi_f(u_k)+
\sum_{F\in B_k(f)}C_F.$$

It remains to show that $\psi_f\in K^{\prime }$. We have to show
that the set
$$\{\psi_{f_i}:i\in\Lambda^f\}\cup\{e_{f_i}^{\ast }:
i\in I^f\cup J^f\}$$
\noindent is ${\cal M}^{\prime}_{q}$--admissible. For $i=1,\ldots ,d$,
let $r_i=\min ({\rm supp}f_i)$. Then, $\{ r_i:i=1,\ldots ,d\}\in {\cal M}_q$.

To each $i\in I^f$ corresponds the vector $e_{f_i}^{\ast }$
with $r_i\leq e_{f_i}^{\ast }<r_{i+1}$.

If $i\in J^f$, then
$i\in\Lambda^f$ also, so to it correspond two vectors
$e_{f_i}^{\ast }$ and $\psi_{f_i}$ with
$r_i\leq e_{f_i}^{\ast }<\psi_{f_i}<r_{i+1}$.

Finally, if $i\in\Lambda^f\backslash J^f$, then to it
corresponds the vector $\psi_{f_i}$ with
$r_i\leq \psi_{f_i}<r_{i+1}$.

It follows from these relations that the family
$$\{ \psi_{f_i}:i\in\Lambda^f\}\cup\{e_{f_i}^{\ast }:
i\in I^f\cup J^f\}$$
\noindent is ${\cal M}_q^{\prime }$--admissible, 
 and since
$\psi_{f_i}$, $e_{f_i}^{\ast }\in K^{\prime }$, we get
$\psi_f\in {\cal A}_q^{\prime }$.

\medskip

\noindent {\bf Case 2:} $w(f)=\frac{1}{m_q}=\frac{1}{2}$.

\medskip

For each $k=2,\ldots ,n$, let $F_1^k=\{ f_i:i\in I_k^f\}$
be the family of type I w.r.t. $\overline{x}_k$
 contained in the decomposition of $f$,
and let $F_2^k=\{ f_i:i\in J_k^f\}$ be the family of type II
 w.r.t. $\overline{x}_k$ contained in the decomposition of $f$. We set
$$\psi_f=\frac{1}{2}\left ( \sum_{i\in\Lambda ^f }\psi_{f_i}+
\sum_{k=2}^n(e_{F_1^k}^{\ast }+e_{F_2^k}^{\ast })\right ).$$
\noindent Then, for each $k$,
$$|f(\overline{x}_k|_{D_f})|=\frac{1}{2}
|\sum_{i\in\Lambda_k^f}f_i(\overline{x}_k|_{D_f})+
\sum_{i\in I_k^f}f_i(\overline{x}_k)+\sum_{i\in J_k^f}f_i(\overline{x}_k)|.$$

We have
$$\frac{1}{2}|\sum_{i\in\Lambda_k^f}f_i(\overline{x}_k|_{D_f})|\leq
\frac{1}{2}\sum_{i\in\Lambda_k^f}|f_i(\overline{x}_k|_{D_f})|
\leq\frac{1}{2}\sum_{i\in\Lambda_k^f}\psi_{f_i}(u_k)+
\sum_{i\in\Lambda_k^f}\sum_{F\in B_k(f_i)}C_F.$$

Also,
$$\frac{1}{2}|\sum_{i\in I_k^f}f_i(\overline{x}_k)|=
|h_{F_1^k}(\overline{x}_k)|e_{F_1^k}^{\ast }(e_{F_1^k})=
\frac{1}{2}e_{F_1^k}^{\ast }(a_{F_1^k}e_{F_1^k})=
\frac{1}{2}e_{F_1^k}^{\ast }(u_k),$$
\noindent and
$$\frac{1}{2}|\sum_{i\in J_k^f}f_i(\overline{x}_k)|=
|h_{F_2^k}(\overline{x}_k)|e_{F_2^k}^{\ast }(e_{F_k})=
\frac{1}{2}e_{F_2^k}(u_k).$$

We conclude that
$$|f(\overline{x}_k|_{D_f})|\leq\frac{1}{2}\left [\sum_{i\in\Lambda_k^f}
\psi_{f_i}(u_k)+e_{F_1^k}^{\ast }(u_k)+e_{F_2^k}^{\ast }(u_k)\right ]
+\sum_{F\in B_k(f)}C_F=\psi_f(u_k)+\sum_{F\in B_k(f)}C_F.$$

It remains to show that $\psi_f$ belongs to $K^{\prime }$. We need to
show that the family
$$B=\{ \psi_{f_i}:i\in\Lambda^f\}\cup\{ e_{F_1^k}^{\ast }: k=2,\ldots ,n\}
\cup\{ e_{F_2^k}^{\ast }:k=2,\ldots ,n\}$$
\noindent is ${\cal S}^{\prime}$--allowable.

We have ${\rm supp}\psi_{f_i}\subset D_{f_i}\subset {\rm supp}f_i$ for each
$i\in\Lambda^f$ and ${\rm supp}e_{F_1^k}=\{ n_{F_1^k}\}\subset
\cup\{ {\rm supp}f_i : f_i\in F_1^k\}\cap {\rm supp}\overline{x}_k$
 and the same 
is true for $e_{F^2_k}$.
 Also, if $f_i$ belongs to a family
 $F_2^k$, then $D_{f_i}\cap {\rm supp}\overline{x}_k=\emptyset $,
while 
$n_{F_2^k}\in {\rm supp}\overline{x}_k$.
Finally, we clearly have $e_{F_1^k}^{\ast }\neq e_{F_2^k}^{\ast }$.

The above remarks imply that the functionals in $B$ are disjointly
supported. Moreover, it is easy to see that
$$\# B\leq 2d=2(\#\{f_i:i=1,\ldots ,d\}).$$
\noindent We conclude that the family $B$ is ${\cal S}^{\prime}$--allowable,
and thus $\psi_f\in K^{\prime }$.

\medskip

This completes the inductive step. Of course, we set $\psi =\psi_{\varphi }$.

Then, $D_{\varphi }={\rm supp}\varphi\cap (\cup_{k=2}^n{\rm supp}\overline{x}_k
)$ (see Remark (1) following Definition 2.11), and by the
inductive assumption (b) we get: For each $k=2,\ldots ,n$,
$$|\varphi (\overline{x}_k)|\leq\psi (u_k)+\sum_{F\in B_{\overline{x}_k}}C_F.$$

To complete the proof of the Lemma it remains to show that,
for each $k=2,\ldots ,n$,
$$\sum_{F\in B_{\overline{x}_k}}C_F<\frac{2}{m_{j_k}}.$$
\noindent We have
$$\sum_{F\in B_{\overline{x}_k}}C_F=
\sum_{F\in B_{\overline{x}_k}}w(g_F)|\sum_{t\in L_F}a_{(k,t)}
\sum_{f\in F}f(\overline{x}_{(k,t)})|.$$
\noindent For each $F\in B_{\overline{x}_k}$, setting $x'_{(k,t)}=x_{(k,t)}|
{\rm supp}x'_k$ and $x^{\prime\prime}_{(k,t)}=x_{(k,t)}|{\rm supp}x^{\prime
\prime}_k$,  we have
$$|\sum_{t\in L_F}a_{(k,t)}\sum_{f\in F}f(\overline{x}_{(k,t)})|\leq
|\sum_{t\in L_F}a_{(k,t)}\sum_{f\in F}f(x_{(k,t)})|$$
$$+\sum_{t\in L_F}\sum_{f\in F}|f(a_{(k,t)}x^{\prime }_{(k,t)})|+
\sum_{t\in L_F}\sum_{f\in F}|f(a_{(k,t)}x^{\prime\prime }_{(k,t)})|.$$
\noindent Using Lemma 2.12 we get
$$|\sum_{t\in L_F}a_{(k,t)}\sum_{f\in F}f(\overline{x}_{(k,t)})|\leq
|\sum_{t\in L_F}a_{(k,t)}\sum_{f\in F}f(x_{(k,t)})|$$
$$+\sum_{t\in L_F}\sum_{f\in F}\frac{1}{m_{j_{k+1}}}
\|a_{(k,t)}x_{(k,t)}|_{A_{k,f}}\|_{\ell_1}+
\sum_{t\in L_F}\sum_{f\in F}\frac{1}{2^{n_k}}\|a_{(k,t)}x_{(k,t)}|_{A^{\prime }_{k,f}}\|_{\ell_1}$$
$$\leq |\sum_{t\in L_F}a_{(k,t)}\sum_{f\in F}f(x_{(k,t)})|+
\frac{1}{m_{j_{k+1}}}\sum_{f\in F}\|x_k|_{A_{k,f}}\|_{\ell_1}+
\frac{1}{2^{n_k}}\sum_{f\in F}\|x_k|_{A^{\prime }_{k,f}}\|_{\ell_1}.$$

To estimate
$$w(g_F)|\sum_{t\in L_F}a_{(k,t)}\sum_{f\in F}f(x_{(k,t)})|,$$
\noindent we use Remark (3) after 2.11. According to this Remark,
$w(g_F)>\frac{1}{m_{j_k}}$ and so, $g_F\in {\cal A}_r$
for some $1\leq r<j_k$. Let $g_F=w(g_F)\sum_{i=1}^lf_i$ where
$f_1<f_2<\ldots <f_l$ and suppose $i_1=\min \{ i:f_i\in F\} $ and
$i_2=\max \{ i:f_i\in F\}.$ We set $\tilde{F}=\{ f_i: i_1\leq i\leq i_2\}.$
The family $\tilde{F}$ contains $F$ but might also contain some functionals
$f_i$ with $f_i(x_k)\neq 0$ but $f_i(\overline{x}_k)=0.$ 
 Since $\tilde{K}$ is closed under projections onto intervals, the functional
$w(g_F)\sum_{f\in \tilde{F}}f$  belongs to ${\cal A}_r\cap \tilde{K}.$ 
Applying  Lemma 2.7 (a) (in fact, since our assumption is
$\| x_{(k,t)}\|_{\tilde{K}}\leq 1,$
 we use the analogue of this Lemma for
the space with norm $\| \cdot \|_{\tilde{K}}$)
we get that
$$w(g_F)|\sum_{t\in L_F}a_{(k,t)}\sum_{f\in \tilde{F}}f(x_{(k,t)})|\leq
\frac{1}{m_{j_k}^4}.$$
Notice that
$C_F:=w(g_F)|\sum_{t\in L_F}a_{(k,t)}\sum_{f\in F}f(\overline{x}_{(k,t)})|=
w(g_F)|\sum_{t\in L_F}a_{(k,t)}\sum_{f\in \tilde{F}}f(\overline{x}_{(k,t)})|$
and also that Lemma 2.12 remains true for $f\in \tilde{F}.$

We conclude that for each $F\in B_{\overline{x}_k}$,
$$C_F=w(g_F)|\sum_{t\in L_F}a_{(k,t)}\sum_{f\in \tilde{F}}f(\overline{x}_{(k,t)})|
\leq\frac{1}{m_{j_k}^4}+\frac{1}{m_{j_{k+1}}}\sum_{f\in \tilde{F}}
\| x_k|_{A_{k,f}}\|_{\ell_1}+
\frac{1}{2^{n_k}}\sum_{f\in \tilde{F}}\| x_k|_{A^{\prime }_{k,f}}\|_{\ell_1}.$$
\noindent Now, we add over all $F\in B_{\overline{x}_k}$. By Lemma
2.13, $\# B_{\overline{x}_k}<m_{j_k}$. Also, by Lemma 2.12 we
have that the sets $A_{k,f}$, $f\in\cup_{F\in B_{\overline{x}_k}}\tilde{F}$
are mutually disjoint, and the same is true for the sets
$A_{k,f}^{\prime }$. We conclude that
$$\sum_{F\in B_{\overline{x}_k}}C_F\leq\frac{m_{j_k}}{m_{j_k}^4}+
\frac{1}{m_{j_{k+1}}}\| x_k\|_{\ell_1}+\frac{1}{2^{n_k}}\| x_k\|_{\ell_1}$$
$$\leq\frac{1}{m_{j_k}^3}+\frac{1}{m_{j_{k+1}-1}}+
\frac{1}{m_{j_k}}<\frac{2}{m_{j_k}}$$
\noindent by Definition 2.8. This completes the proof of the Lemma.   $\Box $

\bigskip

%\newpage

\noindent {\bf Proof of Proposition 2.9.}

\medskip

Recall (Definition 2.11) that for our intermediate lemmas we have
assumed that ${\rm supp}\varphi\cap {\rm supp}\overline{x}_1\neq
\emptyset $. If this is not true, then we can set
$k_0=\min\{ k: {\rm supp}\varphi\cap {\rm supp}\overline{x}_k\neq
\emptyset\}$ and construct in the same way $u_k$'s,
$k=k_0+1,\ldots ,n$, and $\psi $ supported on
$\cup_{k=k_0+1}^n{\rm supp}u_k$, such that
$$|\varphi (\overline{x}_k)|\leq\psi (u_k)+\frac{2}{m_{j_k}},\;\;
k=k_0+1,\ldots ,n.$$
\noindent Setting $u_k=0$, for $k=2,\ldots ,k_0$ we have
$$|\varphi (\sum_{k=1}^n\lambda_k\overline{x}_k)|\leq
|\lambda_{k_0}|\;|\varphi (\overline{x}_{k_0})|+
\psi (\sum_{k=2}^n|\lambda_k|u_k)+\sum_{k=2}^n|\lambda_k|\frac{2}{m_{j_k}}$$
\noindent for any choice of coefficients $(\lambda_k)_{k=1}^n$.

For $\varphi (\sum_{k=1}^n\lambda_kx_k)$ we have
$$|\varphi (\sum_{k=1}^n\lambda_kx_k)|\leq
|\varphi (\sum_{k=1}^n\lambda_k\overline{x}_k)|+
\sum_{k=1}^n|\lambda_k|\;( |\varphi (x^{\prime }_k)|+
|\varphi (x^{\prime\prime }_k)|).$$
\noindent Using the previous estimate and Lemma 2.10 we get
$$|\varphi (\sum_{k=1}^n\lambda_kx_k)|\leq |\lambda_{k_0}|
|\varphi (\overline{x}_{k_0})|+\psi (\sum_{k=2}^n|\lambda_k|u_k)
+4\sum_{k=1}^n|\lambda_k|\frac{1}{m_{j_k}}$$
$$\leq |\lambda_{k_0}|\;(|\varphi (x_{k_0})|+
|\varphi (x^{\prime }_{k_0})|+|\varphi (x^{\prime\prime }_{k_0})|)
+\psi (\sum_{k=2}^n|\lambda_k|u_k)+4\sum_{k=1}^n|\lambda_k|
\frac{1}{m_{j_k}}$$
$$\leq |\lambda_{k_0}|\;\|x_{k_0}\|_{\tilde{K}}
+\psi (\sum_{k=2}^n|\lambda_k|u_k)+
6\sum_{k=1}^n|\lambda_k|\frac{1}{m_{j_k}}$$
$$\leq\max_{1\leq k\leq n}|\lambda_k|+\psi (\sum_{k=2}^n|\lambda_k|u_k)
+6\sum_{k=1}^n|\lambda_k|\frac{1}{m_{j_k}}.\;\;\;\Box $$

\bigskip

\noindent {\bf 2.16 Definition.} Let $j\geq 2$, $\varepsilon >0$.
An $(\varepsilon ,j)$-special convex combination
$\sum_{k=1}^nb_kx_k$ is called an $(\varepsilon ,j)$-R.I.s.c.c.
w.r.t. $\tilde{K}$ if the sequence $(x_k)_{k=1}^n$ is a R.I.S.
w.r.t. $\tilde{K}$ and the corresponding integers
$(j_k)_{k=1}^n$ satisfy
$j+2<j_1<\ldots <j_n.$

\medskip

\noindent {\bf 2.17 Corollary.} {\it If $\sum_{k=1}^nb_kx_k$ is
a $(\frac{1}{m_j^2},j)$--R.I.s.c.c. w.r.t. $\tilde{K}$ and $\varphi\in
\tilde{K}$ with $w(\varphi )=\frac{1}{m_s}$, then}
$$|\varphi (\sum_{k=1}^nb_kx_k)|\leq 2b_1+\frac{16}{m_s}\;\;\;,\;\;\;
{\it if}\;s\geq j. \leqno {\rm (a)}$$
$$|\varphi (\sum_{k=1}^nb_kx_k)|\leq \frac{33}{m_sm_j}\;\;\;,
\;\;\;{\it if}\;s<j.$$

$$\frac{1}{4m_j}\leq\|\sum_{k=1}^nb_kx_k\|_{\tilde{K}}\leq\frac{17}{m_j}.
\leqno {\rm (b)}$$

\medskip

\noindent {\it Proof:} (a) Recall that the sequence $(b_k)_{k=1}^n$
is decreasing. By Proposition 2.9,
$$|\varphi (\sum_{k=1}^nb_kx_k)|\leq b_1+
\psi (\sum_{k=2}^nb_ku_k)+6\sum_{k=1}^n\frac{b_k}{m_{j_k}},$$
\noindent where $\psi\in K^{\prime }$ with $w(\psi )=w(\varphi )=s$
and $\| u_k\|_{\ell_1}\leq 16$. By Lemma 2.4 we get
$$|\varphi (\sum_{k=1}^nb_kx_k)|\leq 2b_1+\frac{16}{m_s}$$
\noindent for $s\geq j$, and
$$|\varphi (\sum_{k=1}^nb_kx_k)|\leq 2b_1+\frac{32}{m_sm_j}<
\frac{33}{m_sm_j}$$
\noindent for $s<j$.

(b) The upper estimate  follows from (a). The lower estimate
is a consequence of the fact that $\| x_k\|_{\tilde{K}}\geq\frac{1}{2}$
and the sequence $(x_k)_{k=2}^n$ is ${\cal M}_j$--admissible.   $\Box $

\bigskip

\noindent {\bf 2.18 Theorem.} {\it The space $X_{M(1),u}$ is
arbitrarily distortable.}

\medskip

\noindent {\it Proof:} It follows from Lemmas 2.2 and 2.6 that for every $j\geq 2$
every block subspace $Y$ contains a $(\frac{1}{m_j^2},j)$--R.I.s.c.c. w.r.t.
$K$.

Fix $i_0\in {\bf N}$ large
 and define an equivalent norm $|\|.\||$ on $X_{M(1),u}$
by
$$|\|x\||=\frac{1}{m_{i_0}}\|x\|+\sup\{ \varphi (x):
\varphi\in {\cal A}_{i_0}\}.$$
\noindent Let $Y$ be a block subspace and let $y=\sum a_ky_k\in Y$
be a $(\frac{1}{m_j^2},j)$--R.I.s.c.c.
for some $j>i_0$, and $z=\sum b_lz_l\in Y$ be a
$(\frac{1}{m_{i_0}^2},i_0)$--R.I.s.c.c. Then, by Corollary 2.16,
$$\||m_jy\||\leq\frac{17}{m_{i_0}}+\frac{33}{m_{i_0}}=\frac{50}{m_{i_0}}
\;\;\;\;{\rm and}\;\;\;\;\|m_jy\|\geq\frac{1}{4}.$$

On the other hand,
$$\||m_{i_0}z\||\geq\frac{1}{4}\;\;\;\;{\rm and}\;\;\;\;\|m_{i_0}z\|\leq
17.$$

This shows that $\||.\||$ is a $\frac{1}{10^3}m_{i_0}$--distortion.
Since $i_0$ was arbitrary, this completes the proof.   $\Box $ 

\bigskip

The following Remarks on the proof of Proposition 2.9 will be used in the
next Section.

\medskip

\noindent {\bf 2.19 Remark:}
Let $\varphi$, $\overline{x}_k$, $\psi$, 
$u_k$ be as in Proposition 2.9.  It follows from the proof of Lemma 2.15 that
the functional $\psi$ which is constructed inductively folowing the analysis
$\{K^s(\varphi )\}$ 
of $\varphi$ satisfies the following properties.

(a) There exists an analysis $\{ K^s(\psi )\}$ of $\psi$ contained in $K'$
such that, for every $g\in\cup K^s(\psi )$ there exists a unique
$f\in\cup K^s(\varphi )$ with $g=\psi_f$; moreover, if $g\not\in K^0$
then $w(f)=w(g)$.

(b) The functional $\psi$ is supported in the set 
$$L=\left\{ e_f: f\in \cup\{ F: F \;
{\rm is\; a\; family \; of\; type \; I \; or\;
II\; w.r.t. \; some}\; \overline{x}_k\}\right\}.$$
Moreover, for $k=2,\ldots , n$ and for every family $F$ of type I or II w.r.t.
$\overline{x}_k$, if we set $V_F=\cup_{f\in F}{\rm supp}f$ and
$W_F=\{ e_f:f\in F\}$ we have 
$$|\varphi |_{V_F}(\overline{x}_k)|\leq \psi |_{W_F}(u_k)+C_F$$
where we have set $C_F=0$ if $F\not\in B_{\overline{x}_k}$.

(c) Let $\varphi_2=\varphi |J$ for some $J\subset {\bf N}$. Assume further
that $\varphi_2$ has the following property:

{\it For every $k=2,\ldots ,n$ and every family $F=\{ f_1,\ldots,f_l\}\subset
\cup K^s(\varphi)$ 
of type I or II w.r.t. $\overline{x}_k$, either $f_i|_J(\overline{x}_k)=0$
for all $i=1,\ldots,l$ or $f_i|_{J}(\overline{x}_k)=f_i(\overline{x}_k)$
for all $i=1,\ldots ,l$.} 

For $k=2,\ldots ,n$, we let 
$$L_k=\{ e_f: f \; {\rm belongs\; to \; some \; family\; of \; type\; I\;
or\; II\; w.r.t. }\; \overline{x}_k\; {\rm and}\; {\rm supp}f\cap J\neq
\emptyset\}$$
and we set $\psi_2=\psi |\cup_{k=2}^n L_k.$ Then it follows from (b) that
$$|\varphi_2(\overline{x}_k)|\leq \psi_2(u_k)+\frac{1}{m_{j_k}},\;\;\;
k=2,\ldots, n.$$

\bigskip

\bigskip

\newpage

\centerline { {\sc 3. The space $X$.}}

\medskip

We pass now to the construction of a space $X$ not containing
any unconditional basic sequence. It is based on the modification
$X_{M(1),u}$. Let $K=\cup_n\cup_jK^n_j$ be
 the norming set of the space $X_{M(1),u}$.
Consider the countable set
$$G=\{ (x_1^{\ast }, x_2^{\ast },\ldots , x_k^{\ast }):
k\in {\bf N}, x_i^{\ast }\in K, i=1,\ldots ,k\;{\rm and}
\;x_1^{\ast }<x_2^{\ast }<\ldots <x_k^{\ast }\}.$$

There exists a one to one function $\Phi :G\rightarrow
\{ 2j\}_{j=1}^{\infty }$ such that for every
$(x_1^{\ast },\ldots ,x_k^{\ast })\in G$, if $j_1$
is minimal such that $x_1^{\ast }\in {\cal A}_{j_1}$
and $j_l=\Phi (x_1^{\ast },\ldots ,x_{l-1}^{\ast })$, $l=2,3,\ldots ,k$,
then
$$\Phi (x_1^{\ast },\ldots ,x_k^{\ast })>
\max\{ j_1,\ldots ,j_k\} .$$

\noindent {\bf Definition of the space $X$.}
For $n=0,1,2,\ldots ,$ we define by induction sets
$\{ L_j^n\}_{j=1}^{\infty }$ such that $L_j^n$ is a subset
of $K_j^n$.

For $j=1,2,\ldots $, we set $L_j^0=\{ \pm e_n:n\in {\bf N}\}$.
Suppose that $\{ L_j^n\}_{j=1}^{\infty }$ have been defined. 
We set $L^n=\cup_{j=1}^{\infty}L^n_j$ and
$$L_1^{n+1}=\pm L_1^n\cup\{\frac{1}{2}(x_1^{\ast }+\ldots +x_d^{\ast }):
d\in {\bf N},x_i^{\ast }\in L^n,$$
$$d\leq\min {\rm supp}x_1^{\ast }<\ldots <\min {\rm supp}x_d^{\ast },
\;{\rm supp}x_i^{\ast }\cap {\rm supp}x_l^{\ast }=\emptyset \;
{\rm for}\; i\neq l\},$$
\noindent and for $j\geq 1$,
$$L_{2j}^{n+1}=\pm L_{2j}^n\cup\{\frac{1}{m_{2j}}(x_1^{\ast }+\ldots
+x_d^{\ast }):d\in {\bf N},x_i^{\ast }\in L^n,$$
$$({\rm supp}x_1^{\ast },\ldots ,{\rm supp}x_d^{\ast })\;{\rm is}\;
{\cal M}_{2j}-{\rm admissible }\},$$

$$L_{2j+1}^{\prime\;n+1}=\pm L_{2j+1}^n\cup
\{\frac{1}{m_{2j+1}}(x_1^{\ast }+\ldots +x_d^{\ast }):d\in {\bf N},
x_1^{\ast }\in L_{2k}^n\;{\rm for}\;{\rm some}\;k>2j+1,$$
$$x_i^{\ast }\in L_{\Phi (x_1^{\ast },\ldots ,x_{i-1}^{\ast })}^n
\;{\rm for}\;1<i\leq d\;{\rm and}\;({\rm supp}x_1^{\ast },\ldots ,
{\rm supp}x_d^{\ast })\;{\rm is}\;{\cal M}_{2j+1}-{\rm admissible}\},$$

$$L_{2j+1}^{n+1}=\{ E_sx^{\ast }:x^{\ast }\in L_{2j+1}^{\prime\;n+1},
s\in {\bf N}, E_s=\{ s,s+1,\ldots \}\}.$$

\noindent This completes the definition of $L_j^n$, $n=0,1,2,\ldots $,
$j=1,2,\ldots $ It is obvious that each $L_j^n$ is a subset of the
corresponding set $K_j^n$.

We set ${\cal B}_j=\cup_{n=1}^{\infty }(L_j^n\setminus L^0)$
 and we consider the norm on $c_{00}$
defined by the set $L=L^0\cup (\cup_{j=1}^{\infty }{\cal B}_j)$.
The space $X$ is the completion of $c_{00}$ under this norm.
It is easy to see that $\{ e_n\}_{n=1}^{\infty }$ is a
bimonotone basis for $X$.

\medskip

\noindent {\bf Remark:} The norming set $L$ is closed under 
projections onto {\it intervals}, and has the property that for every $j$
and every ${\cal M}_{2j}$--admissible family
$f_1, f_2, \ldots f_d$ contained in $L$, $\frac{1}{m_{2j}}(f_1+\cdots +f_d)$
belongs to $L$. It follows that for every $j=1,2,\ldots$ and every
 ${\cal M}_{2j}$--admissible family $x_1<x_2<\ldots<x_n$ in $c_{00}$, 
$$\|\sum_{k=1}^nx_k\|\geq\frac{1}{m_{2j}}\sum_{k=1}^n\| x_k\|.$$
For the same reason, for ${\cal S}$--admissible families $x_1<x_2<\ldots
<x_n$, we have
$$\|\sum_{k=1}^nx_k\|\geq\frac{1}{2}\sum_{k=1}^n\| x_k\|.$$
We note however that such a relation is {\it not} true for ${\cal S}$--{\it 
allowable} families $(x_i)$. Of course, if it were true, it would 
immediately imply that the basis $\{ e_n\}$ is unconditional.

\medskip

For $\varepsilon >0$, $j=2,\ldots $, $(\varepsilon ,j)$--special
convex combinations are defined in $X$ exactly as in
 $X_{M(1),u}$ (Definition 2.5). Rapidly increasing sequences and
$(\varepsilon , j)$--R.I. special convex combinations in $X$ are
defined by Definitions 2.8 and 2.16 respectively, with $\tilde{K}=L.$

\medskip

By the previous Remark we get the following.

\medskip

\noindent {\bf 3.1 Lemma.} {\it For $j=2,3,\ldots $ and every
normalized block sequence $\{ x_k\}_{k=1}^{\infty }$ in $X$,
there exists a finite normalized block sequence $\{ y_s\}_{s=1}^n$
of $\{ x_k\}$ such that $\sum_{s=1}^na_sy_s$ is a
seminormalized $(\frac{1}{m_{j}^4},j)$--s.c.c.   $\Box $}

\medskip

By  Corollary 2.17, we have:

\medskip

\noindent {\bf 3.2 Proposition.} {\it Let $\sum_{k=1}^nb_ky_k$
be a $(\frac{1}{m_j^2},j)$--R.I.s.c.c. in $X$. Then, for $i\in {\bf N}$,
$\varphi\in {\cal B}_i$, we have the following:

$$|\varphi (\sum_{k=1}^nb_ky_k)|\leq\frac{33}{m_im_j}\;\;\;,\;\;\;
{\rm if}\;i<j\leqno {\rm (a)}$$
$$|\varphi (\sum_{k=1}^nb_ky_k)|\leq\frac{16}{m_i}+2b_1\;\;\;,
\;\;\;{\rm if}\;i\geq j.\leqno {\rm (b)}$$
\noindent In particular, $\|\sum_{k=1}^nb_ky_k\|\leq\frac{17}{m_j}$.
   $\Box $}

\medskip

\noindent {\bf 3.3 Proposition.} {\it Let $j>100$. Suppose that
$\{ j_k\}_{k=1}^n$, $\{ y_k\}_{k=1}^n$, $\{ y_k^{\ast }\}_{k=1}^n$
and $\{\theta_k\}_{k=1}^n$ are such that

\medskip

{\rm (i)} There exists a rapidly increasing sequence (w.r.t. $X$)
$$\{ x_{(k,i)}:\; k=1,\ldots ,n,\; i=1,\ldots ,n_k\} $$ with 
$x_{(k,i)}<x_{(k,i+1)}<x_{(k+1,l)}$ for all $k<n$, $i<n_k$, $l\leq n_{k+1},$
such that:  

\noindent {\rm (a)} 
Each $x_{(k,i)}$ is a seminormalized $(\frac{1}{m^4_{j_{(k,i)}}},
j_{(k,i)})$--s.c.c. where, for each $k$, $2j_k+2<j_{(k,i)},\; i=1,\ldots n_k.$

\noindent {\rm (b)} Each $y_k$ is a $(\frac{1}{m^4_{2j_k}},2j_k)$--
 R.I.s.c.c. of  $\{ x_{(k,i)}\}_{i=1}^{n_k}$ of the form
$y_k=\sum _{i=1}^{n_k}b_{(k,i)}x_{(k,i)}.$

\noindent {\rm (c)} There exists a decreasing
sequence $\{ b_k\}_{k=1}^n$ such that $\sum _{k=1}^nb_ky_k$ is a $(\frac{1}
{m^4_{2j+1}}, 2j+1)$--s.c.c.

\medskip

{\rm (ii)} $y_k^{\ast }\in L_{2j_k}$, $y_k^{\ast }(y_k)\geq
\frac{1}{4m_{2j_k}}$ and ${\rm supp}y_k^{\ast }\subset
[\min {\rm supp}y_k,\max {\rm supp}y_k]$.

\medskip

{\rm (iii)} $\frac{1}{17}\leq\theta_k\leq 4$ and
$y_k^{\ast }(m_{2j_k}\theta_ky_k)=1$.

\medskip

{\rm (iv)} $j_1>2j+1$ and $2j_k=\Phi (y_1^{\ast },\ldots ,y_{k-1}^{\ast })$,
$k=2,\ldots ,n$.

\medskip

Let  $\varepsilon_k=(-1)^{k+1}$, $k=1,\ldots n.$
Then,}
$$\|\sum_{k=1}^n\varepsilon_kb_k\theta_km_{2j_k}y_k\|\leq
\frac{300}{m_{2j+1}^2}.$$

\medskip

Before presenting the proof of Proposition 3.3 let us show how
from it the main result of this section follows.

\medskip

\noindent {\bf 3.4 Corollary.} {\it 
The space $X$ is Hereditarily Indecomposable.}

\medskip

\noindent {\it Proof:} It is clear by the choice of the sequences
$\{ y_k\}_{k=1}^n$, $\{ y^*_k\}_{k=1}^n$ in Proposition 3.3 that 
the functional $\psi =\frac{1}{m_{2j+1}}\sum_{k=1}^ny^*_k $ belongs 
to $L$ and that $\psi (\sum_{k=1}^n b_km_{2j_k}\theta_ky_k)=\frac{1}{m_{2j+1}}.$
It follows that 
$$\| \sum_{k=1}^nb_km_{2j_k}\theta_ky_k\|\geq \frac{1}{m_{2j+1}}.$$
To conclude that $X$ is Hereditarily Indecomposable it remains to
show that, for every $j>100$ and every block subspaces $U$ and $V$ of
$X$, one can choose $\{ y_k\}$ and $\{ y^*_k\}$ satisfying the assumptions 
of Proposition 3.3 and such that $y_k\in U$ if $k$ is odd, $y_k\in V$
if $k$ is even. The proof of this is the same as that of Proposition
3.12 [3], so we omit it. $\Box$

\medskip

\noindent {\bf Proof of Proposition 3.3.}

  Our aim is to show that
for every $\varphi\in\cup_{i=1}^{\infty }{\cal B}_i$,
$$\varphi (\sum_{k=1}^n\varepsilon_kb_k\theta_km_{2j_k}y_k)\leq
\frac{300}{m_{2j+1}^2}.$$
The proof is given in several steps. We give a brief description:

We consider separately three cases for $\varphi$:

{\tt 1st Case}: $w(\varphi)=\frac{1}{m_{2j+1}}$. Then $\varphi$ has the
form $\varphi =\frac{1}{m_{2j+1}}(Ey^*_{k_1-1}+y^*_{k_1}+\cdots
+y^*_{k_2}+z^*_{k_2+1}+\cdots z^*_d)$ and for the part $\frac{1}{m_{2j+1}}
(y^*_{k_1}+\cdots +y^*_{k_2})$ acting on 
$\sum _{k=k_1}^{k_2}\varepsilon_k
b_k\theta_k m_{2j_k}y_k$ 
we have an obvious conditional (i.e. depending
on the signs) estimate. For the remainining part we get an unconditional
estimate using Proposition 3.2 (Lemmas 3.5, 3.6).

{\tt 2nd Case}: $w(\varphi )\leq \frac{1}{m_{2j+2}}.$ Then we get an
unconditional estimate for 
$\varphi (\sum _{k=1}^{n}\varepsilon_k
b_k\theta_k m_{2j_k}y_k)$ directly, applying Proposition 3.2 (Lemma 3.7).

{\tt 3rd Case}: $w(\varphi)>\frac{1}{m_{2j+1}}$. We
fix an analysis $\{ K^s(\varphi )\}$ of $\varphi$. By Proposition 2.9
we get that there exist $\psi\in co(K')$ and a block sequence 
$u_k=\sum_{i=1}^{n_k}b_{(k,i)}u_{(k,i)}$, $k=1,\ldots ,n$
 of subconvex combinations
of the basis with 
$\varphi (\sum _{k=1}^{n}\varepsilon_k
b_k\theta_k m_{2j_k}y_k)\leq 
\psi(\sum _{k=1}^{n}\varepsilon_k
b_k\theta_k m_{2j_k}u_k) +\frac{1}{m_{2j+2}}$. However, since the 
estimate that we get in this way is unconditional, it is insufficient.
So, we partition $\varphi$ into two disjointly supported functionals
$\varphi_1$ and $\varphi_2$, where $\varphi_1$ is the restriction of
$\varphi$ which contains in its analysis certain projections of
the functionals of the  
form $f =\frac{1}{m_{2j+1}}(Ey^*_{k_1-1}+y^*_{k_1}+\cdots
+y^*_{k_2}+z^*_{k_2+1}+\cdots z^*_d)$ in $\cup K^s(\varphi ).$
For $\varphi_1 (\sum _{k=1}^{n}\varepsilon_k
b_k\theta_k m_{2j_k}y_k)$ we give  a conditional estimate (Lemma
3.12(b)). To get an estimate for 
 $\varphi_2 (\sum _{k=1}^{n}\varepsilon_k
b_k\theta_k m_{2j_k}y_k)$  we show that 
 $\varphi_2 (\sum _{k=1}^{n}\varepsilon_k
b_k\theta_k m_{2j_k}y_k)$  is dominated by 
 $\psi_2 (\sum _{k=1}^{n}\varepsilon_k
b_k\theta_k m_{2j_k}y_k)$  where $\psi_2$ is the restriction of $\psi$
corresponding to $\varphi_2$ (Lemma 3.10). Then we estimate the action
of $\psi_2$ on  
 $\sum _{k=1}^{n}\varepsilon_k
b_k\theta_k m_{2j_k}y_k$ (Lemma 3.11(a)). 

\medskip

\noindent {\bf 3.5 Lemma.} {\it Let $j,\{j_k\}_{k=1}^n$ and $\{ y_k\}_{k=1}^n$
be as in Proposition 3.3. Suppose that $2j+1<t_1<\ldots <t_d$ and
let $\{ z_s^{\ast }\}_{s=1}^d$ be such that $z_1^{\ast }<\ldots <
z_d^{\ast }$, $z_s^{\ast }\in {\cal B}_{2t_s}$ and
$\frac{1}{m_{2j+1}}(z_1^{\ast }+\ldots +z_d^{\ast })\in {\cal B}_{2j+1}$.
Assume that for some $k=1,2,\ldots ,n$, $j_k\notin\{ t_1,\ldots ,t_d\}$.
Then,}
$$|(\sum_{s=1}^dz_s^{\ast })(m_{2j_k}y_k)|\leq
\frac{1}{m_{2j+2}^2}.$$

\medskip

\noindent {\it Proof:} Each $y_k$ is a $(\frac{1}{m_{2j_k}^4},2j_k)$--R.I.s.c.c. of the form $y_k=\sum_{i=1}^{n_k}b_{(k,i)}x_{(k,i)}$.
Let $s_1\leq d$ be such that $s_1=\max\{ s\in\{1,\ldots ,d\}:
t_{s}<j_k\}$.

If $s\leq s_1$, by Proposition 3.2(a) we get $|z_s^{\ast }(y_k)|\leq
\frac{33}{m_{2t_s}m_{2j_k}}$ and so, using that $2j+1<t_1<\ldots <t_d$
and that the sequence $\{ m_j\}$ is increasing sufficiently fast, we get
$$|(\sum_{s=1}^{s_1}z_s^{\ast })(y_k)|\leq
\frac{33}{m_{2j_k}}\sum_{s=1}^{s_1}\frac{1}{m_{2t_s}}\leq
\frac{1}{2m_{2j+2}^2m_{2j_k}}.\leqno (\ast )$$
\noindent For every $s\geq s_1+1$ set
$$D_s=\{ i: {\rm supp}x_{(k,i)}\cap {\rm supp}z_s^{\ast }=
{\rm supp}x_{(k,i)}\cap {\rm supp}\sum_{t=s_1+1}^dz_t^{\ast }\}.$$
\noindent The sets $D_s$ are disjoint. Put $I=\{ s\geq s_1+1:
D_s\neq\emptyset\}$ and 
$$T=\{ r:1\leq r\leq n_k, {\rm supp}x_{(k,r)}\cap
{\rm supp}\sum_{t=s_1+1}^dz_t^{\ast }\neq\emptyset\}\backslash
\cup_{s\in I}D_s.$$
 Then,
$$|(\sum_{s=s_1+1}^dz_s^{\ast })(y_k)|\leq
\sum_{s\in I}|z_s^{\ast }(\sum_{r\in D_s}b_{(k,r)}x_{(k,r)})|+
|\sum_{s=s_1+1}^d z_s^{\ast }
(\sum_{r\in T}b_{(k,r)}
x_{(k,r)})|
.\leqno (1)$$
\noindent It follows from Proposition 3.2(b) that for every $s\in I$,
$$|z_s^{\ast }(\sum_{r\in D_s}b_{(k,r)}x_{(k,r)})|\leq
\frac{16}{m_{2t_s}}+2b_{(k,r_s)},\leqno (2)$$
\noindent where $r_s=\min D_s$. Since by the definition of $D_s$ we
have that $\{\max{\rm supp}x_{(k,r_s)}\}_{s\in I}\in {\cal M}_{2j+1}$,
then
$$\sum_{s\in I}b_{(k,r_s)}\leq\frac{1}{m_{2j_k}^4}.\leqno (3)$$

Since $\frac{1}{m_{2j+1}}(z_1^{\ast }+\ldots +z_d^{\ast })\in {\cal B}_{2j+1}$,
as in Lemma 2.7(a) we have
$$
|(\sum_{s=s_1+1}^dz_s^{\ast })(\sum _{r\in T}b_{(k,r)}x_{(k,r)})|
(x_{(k,r)})|
\leq\frac{m_{2j+1}}{m_{2j_k}^4}<\frac{1}{m^3_{2j_k}}.\leqno (4)$$
\noindent By (1), (2), (3), (4), using that $j_k<t_{s_1+1}$,
 $m_{i+1}\geq m^i_i$ and that $2j+2<2j_1$,
we have that
$$|(\sum_{s=s_1+1}^dz_s^{\ast })(y_k)|\leq
16\sum_{s=s_1+1}^d\frac{1}{m_{2t_s}}+\frac{2}{m_{2j_k}^4}+\frac{1}{m_{2j_k}^3}
\leq
\frac{1}{2m_{2j+2}^2m_{2j_k}}.\leqno (\ast\ast )$$
\noindent Therefore, by ($\ast $) and ($\ast\ast $), we get
$$|(\sum_{s=1}^dz_s^{\ast })(m_{2j_k}y_k)|\leq\frac{1}{m_{2j+2}^2}.$$

\medskip

\noindent {\bf 3.6 Lemma.} {\it Let $j,\{j_k\}_{k=1}^n,\{ y_k\}_{k=1}^n,
\{ y_k^{\ast }\}_{k=1}^n,\{\theta_k\}_{k=1}^n$
and $\{\varepsilon_k\}_{k=1}^n$ be as in Proposition 3.3.
For every $\varphi\in {\cal B}_{2j+1}$ we have}
$$|\varphi (\sum_{k=1}^n\varepsilon_kb_k\theta_km_{2j_k}y_k)|\leq
\frac{1}{m_{2j+1}^2}.$$

\medskip

\noindent {\it Proof:} Let $\varphi =\frac{1}{m_{2j+1}}(Ey_{k_1}^{\ast }+
y_{k_1+1}^{\ast }+\ldots +y_{k_2}^{\ast }+z_{k_2+1}^{\ast }+
\ldots +z_d^{\ast })$, where $E=E_s$ for some $s$ and
$z_{k_2+1}^{\ast }\neq y_{k_2+1}^{\ast }$.

For $k=1,2,\ldots ,n$ we set $z_k=\theta_km_{2j_k}y_k$, hence
$y_k^{\ast }(z_k)=1$. Since $\{b_k\}$ is decreasing,
$$|\varphi (\sum_{k=k_1+1}^{k_2-1}\varepsilon_kb_kz_k)|
\leq\frac{1}{m_{2j+1}}|\sum_{k=k_1+1}^{k_2-1}\varepsilon_kb_k
y_k^{\ast }(z_k)|\leqno ({\rm a})$$
$$=\frac{1}{m_{2j+1}}|\sum_{k=k_1+1}^{k_2-1}\varepsilon_kb_k|\leq
\frac{1}{m_{2j+1}}b_{k_1+1},$$
\noindent and
$$|\varphi (z_{k_1})|=\frac{1}{m_{2j+1}}|Ey_{k_1}^{\ast }(z_{k_1})|
\leq\frac{1}{m_{2j+1}}\| z_{k_1}\|\leq\frac{68}{m_{2j+1}}.\leqno ({\rm b})$$
\noindent For $z_{k_2}$ we have
$$|\varphi (z_{k_2})|\leq\frac{1}{m_{2j+1}}|y_{k_2}^{\ast }(z_{k_2})|
+\frac{1}{m_{2j+1}}|(\sum_{k=k_2+1}^dz_k^{\ast })(z_{k_2})|.$$
\noindent If $k\geq k_2+1$, then $z_k^{\ast }\in B_{2t_k}$ where
$2t_k=\Phi (y_1^{\ast },\ldots ,y_{k_1}^{\ast },\ldots ,z_{k-1}^{\ast })$.
Since $\Phi $ is one to one, $2t_k\neq\Phi (y_1^{\ast },\ldots ,
y_{k_2-1}^{\ast })=2j_{k_2}$. Thus, by Lemma 3.5,
$$\frac{1}{m_{2j+1}}|\sum_{k=k_2+1}^dz_k^{\ast }(z_{k_2})|\leq
\frac{1}{m_{2j+1}}\frac{\theta_k}{m_{2j+2}^2}<\frac{1}{m_{2j+1}},$$
\noindent and so,
$$|\varphi (z_{k_2})|\leq\frac{2}{m_{2j+1}}.\leqno ({\rm c})$$
In a similar way, for $z_{k_2+1}$ we have
$$|\varphi (z_{k_2+1})|\leq \frac{1}{m_{2j+1}}|z^*_{k_2+1}(z_{k_2+1})|+
\frac{1}{m_{2j+1}}|(\sum_{k>k_2+1}z^*_k)(z_{k_2+1})|
<\frac{69}{m_{2j+1}}. \leqno ({\rm d})$$
\noindent If $k<k_1$, then $\varphi (z_k)=0$. By Lemma 3.5, for
$k>k_2+1$ we have
$$|\varphi (z_k)|=\frac{1}{m_{2j+1}}|\sum_{p=k_2+1}^dz_p^{\ast }(z_k)|
\leq\frac{1}{m_{2j+1}}\frac{\theta_k}{m_{2j+2}^2}<\frac{1}{m_{2j+2}^2}.
\leqno ({\rm e})$$
\noindent Putting (a), (b), (c),  (d) and (e) together and using
that, since $\sum b_ky_k$ is a $(\frac{1}{m^4_{2j+1}},2j+1)
$--s.c.c., $b_k<\frac{1}{m^4_{2j+1}}$,  we get the result.

\medskip

\noindent {\bf 3.7 Lemma.} {\it Under the assumptions of
Proposition 3.3, let $\varphi\in {\cal B}_r$ for $r\geq 2j+2$. Then,}
$$|\varphi (\sum_{k=1}^n\varepsilon_kb_k\theta_km_{2j_k}y_k)|
\leq\frac{134}{m_{2j+2}}.$$

\medskip

\noindent {\it Proof:} If $2j+1<r<2j_1$, it follows from
Proposition 3.2(a). The case $2j_{k_0}\leq r<2j_{k_0+1}$ follows
from Proposition 3.2(a), (b), and the lacunarity of the sequence
$\{ m_j\}_{j=1}^{\infty }$. The case $r>2j_n$ is similar.   $\Box $

\medskip

\noindent {\bf 3.8 Proposition.} {\it Let $j,\{j_k\}_{k=1}^n,\{ y_k\}_{k=1}^n,
\{ y_k^{\ast }\}_{k=1}^n,\{\theta_k\}_{k=1}^n,\{\varepsilon_k\}_{k=1}^n$
be as in Proposition 3.3. For every $\varphi\in {\cal B}_r$,
$r<2j+1$, we have}
$$|\varphi (\sum_{k=1}^n\varepsilon_kb_k\theta_km_{2j_k}y_k)|\leq
\frac{262}{m_{2j+1}^2}.$$

\medskip

\noindent The proof is based on
Proposition 2.9. We first need
to introduce new notation and establish several
Lemmas. We have $y_k=\sum_{i=1}^{n_k}b_{(k,i)}x_{(k,i)}$
and the sequence $\{ x_{(k,i)},k=1,\ldots n,i=1,\ldots n_k\}$ is a R.I.S.
w.r.t. $L$.
By Proposition 2.9 there exist a functional $\psi\in K^{\prime }$ and
blocks of the basis $u_{(k,i)}$, $k=1,\ldots ,n$, $i=1,\ldots ,n_k$ with
 $\psi\in {\cal A}_r^{\prime }$,
 ${\rm supp}u_{(k,i)}\subset {\rm supp}x_{(k,i)}$,
$\| u_k\|_{\ell_1}\leq 16$ and such that
$$|\varphi (\sum_{k=1}^n\varepsilon_kb_k\theta_km_{2j_k}
(\sum_{i=1}^{n_k}b_{(k,i)}x_{(k,i)}))|\leq
\theta_1m_{2j_1}b_1b_{(1,1)}+\psi (\sum_{k=1}^nb_k\theta_k
m_{2j_k}(\sum_{i=1}^{k_n}b_{(k,i)}u_{(k,i)}))+\frac{1}{m_{2j+2}^2}$$
$$\leq\psi (\sum_{k=1}^nb_k\theta_km_{2j_k}(\sum_{i=1}^{k_n}
b_{(k,i)}u_{(k,i)}))+\frac{1}{m_{2j+2}}.$$

Recall that the construction of $\psi $ and $u_{(k,i)}$ is done
via some analysis $\{ K^s(\varphi )\}$ of $\varphi $ and some
restriction on the support of $x_{(k,i)}$ which we denote by
$\overline{x}_{(k,i)}$. Let $\{ K^s(\varphi )\}$ be the analysis
of $\varphi $ which we use to construct $\psi $. Let $f\in\cup
K^s(\varphi )$ be of the form
$f=\frac{1}{m_{2j+1}}(Ey_{k_1}^{\ast }+y_{k_1+1}^{\ast }+
\ldots +y_{k_2}^{\ast }+z_{k_2+1}^{\ast }+\ldots +z_d^{\ast })$,
where $E$ is an interval of integers $\{ p,p+1,\ldots \}$.
Put
$$k^f=\min\{ t\in\{k_1,\ldots ,k_2-2\}:{\rm supp}Ey_t^{\ast }\cap
{\rm supp}\overline{x}_{(t,i)}\neq\emptyset\;{\rm for}\;{\rm some}\;
i\in\{ 1,2,\ldots ,n_t\}\}.$$
Set $$If=\frac{1}{m_{2j+1}}(y_{k^f+2}^{\ast }+\ldots
+y_{k_2}^{\ast }),$$
\noindent while for the other functionals in $\cup K^s(\varphi )$ set $If=0$.

\noindent We set
 $$\varphi_1=\varphi |_{\cup{\rm supp}If}\;\;\;{\rm and }\;\;\; \varphi_2=
\varphi -\varphi_1.$$

Recall that, for $f\in\cup K^s(\varphi )$ which is a member of a family
of type-I or type-II w.r.t. $\overline{x}_{(k,i)},$ we have defined 
$e_f=\min\{ m:m\in {\rm supp}f\cap {\rm supp}
\overline{x}_{(k,i)}\}$. Let 

\centerline {
$P=\cup \{ F\subset \cup K^s(\varphi ):
F $ is a family of type-I or type-II w.r.t. some $\overline{x}_{(k,i)}\} .$}

The functional $\psi $ is  
 supported in the set 
$\{ e_f: f\in P\}.$
We set
$$\psi_1=\psi |_{\{e_f:f\in P\; {\rm and}\; f\;{\rm is}\;{\rm in}\;{\rm the}
\;{\rm analysis}\;{\rm of}\;\varphi_1\}} \;\;\; {\rm and}\;\;\;
\psi_2=\psi -\psi_1.$$
As in the previous section without loss of generality
we assume that ${\rm supp}\varphi\cap {\rm supp}\overline{x}_{(1,1)}\neq
\emptyset $.

\medskip

\noindent {\bf 3.9 Lemma.} {\it {\rm (a)} For every
$f,g\in\cup K^s(\varphi )$ with $f\neq g$ and
$If\neq 0$, $Ig\neq 0$, we have ${\rm supp}If\cap
{\rm supp}Ig=\emptyset $.

{\rm (b)} Let $F=\{ f_1,\ldots ,f_l\}\subset \cup K^s(\varphi)$
 be a family of type--I or
type--II w.r.t.  	$\overline{x}_{(k,i)}.$ 
Suppose that for  some $p\in\{ 1,\ldots ,l\}$, ${\rm supp}f_p\subseteq
{\rm supp}\varphi_1$. Then,
${\rm supp}f_r\subseteq {\rm supp}\varphi_1$ for every
$r\in\{ 1,\ldots ,l\}$.

{\rm (c)} Let $F=\{ f_1,\ldots ,f_l\}\subset \cup K^s(\varphi)$
 be a family of type--I or
type--II w.r.t.  	$\overline{x}_{(k,i)}.$ Suppose that for some
$p=1,\ldots ,l$, ${\rm supp }f_p\not\subseteq {\rm supp}\varphi_1.$
Then $f_p|_{{\rm supp}\varphi_2}(\overline{x}_{(k,i)})=f_p(\overline{x}_{(k,i)}).$

{\rm (d)} Let $F=\{ f_1,\ldots ,f_l\}\subset \cup K^s(\varphi)$
 be a family of type--I or
type--II w.r.t.  	$\overline{x}_{(k,i)}.$ If ${\rm supp}f_p\not\subseteq
{\rm supp}\varphi_1$ for some $p=1,\ldots,l$, then, for all $r=1,\ldots ,l$,
 $f_r|_{{\rm supp}\varphi_2}(\overline{x}_{(k,i)})=f_r(\overline{x}_{(k,i)}).$}

\medskip

\noindent {\it Proof:} (a) Let 
$f=\frac{1}{m_{2j+1}}(Ey^*_{k_1}+\cdots
+y^*_{k_2}+z^*_{k_2+1}+\cdots +z^*_{k_3})$ and
$g=\frac{1}{m_{2j+1}}(Ey^*_{t_1}+\cdots
+y^*_{t_2}+z^*_{t_2+1}+\cdots +z^*_{t_3})$. 
 If ${\rm supp}f\cap {\rm supp}g\neq\emptyset$,
then either ${\rm supp} f\subset {\rm supp} g$ or ${\rm supp}g\subset
{\rm supp}f.$ Suppose that the first is true.
Since 
${\rm supp} y^{*}_l\subseteq [\min {\rm supp}y_l,\max {\rm supp}y_l]$,
it is impossible to have
  ${\rm supp}f\subseteq
{\rm supp}y^{*}_l$ for any $t_1\leq l\leq t_2$. It follows that 
${\rm supp}f \subseteq {\rm supp}z^*_t$ for some $t_2+1\leq t\leq t_3$.
 This implies that
${\rm supp} If \cap {\rm supp} Ig=\emptyset.$ 

(b)
 Let $F=\{ f_1,\ldots ,f_l\}$ be a family of type-I or type-II w.r.t.
$\overline{x}_{(k,i)}$ and suppose that ${\rm supp}f_p\subset {\rm supp}
\varphi_1$
for some $p$. If $\# F=1$ there is nothing to prove. So assume that
$\# F\geq 2$. Let $g_F$ be the functional in $\cup K^s(\varphi)$ which
contains $F$ in its decomposition. Since $f_p\in\cup K^s(\varphi_1)$, we
have that $f_p$ belongs to the analysis of $If$ for some 
$If=\frac{1}{m_{2j+1}}(y^*_{k^f+2}+\cdots y^*_{k_2})$. It follows that
$k^f+2\leq k\leq k_2$ and $f_p$ belongs to the analysis of $y^*_k$.
We have to show that ${\rm supp}g_F\subseteq {\rm supp}y^*_k$ or 
equivalently that $g_F$ does not coincide with $f$.
If $w(g_F)=\frac{1}{2}$
then we get ${\rm supp}g_F\subseteq {\rm supp}y^*_k$, since $w(f)<\frac{1}{2}$.
If $w(g_F)<\frac{1}{2}$ then, since $\# F\geq 2$,
 $F$ is of type-I and again we get 
${\rm supp }g_F\subseteq {\rm supp}y^*_k$, since $\cup_{f\in F}{\rm supp} f$
intersects only  ${\rm supp} \overline{x}_{(k,i)}.$

(c) Suppose that ${\rm supp}f_p\cap {\rm supp}Ig\neq \emptyset$
for some $g=\frac{1}{m_{2j+1}}(Ey^*_{k_1}+\cdots +y^*_{k_2}+z^*_{k_2+1}
+\cdots +z^*_{k_3}) \in \cup K^s(\varphi ).$ Then either ${\rm supp}f_p
\subset {\rm supp} g$ strictly or ${\rm supp} g\subseteq {\rm supp} f_p.$
In the first case we get that ${\rm supp}f_p\subseteq {\rm supp} y^*_l$
for some $k^g+2\leq l\leq k_2$ and so ${\rm supp}f_p\subseteq {\rm supp}
\varphi_1$, a contradiction. In the case ${\rm supp}g\subseteq {\rm supp}
f_p$, since ${\rm supp} g\cap {\rm supp}\overline{x}_{(k^g,q)}\neq \emptyset$
for some $q$,
we get by the definition of families of type I and type II w.r.t.
$\overline{x}_{(k,i)}$ that $k\leq k^g$.
So $Ig=\frac{1}{m_{2j+1}}(
y^*_{k^g+2}+\cdots +y^*_{k_2})$ 
does not intersect $\overline{x}_{(k,i)}.$ It follows that
$(f_p-f_p|_{{\rm supp}Ig})(\overline{x}_{(k,i)})=f_p(\overline{x}_{(k,i)}).$
Since ${\rm supp}\varphi_1=\cup_g {\rm supp}Ig$, we conclude that
$(f_p|_{{\rm supp}\varphi_2})(\overline{x}_{(k,i)})=f_p(\overline{x}_{(k,i)}).$

(d) It follows from (b) and (c).    $\Box$

\bigskip

\noindent {\bf 3.10 Lemma.} {\it For $\varphi_2$ we have}
$$|\varphi_2(\sum_{k=1}^n\varepsilon_kb_k\theta_km_{2j_k}
(\sum_{i=1}^{n_k}b_{(k,i)}x_{(k,i)}))|\leq\psi_2(\sum_{k=1}^n
b_k\theta_km_{2j_k}(\sum_{i=1}^{n_k}b_{(k,i)}u_{(k,i)}))+
\frac{1}{m_{2j+2}}.$$

\medskip

\noindent {\it Proof:} By Lemma 3.9(d) 
we have that  $\varphi_2$ satisfies the assumptions of Remark 2.19(c).
The proof follows from this Remark.
   $\Box $

\medskip

\noindent {\bf 3.11 Lemma. }
$$|\varphi_2(\sum_{k=1}^n\varepsilon_kb_k\theta_km_{2j_k}y_k)|\leq
\frac{257}{m_{2j+1}^2},\leqno ({\rm a})$$
$$|\varphi_1(\sum_{k=1}^n\varepsilon_kb_k\theta_km_{2j_k}y_k)|\leq
\frac{4}{m_{2j+1}^2}.\leqno ({\rm b})$$

\medskip

\noindent {\it Proof:} (a) By Lemma 3.10 it suffices to estimate
$$\psi_2(\sum_{k=1}^nb_k\theta_km_{2j_k}(\sum_{i=1}^{n_k}
b_{(k,i)}u_{(k,i)})).$$
\noindent Recall that $u_{(k,i)}$ is of the form
$u_{(k,i)}=\sum_{m\in A_{(k,i)}}a_{m}e_{m}$, where $a_m >0$
 and $\sum_{m\in A_{(k,i)}}a_m\leq 16$.
Let $\{ K^s(\psi_2)\}$ be the corresponding analysis of $\psi_2$.
For $k=1,2,\ldots ,n$ set
$$D_1^k=\{ 
m\in\cup_{i=1}^{n_k}A_{(k,i)}
:{\rm for}\;{\rm all}\;f\in\cup_sK^s(\psi_2)
\;{\rm such}\;{\rm that}\;m\in {\rm supp}f,\;
w(f)>\frac{1}{m_{2j_k}}\},$$

$$D_2^k=\{ 
m\in\cup_{i=1}^{n_k}A_{(k,i)}
:{\rm there}\;{\rm exists}\;f\in\cup_sK^s(\psi_2)
\;{\rm such}\;{\rm that}\;m\in {\rm supp}f\;{\rm and}\;
w(f)<\frac{1}{m_{2j_k}}\},$$

$$D_3^k=\{ 
m\in\cup_{i=1}^{n_k}A_{(k,i)}
: m\notin D_2^k,\;{\rm there}\;{\rm exists}\;
f\in\cup_sK^s(\psi_2)\; {\rm with}\;m\in{\rm supp}f, w(f)=\frac{1}{m_{2j_k}}$$
$${\rm and}\;{\rm there}\;{\rm exists}\;g\in\cup_sK^s(\psi_2)\; {\rm with}\;
{\rm supp}f\subset {\rm supp}g\; {\rm strictly}\;
{\rm and}\;w(g)\leq\frac{1}{m_{2j+2}}\},$$

$$D_4^k=\{ 
m\in\cup_{i=1}^{n_k}A_{(k,i)}
: m\notin D_2^k,\;{\rm there}\;{\rm exists}\;
f\in\cup_sK^s(\psi_2)\;{\rm with}\;
m\in {\rm supp}f, w(f)=\frac{1}{m_{2j_k}}$$
$${\rm and}\;{\rm for}\;{\rm every}\;g\in\cup_sK^s(\psi_2)\; {\rm with}\;
{\rm supp}f
\subset {\rm supp}g, w(g)\geq\frac{1}{m_{2j+1}}\}.$$

Then, $\cup_{p=1}^4D_p^k=\cup _{i=1}^{n_k}
{\rm supp}u_{(k,i)}\cap {\rm supp}\psi_2$. For
every $k$,
$$\psi_2|_{D_2^k}(b_k\theta_km_{2j_k}(\sum_ib_{(k,i)}u_{(k,i)}))
\leq b_k\theta_km_{2j_k}\frac{16}{m_{2j_k+1}}<\frac{1}{m_{2j_k}},$$
\noindent thus
$$\psi_2|_{\cup_kD_2^k}(\sum_kb_k\theta_km_{2j_k}(\sum_ib_{(k,i)}u_{(k,i)}))
\leq\sum_k\frac{1}{m_{2j_k}}<\frac{1}{m_{2j+2}}. \leqno (1)$$
\noindent Also,
$$\psi_2|_{\cup_kD_3^k}(\sum_kb_k\theta_km_{2j_k}(\sum_ib_{(k,i)}u_{(k,i)}))
\leq\sum_kb_k\theta_k\frac{16}{m_{2j+2}}\leq\frac{64}{m_{2j+2}}. \leqno (2)$$

For $k=1,2,\ldots ,n$, $|\psi_2|_{D_1^k}|^{\ast }_{2j_k-1}\leq 1$
(see Notation after Lemmma 2.2).
So, by Lemma 2.4(b),
$$\psi_2|_{D_1^k}(b_k\theta_km_{2j_k}(\sum_ib_{(k,i)}u_{(k,i)}))
\leq b_k\theta_km_{2j_k}\frac{32}{m_{2j_k}^2}\leq
b_k\frac{128}{m_{2j_k}}.$$
\noindent Hence,
$$\psi_2|_{\cup_kD_1^k}(\sum_kb_k\theta_km_{2j_k}(\sum_ib_{(k,i)}
u_{(k,i)}))\leq\sum_kb_k\frac{128}{m_{2j_k}}<\frac{1}{m_{2j+2}}. \leqno (3)$$
 
 For every $k=1,\ldots ,n$, $i=1,\ldots ,n_k$
and every $m\in {\rm supp}u_{(k,i)}\cap D_4^k,$ there exists a unique 
functional $f^{(k,i,m)}\in \cup _sK^s(\psi_2)$ with $m\in {\rm supp}f$,
$w(f)=\frac{1}{m_{2j_k}}$ and such that, for all $g\in \cup_sK^s(\psi_2)$
with ${\rm supp}f \subset {\rm supp}g$ strictly, $w(g)\geq
\frac{1}{m_{2j+1}}.$ By definition, for $k\neq p$ and $i=1,\ldots ,n_k$,
$m\in {\rm supp}u_{(k,i)}$, we have ${\rm supp}f^{(k,i,m)}\cap
D^p_4=\emptyset.$ Also, if $f^{(k,i,m)}\neq f^{(k,r,n)},$ then 
${\rm supp}f^{(k,i,m)}\cap {\rm supp}f^{(k,r,n)}=\emptyset.$

For each $k=1,\ldots ,n$, let $\{ f^{k,t}\}_{t=1}^{r_k}\subset \cup
K^s(\varphi )
$ be a selection of mutually disjoint such functionals with
$D^k_4=\cup _{t=1}^{r_k} {\rm supp}f^{k,t}.$ For each such functional
$f^{k,t}$, we set $H^k_t={\rm supp}f^{k,t}$ and
$$a_{f^{k,t}}=\sum_{i=1}^{n_k}b_{(k,i)}\sum_{m\in H^k_t} a_m.$$
Then,
$$f^{k,t}(b_k\theta_km_{2j_k}(\sum_ib_{(k,i)}u_{(k,i)}))\leq
b_k\theta_ka_{f^{k,t}}.\leqno (\ast )$$

\noindent {\it Claim:} Let $D_4=\cup_{k=1}^nD_4^k.$ Then
 $\psi_2|_{D_4}(\sum_kb_k\theta_km_{2j_k}
(\sum_ib_{(k,i)}u_{(k,i)}))\leq\frac{256}{m_{2j+1}^2}.$

\medskip

\noindent {\it Proof of the claim:} We shall define a functional
$g\in K^{\prime }$ with $|g|_{2j}^{\ast }\leq 1$
and blocks $u_k$ of the basis so that $\| u_k\|_{\ell_1}\leq 16$,
${\rm supp}u_k\subseteq \cup_i{\rm supp}u_{(k,i)}$
and
$$\psi_2|_{D_4}(\sum_kb_k\theta_km_{2j_k}(\sum_ib_{(k,i)}u_{(k,i)}))
\leq g(2\sum_kb_k\theta_ku_k),$$
\noindent hence by Lemma 2.4(b) we shall have the result.

For $f=\frac{1}{m_q}\sum_{p=1}^df_p\in\cup_sK^s(\psi_2|_{D_4})$
we set
$$J=\{ 1\leq p\leq d: f_p=f^{k,t}\;{\rm for}\;{\rm some}\; k=1\ldots ,n, \;
t=1,\ldots , r_k\},$$
$$T=\{ 1\leq p\leq d:\;{\rm there}\;{\rm exists}\;f^{k,t}\;{\rm with}\;
{\rm supp}f^{k,t} \subset {\rm supp}f_p \; {\rm strictly}\}.$$

For every $f\in\cup_sK^s(\psi_2|_{D_4})$ such that $J\cup T=\emptyset $
we set $g_f=0$, while if $J\cup T\neq\emptyset $ we shall define a
functional $g_f$ with the following properties:

 \noindent Let $D_f=
\cup_{p\in J\cup T}{\rm supp}f_p$ and $u_k=\sum a_{f^{k,t}}e_{f^{k,t}}$,
where $e_{f^{k,t}}=e_{\min H_t^k}$.

\noindent Then,

(a) ${\rm supp}g_f\subseteq {\rm supp}f$.

(b) $g_f\in K^{\prime }$ and $w(g_f)\geq w(f)$,

(c) $f|_{D_f}(\sum_kb_k\theta_km_{2j_k}(\sum_ib_{(k,i)}u_{(k,i)}))
\leq g_f(2\sum_kb_k\theta_ku_k)$.

\noindent Let $s>0$ and suppose that $g_f$ have been defined for all
$f\in\cup_{t=0}^{s-1}K^t(\psi_2|_{D_4})$ and let
$f=\frac{1}{m_q}(f_1+\ldots +f_d)\in K^s(\psi_2|_{D_4})\backslash
K^{s-1}(\psi_2|_{D_4})$ where the family $(f_p)_{p=1}^d$
is ${\cal M}'_q$-admissible if $q>1$, or ${\cal S}'$-allowable
if $q=1$. We consider three cases:

\medskip

\noindent {\tt Case (i)}: $\frac{1}{m_q}=\frac{1}{m_{2j_{k_0}}}$ for some
$k_0$, $1\leq k_0\leq n$. Then $f=f^{k_0,t}$ for some
$t$ and we set $g_f=e^{\ast }_{f^{k_0,t}}$. By ($\ast $) we get
$$f(\sum_kb_k\theta_km_{2j_k}(\sum_ib_{(k,i)}u_{(k,i)}))=
b_{k_0}\theta_{k_0}m_{2j_{k_0}}f(\sum_ib_{(k_0,i)}u_{(k_0,i)})
\leq  b_{k_0}\theta_{k_0}a_{f^{k_0,t}}$$
$$= b_{k_0}\theta_{k_0}a_{f^{k_0,t}}e^{\ast }_{f^{k_0,t}}
(e_{f^{k_0,t}})=g_f(b_{k_0}\theta_{k_0}u_{k_0}).$$

\medskip

\noindent {\tt Case (ii)}: $\frac{1}{m_q}>\frac{1}{m_{2j+1}}$. Then if
$J\cup T\neq\emptyset $, set
$$g_f=\frac{1}{m_q}(\sum_{p\in J}e^{\ast }_{f_p}+\sum_{p\in T}g_{f_p}).$$
\noindent For $p\in J$, $f_p=f^{k_p,t}$ for some $(k_p,t)$ and
by ($\ast$ ),
$$f_p(\sum_kb_k\theta_km_{2j_k}(\sum_ib_{(k,i)}
u_{(k,i)}))\leq b_{k_p}\theta_{k_p}a_{f_p}e^{\ast }_{f_p}(e_{f_p}).$$
\noindent For $p\in T$ we obtain by the inductive hypothesis
$$f_p(\sum_kb_k\theta_km_{2j_k}(\sum_ib_{(k,i)}u_{(k,i)}))\leq
2g_{f_p}(\sum_kb_k\theta_ku_k).$$
\noindent Therefore,
$$f(\sum_kb_k\theta_km_{2j_k}(\sum_ib_{(k,i)}u_{(k,i)}))=
\frac{1}{m_q}\sum_{p\in J\cup T}f_p(\sum_kb_k\theta_km_{2j_k}
(\sum_ib_{(k,i)}u_{(k,i)}))$$
$$\leq g_f(2\sum_kb_k\theta_ku_k).$$
\noindent Since ${\rm supp}g_{f_p}\subseteq {\rm supp}f_p$, $e_{f_p}
\in{\rm supp}f_p$ and $J\cap T=\emptyset ,$ 
 we have that  the family $\{ e^*_{f_p}: p\in J\}\cup
\{ g_{f_p} : p\in T\} $ is ${\cal M}_q^{\prime}$ - admissible if $q>1$,
or ${\cal S}'$ - allowable if $q=1$,
 therefore $g_f\in {\cal A}_q^{\prime }$.

\medskip

\noindent {\tt Case (iii)}: $\frac{1}{m_q}=\frac{1}{m_{2j+1}}$.

Suppose that $f_p\in T$.
Then, by the definition of $f^{k,t}$ and $T$,  $w(f_p)\geq
\frac{1}{m_{2j+1}}.$ On the other hand, recall  (Remark 2.19(a)) that $\psi$ is 
defined through $\varphi$, so that every functional in 
$\cup K^s(\psi )$ has the same weight as the corresponding functional
in $\cup K^s(\varphi )$. So, in this case, by the definition of
$L'_{2j+1}$, we get that $w(f_p)<\frac{1}{m_{2j+1}}$  for every $p$. It
follows that $T=\emptyset.$

Recalling also the definition of $If$ and $\psi_2$, we get that 
in this case $\# J\leq 3.$
Let $J=\{ p_1,p_2,p_3\}$ and $f_{p_{\lambda }}=
f^{k_{\lambda },t_{\lambda }},\lambda =1,2,3$. Set
$g_f=\frac{1}{2}(e^{\ast }_{f_{p_1}}+e^{\ast }_{f_{p_2}}+
e^{\ast }_{f_{p_3}})$. By ($\ast $),
$f_{p_{\lambda }}(\sum_kb_k\theta_km_{2j_k}(\sum_i
b_{(k,i)}u_{(k,i)}))$\\
$\leq b_{k_{\lambda }}
\theta_{k_{\lambda }}a_{f_{p_{\lambda }}},\lambda =1,2,3$. Thus,
$$f|_{D_f}(\sum_kb_k\theta_km_{2j_k}(\sum_ib_{(k,i)}u_{(k,i)}))\leq
\sum_{\lambda =1}^3b_{k_{\lambda }}\theta_{k_{\lambda }}
a_{f_{p_{\lambda }}}$$
$$=\sum_{\lambda =1}^3b_{k_{\lambda }}\theta_{k_{\lambda }}
a_{f_{p_{\lambda }}}e^{\ast }_{f_{p_{\lambda }}}(e_{f_{p_{\lambda}}})=
g_f(2\sum_kb_k\theta_ku_k).$$
\noindent This completes the proof of the Claim. By the Claim and 
relations (1), (2), (3), statement (a) follows.

\medskip

(b) We have from Lemma 3.9(a) that for $f,f^{\prime }\in\cup_sK^s(\varphi )$,
$f\neq f^{\prime }$,
$${\rm supp}If\cap {\rm supp}If^{\prime }=\emptyset .\leqno (\ast \ast )$$
\noindent For $f$ with $If\neq 0$, let $If=\frac{1}{m_{2j+1}}(y^{\ast }_p+
\ldots +y^{\ast }_{p+q})$. Since $\{b_k\}$ is decreasing,
$$|If(\sum_k\varepsilon_kb_k\theta_km_{2j_k}y_k)|\leq
\frac{b_p}{m_{2j+1}}.\leqno (\ast\ast\ast )$$
\noindent Set
$$I_1=\{ If:{\rm there}\;{\rm exists}\;h\in\cup_sK^s(\varphi )\;{\rm with}
\;{\rm supp}If 
\subset {\rm supp}h \; {\rm strictly}\;{\rm and}\;
w(h)\leq\frac{1}{m_{2j+1}}\},$$
$$I_2=\{ If:{\rm for}\;{\rm every }\;h\in\cup_sK^s(\varphi )\;
{\rm with}\;{\rm supp}If\subset
 {\rm supp}h {\rm \; strictly},\;w(h)\geq
\frac{1}{m_{2j}}\}.$$
Set also 
$$A_1=\cup_{If\in I_1}{\rm supp}If \; \; {\rm and} \;\;
A_2=\cup_{If\in I_2}{\rm supp}If $$

\noindent Then, by ($\ast \ast$) and ($\ast\ast\ast $),
$$|\varphi_1|_{A_1}(\sum_k\varepsilon_kb_k\theta_k
m_{2j_k}y_k)|\leq\frac{1}{m_{2j+1}^2}.$$
For $If\in I_2$, we set\\
\centerline {$k(f)=\min\{ l: y^*_l \;{\rm is \; in \; the \;
decomposition \; of }\; If\},$}

\centerline {$T=\{k=1,\ldots ,n: k=k(f) \; {\rm for \; some \; }If\in I_2\}$}

\noindent and, for $k=k(f)\in T$, $l_k=\min ({\rm supp}y_k\cap {\rm supp}If).$

\noindent Using ($\ast \ast $) and ($\ast\ast \ast$) we construct in a similar
way as in part (a) a functional $g\in K^{\prime }$,
$|g|^{\ast }_{2j}\leq 1$ such that
$$|\varphi_1|_{A_2}(\sum_k\varepsilon_k
b_k\theta_km_{2j_k}y_k)|\leq g(\sum_{k\in T}b_ke_{l_k}).$$

\noindent
Then by Lemma 2.4(b) we have the result. This completes the proof of the
Lemma. Proposition 3.8 follows. $\Box$

\medskip

Proposition 3.3 follows from Lemmas 3.6, 3.7 and Proposition 3.8.

\bigskip

\noindent {\bf 3.12 Remark.} {\it The space $X$ is reflexive.}

The proof of this is similar to the proof of Theorem 1.27. We need to prove
that: (a) The basis $(e_n)_n$ is boundedly complete. (b) The basis $(e_n)_n$
is shrinking. The proof of (a) is exactly the same as that of Theorem 1.27(a).
For (b) we also follow the proof of Theorem 1.27(b). We just need to notice
that the norming set $L$ of $X$ satisfies the properties of the set $K$
which are used in that proof.

\bigskip

\bigskip

\medskip
%\newpage

\bigskip

\bigskip

\noindent {\sc S.A. Argyros: Athens University, Department of
Mathematics, Athens 15784, Greece}

\noindent e--mail: {\tt sargyros@atlas.uoa.gr}

\medskip

\noindent {\sc I. Deliyanni: Oklahoma State University, Department
of Mathematics, Stillwater,
 OK 74078, USA}

\noindent e--mail: {\tt irene@math.okstate.edu}

\medskip

\noindent {\sc D.N. Kutzarova: Institute of Mathematics, Bulgarian
Academy of Sciences, 1113 Sofia, Bulgaria}

\noindent e--mail: {\tt denka@math.acad.bg}

\medskip 

\noindent {\sc A. Manoussakis: Athens University, Department of
Mathematics, Athens 15784, Greece}

\noindent e--mail: {\tt amanous@eudoxos.dm.uoa.gr}

\end{document}